\documentclass[12pt]{article}
\usepackage{amsmath, amssymb}
\usepackage[all]{xy}
\title{Universal degeneration of Riemann surfaces \\  
and fibered complex surfaces
}
\date{}
\author{Tadashi Ashikaga and Yukio Matsumoto}

\usepackage{graphicx}
\usepackage{pict2e}
\usepackage{amsthm}

\newtheorem{definition}{Definition}[section]
\newtheorem{lemma}[definition]{Lemma}
\newtheorem{theorem}[definition]{Theorem}
\newtheorem{prop}[definition]{Proposition}
\newtheorem{coro}[definition]{Corollary}
\newtheorem{remark}[definition]{Remark}
\newtheorem{example}[definition]{Example}

\newcommand{\wU}{\widetilde{U}}
\newcommand{\wJ}{\widetilde{J}}
\newcommand{\wM}{\widetilde{M}}

\newcommand{\wbB}{\widetilde{{\bf B}}}
\newcommand{\wbM}{\widetilde{{\bf M}}}
\newcommand{\wbf}{\widetilde{{\bf f}}}
\newcommand{\wbJ}{\widetilde{{\bf J}}}
\newcommand{\wf}{\widetilde{f}}
\newcommand{\wS}{\widetilde{S}}
\newcommand{\wy}{\widetilde{y}}
\newcommand{\wF}{\widetilde{F}}

\newcommand{\wW}{\widetilde{W}}
\newcommand{\wV}{\widetilde{V}}
\newcommand{\wT}{\widetilde{T}}

\newcommand{\ww}{\widetilde{w}}

\newcommand{\wk}{\widetilde{k}}
\newcommand{\wx}{\widetilde{x}}

\newcommand{\wDelta}{\widetilde{\Delta}}
\newcommand{\wvarphi}{\widetilde{\varphi}}

\newcommand{\wpi}{\widetilde{\pi}}
\newcommand{\wmu}{\widetilde{\mu}}
\newcommand{\wg}{\widetilde{g}}
\newcommand{\wh}{\widetilde{h}}
\newcommand{\whf}{\widehat{f}}
\newcommand{\whh}{\widehat{h}}
\newcommand{\whw}{\widehat{w}}
\newcommand{\whS}{\widehat{S}}

\newcommand{\whT}{\widehat{T}}
\newcommand{\whW}{\widehat{W}}
\newcommand{\whQ}{\widehat{Q}}
\newcommand{\whpi}{\widehat{\pi}}
\newcommand{\whPhi}{\widehat{\Phi}}
\newcommand{\be}{{\bf e}}

\newcommand{\bs}{{\bf s}}
\newcommand{\bo}{{\bf o}}

\newcommand{\bB}{{\bf B}}
\newcommand{\bC}{{\bf C}}
\newcommand{\bGamma}{{\bf \Gamma}}
\newcommand{\bZ}{{\bf Z}}
\newcommand{\bI}{{\bf I}}
\newcommand{\bP}{{\bf P}}
\newcommand{\bJ}{{\bf J}}
\newcommand{\bN}{{\bf N}}
\newcommand{\bA}{{\bf A}}

\newcommand{\bQ}{{\bf Q}}

\newcommand{\bK}{{\bf K}}
\newcommand{\bR}{{\bf R}}

\newcommand{\bM}{{\bf M}}
\newcommand{\bE}{{\bf E}}
\newcommand{\bL}{{\bf L}}
\newcommand{\bnabla}{{\bf \nabla}}
\newcommand{\bff}{{\bf f}}

\newcommand{\bmu}{{\bf \mu}}

\newcommand{\cO}{{\cal O}}
\newcommand{\cH}{{\cal H}}
\newcommand{\cC}{{\cal C}}
\newcommand{\cB}{{\cal B}}

\newcommand{\cS}{{\cal S}}
\newcommand{\cU}{{\cal U}}

\newcommand{\cG}{{\cal G}}
\newcommand{\cW}{{\cal W}}

\newcommand{\cA}{{\cal A}}

\newcommand{\cP}{{\cal P}}
\newcommand{\cQ}{{\cal Q}}

\begin{document}

\maketitle


\begin{center}
Dedicated to Professor Norbert A'Campo on his 80th birthday
\end{center}

\bigskip

\begin{footnotesize}

\noindent{\bf Contents}

\noindent \S1 {\bf Introduction}

\S1.1 Precise orbifold structure of  $\overline{M}_g^{orb}$

\S1.2 The results of the present paper

\noindent \S2 {\bf Mapping class groups of stable curves}

\S2.1 Lifting of stable curves and Fenchel-Nielsen twist at infinity

\S2.2 Obstruction to lifting homeomorphisms

\S2.3 Weyl groups as mapping class groups

\noindent \S3 {\bf Kuranshi families of stable curves and log structures}

\S3.1 Deformations of stable curves and Kuranishi families

\S3.2 Log structure and real blow up of Kuranishi families

\noindent \S4 {\bf Construction of the universal degenerating family}

\S4.1 Weyl marking and controlled deformation spaces

\S4.2 Kuranishi families over controlled deformation spaces

\S4.3 Orbifold fiber space over the Deligne-Mumford compactification

\noindent \S5 {\bf Automorphisms of stable curves and cyclic equisymmetric strata  
on $\overline{M}_g^{orb}$}

\S5.1 Automorphisms of Riemann surfaces and equisymmetric strata

\S5.2 Logarithmic quadratic representation of automorphisms

\S5.3 Little Teichm\"{u}ller space 
in an orbifold chart of $\overline{M}_g^{orb}$

\S5.4 Automorphisms and cyclic branched coverings of stable curves 

\S5.5 Equisymmetric strata at the boundary charts of $\overline{M}_g^{orb}$

\S5.6 Harris-Mumford coordinates around equisymmmetric strata

\noindent \S6 {\bf Monodromy and orbifold moduli maps of degenerations of Riemann surfaces}

\S6.1 Pseudo-periodic maps and automorphisms of stable curves

\S6.2 Orbifold structures of degenerations of Riemann surfaces

\S6.3 Local orbifold moduli maps and Kodaira-periodicity

\S6.4 Examples of degenerations and their invariants

\noindent \S7 {\bf Recovery of fibered complex surfaces from the universal degenerating family}

\S7.1 Local recovery of degenerations from the universal family

\S7.2 Global orbifold moduli maps for fibered complex surfaces

\S7.3 Global recovery of basic members of fibered complex surfaces

\noindent{\bf References}

\end{footnotesize}


\section{Introduction}
\label{Introduction}
Kodaira \cite{Kodaira} constructed elliptic surfaces in a canonical way 
for given monodromies and $J$-invariants, which he called the basic members of elliptic surfaces.
Our aim is to extend this result to fibered complex surfaces of genus $g \geq 2$.
\par
In our previous paper \cite{M2020}, we constructed a new orbifold structure 
$\overline{M}_g^{orb}$ over the Deligne--Mumford compactification 
by using a certain bordification of Teichm\"{u}ller space.
In this paper, we construct an orbifold fiber space
$$\overline{\pi}: \overline{Y}_g^{orb} \to \overline{M}_g^{orb}$$
such that any fibered complex surface admitting unstable fibers 
can be pulled back from $\overline{\pi}$ 
by the orbifold moduli map $J^{orb}$ which we will construct.
The fiber space $\overline{\pi}:\overline{Y}_g^{orb} \to \overline{M}_g^{orb}$ will be constructed by patching Kuranishi families of stable curves.
The map $J^{orb}$ has a nature similar to Kodaira's $J$-invariant in the sense that 
it  is delicately influenced by the monodromy and has the property of {\it pseudo-periodicity} 
which we will explain below.

In \S1.1, we will fix the notation from Teichm\"{u}ller theory, and will review the structure 
of $\overline{M}_g^{orb}$. 
Then we will state the results of the present paper in \S1.2.

\subsection{Precise orbifold structure of  $\overline{M}_g^{orb}$}
\label{subsec0101}
Throughout the paper, we will fix a closed oriented topological surface $\Sigma_g$ of genus $g ( >1)$. A Riemann surface $S$ 
is usually considered to be {\it a marked Riemann surface} $(S, w)$, in other words, a Riemann surface for which 
an orientation-preserving homeomophism $w : \Sigma_g \to S$ (called a {\it marking})
is fixed. Two marked surfaces $(S_1, w_1)$ and $(S_2, w_2)$ are {\it equivalent} if there exists a 
holomorphic map $f: S_1 \to S_2$ such that $f\circ w_1 \simeq w_2$, where $\simeq$ means ``is isotopic to'' 
(or equivalently, in the case of  closed surfaces, ``is homotopic to''). The set of  equivalence classes is nothing 
but the {\it Teichm\"uller space} modeled on $\Sigma_g$, and it is denoted by $T_g$. By Teichm\"uller's theorem, 
$T_g$ is a metric space  homeomorphic to $\bR^{6g -6}$ (see \cite{IT}, Chapter 5).\par
Teichm\"uller \cite{Teich44},\cite{Teich40} constructed the space $T_g$ and the universal curve on it. He 
found that $T_g$ is a complex manifold  of dimension $3g-3$. Grothendieck \cite{Grothen} had a deep insight into $T_g$ from functorial viewpoint. 
For a review of Teichm\"uller--Grothendieck theory, see \cite{4math}, \cite{Ac16}. Ahlfors--Bers \cite{Ah}, \cite{AhBe} rediscovered the 
complex structure on $T_g$ from analytic viewpoint.\par 

The {\it mapping class group} of genus $g$, $\Gamma_g$, is defined to be 
$$ \Gamma_g =\{ \varphi: \Sigma_g \to \Sigma_g\; |\; \text{orientation preserving homeomorphisms}\}/\simeq.$$
The group structure of $\Gamma_g$ is given by the usual composition of maps: for $[\varphi], [\psi] \in \Gamma_g$, 
we have $[\varphi][\psi]=[\varphi\circ \psi]$. 
The mapping class group $\Gamma_g$ acts on the Teichm\"uller space $T_g$ by the rule: for 
$[\varphi] \in \Gamma_g$ and $[S, w] \in T_g$, we define 
$$\varphi_* ([S, w]) = [S, w\circ \varphi^{-1}].$$
This action is properly discontinuous, and preserves  the Teichm\"uller metric and the complex structure (see \cite{IT}). 
The quotient space $M_g =T_g /\Gamma_g$ is the {\it moduli space}. It is a normal complex variety (see \cite{Car}). 
The moduli space $M_g$ parametrizes all the isomorphism classes of closed Riemann surfaces of genus $g$, but it is 
not compact. By adding frontier points corresponding to ``stable curves'' (i.e. Riemann surfaces with nodes 
and with finite automorphisms), it can be 
compactified. This is the {\it Deligne--Mumford compactification} $\overline{M}_g$ (DM-compactification for short) 
of the moduli space $M_g$ (see \cite{DM}). 
Many authors tried to reconstruct $\overline{M}_g$ from the viewpoint of Teichm\"uller theory. 
Bers \cite{Bers1973} began this project, and Harvey \cite{H1974} topologically reconstructed  $\overline{M}_g$ by 
using discrete group theory.  Hubbard and Koch's paper \cite{HK} contains a bit of history 
concerning the DM-compactification of the moduli space. For recent developments of \lq\lq compactifications of Teichm\"uller space\rq\rq, 
see Ohshika \cite{Ohshika}, Miyachi \cite{Miyachi} or Masai \cite{Masai}.  \par
According to Kra's overview \cite{Kra2012} in this field, an analytic construction of the DM-compactification
 had to await the work of Hubbard and 
Koch \cite{HK} in the twenty-first century. 
But as we explained in \cite{M2020} (by a few sentences after Theorem 6.5 and  by Remark 6.3), 
Hubbard--Koch's orbifold charts contain certain insufficient points as orbifold charts of $\overline{M}_g$.  
On the other hand, we gave a complete set of orbifold charts to the DM-compactification of 
the moduli space $\overline{M}_g$ (\cite[\S 6]{M2020}). \par
\vskip 3mm
\noindent{\bf Harvey's curve complex}\par
\vskip 3mm
For our construction, {\it Harvey's curve complex} $\cC_g$ plays an important role. 
Given a surface $\Sigma_g$, Harvey \cite{Harvey} introduced an abstract simplicial complex called the 
{\it curve complex} $\cC_g =\cC(\Sigma_g)$.  \par
By definition, a {\it vertex} of $\cC_g$ is an isotopy class of an essential (i.e. not null-homotopic) 
simple closed curve on $\Sigma_g$. A {\it simplex} $\sigma \in \cC_g$ is a collection of disjoint, 
mutually non-isotopic essential simple closed curves on $\Sigma_g$:
$$\sigma=\langle C_1, \ldots, C_k\rangle.$$
The number $k$ of simple closed curves contained in $\sigma$ will be denoted by $|\sigma|$. It is known 
that $|\sigma| \leq 3g-3$. We have ${\rm dim}\; \sigma =|\sigma| -1$. \par 
Let $[S, w]$ be a point of $T_g$, and let $\sigma \in \cC_g$ be a simplex: $\sigma =\langle C_1, C_2, \ldots, C_k\rangle$. 
We represent each simple closed curve $w(C_i)$ on the Riemann surface $S$ by a geodesic, and 
contract each of them to a point. Then we have a stable curve. 
  Thus the topological type of a stable curve  is described by a 
simplex $\sigma$ of the curve complex $\cC_g$, and the topological model of the 
stable curve is denoted by $\Sigma_g(\sigma)$. \par
\vskip 3mm
\noindent{\bf Weyl groups}
\vskip 3mm
Another important object associated with a simplex $\sigma \in \cC_g$ is the {\it Weyl group} $W(\sigma)$,  
which is defined as follows:\par
Let $\Gamma(\sigma)$ be a subgroup of the mapping class group $\Gamma_g$ generated by the 
Dehn twists about the simple closed curves $C_i, i=1, \ldots, k$ on $\Sigma_g$. Let $N\Gamma(\sigma)$ be 
the {\it normalizer} of $\Gamma(\sigma)$ in $\Gamma_g$. Then we can prove that 
{\it a mapping class $\varphi \in \Gamma_g$ belongs to the normalizer $N\Gamma(\sigma)$ if and 
only if $\varphi$ permutes the isotopy classes of the curves $C_i,\; i=1, \ldots, k$ of $\sigma$} (Theorem 4.5 in \cite{M2020}).\par
Now the definition of the {\it Weyl group} $W(\sigma)$ is the following:
$$W(\sigma):= N\Gamma(\sigma)/\Gamma(\sigma).$$
Note that the Weyl group $W(\sigma)$ in our sense is not necessarily a finite group. In \S2, we will prove that
{\it the Weyl group $W(\sigma)$ is the mapping class group of a (topological) surface with nodes $\Sigma_g(\sigma)$}.\par
\vskip 3mm
\noindent {\bf Little Tichm\"uller space $T(\sigma)$}\par
\vskip 3mm
Just as the Teichm\"uller space $T_g$ was constructed by using  marked Riemann surfaces $(S, w)$ 
($w$ being an orientation-preserving homeomorphism $w: \Sigma_g \to S$), the {\it little 
Teichm\"uller space} $T(\sigma)$ is constructed using marked stable curves $(S, w)$ (where $w: \Sigma_g(\sigma) \to 
S$ is an orientation-preservig homeomorphism, $S$ being a stable curve). $T(\sigma)$ is a bounded 
domain of $\bC^{3g-3-|\sigma|}$.\par
\vskip 3mm
\noindent{\bf Controlled deformation space $D_\varepsilon(\sigma)$.}
\vskip 3mm
Let $M$ be a (2-dimensional) {\it Margulis constant}: two distinct simple closed geodesics on any Riemann surface $S$ 
are disjoint if their lengths are smaller than $M$. Of course, any positive number smaller than $M$ has again the same property. 
Thus the Margulis constant is not unique.
Let $\varepsilon$ be a positive number smaller than a Margulis constant $M$. We will fix such an $\varepsilon$ 
throughout the present paper.\par
Let $\sigma$ be a simplex of $\cC_g$. In \S4.1, we will give the definition of the {\it controlled deformation space} $D_\varepsilon(\sigma)$ 
(Def. 4.3). 
$D_\varepsilon(\sigma)$ is a complex manifold of complex dimension $3g-3$, and is homeomorphic to 
an open $(6g-6)$-cell (see Lemma 6.7 in \cite{M2020}).  $D_\varepsilon(\sigma)$ contains $T(\sigma)$. The Weyl 
group $W(\sigma)$ acts on $D_\varepsilon(\sigma)$. This action is properly discontinuous and holomorphic 
(see Lemmas 6.4 and 6.6 of \cite{M2020}). We can prove that the compactified moduli space 
$\overline{M}_g$ is covered by $D_\varepsilon(\sigma)/ W(\sigma)$, where $\sigma$ runs over the simplexes of 
$\cC_g$, and can be empty (see Lemma 6.9 in \cite{M2020}). Thus we have \par
\vskip 3mm
\noindent
{\bf Theorem}(Theorem 6.11 \cite{M2020}) \;In the Deligne--Mumford compactification $\overline{M}_g$, the finite family 
$$\{(D_\varepsilon(\sigma), W(\sigma))\}_{\sigma \in \cC_g/\Gamma_g \cup \emptyset}$$
forms an atlas of orbifold charts of  a complex $(3g-3)$-orbifold.\par
\vskip 3mm 
We will denote the compactified moduli space $\overline{M}_g$ with this  atlas of orbifold charts by 
$\overline{M}_g^{orb}$.\par

\subsection{The results of the present paper}
\label{0102}

The main results consist of the following  three parts (A) $\sim$ (C).

\medskip

{\bf (A)} \  Arbarello-Cornalba \cite{AC2009} reconstructed the universal family of Riemman surfaces 
due to Teichm\"{u}ller \cite{{Teich44}} and Bers \cite{Bers1973} 
over Teichm\"{u}ller space by patching Kuranishi families of 
Riemann surfaces.
Inspired by their work, we construct a family of stable curves over 
the controlled deformation spaces by  patching Kuranishi families of 
stable curves (Th.4.4), which may be considered as a reconstruction of Hubbard--Koch's family \cite[Th.10.1]{HK}.
By patching these families over all the orbifold charts of $\overline{M}_g^{orb}$, 
we have; 
\vskip 3mm
\noindent
{\bf Main Theorem I} (Th.4.8 and Def. 4.9) There exists a (strong) orbifold fibration 
$$\overline{\pi}: \overline{Y}_g^{orb} \to \overline{M}_g^{orb}$$
which is obtained by patching standard Kuranishi families of stable curves.\par
\vskip 3mm
We call $\overline{\pi}$ {\it the universal degenerating family of Riemann surfaces}.
The naming comes from the universality in the sense that any fibered complex surface is 
obtained by pulling back from $\overline{\pi}$.
(We will  explain this point in {\bf (C)}.)
We are also inspired by Arbarello--Cornalba--Griffiths' description  \cite[Chap.XV \S8]{ACG} of the bordification of Teichm\"{u}ller space 
using  real blow-ups and the method of log geometry by Kato-Nakayama \cite{KN} and Usui \cite{U2001}.
We also refer to Hinich-Vaintrob \cite{HV} for a related work. Philosophically we are inspired by the project of 
Catanese \cite{Cat} for studying the relation between the Kuranishi families and the Teichm\"uller theory for 
many kinds of algebraic varieties.\par
\bigskip

{\bf (B)} \ For a mapping class $\varphi \in \Gamma_g$, let us denote by $[\varphi]$ the set of conjugate elements of $\varphi$. 
Harvey \cite{H1971} described the  fixed point locus $T_g^{[\varphi]}$ 
in $T_g$ for a finite (elliptic, or periodic) element $\varphi \in \Gamma_g$ 
in terms of the space of the base Riemann surfaces pointed by the branch points 
for the cyclic covering associated with $\varphi$.
$T_g^{[\varphi]}$  is called the {\it equisymmetric strata} for $\varphi$.
Broughton \cite{Broughton1990} described the similar set $M_g^{[\varphi]}$  in $M_g$.
We extend these results to the boundaries of $T_g$ and $M_g$, i.e. to the 
little Teichm\"{u}ller space $T(\sigma)$ and its image on $\overline{M}_g$.
\vskip 3mm
\noindent
{\bf Main Theorem II} (Th.5.17 and Def. 5.15) 
Let $T^{[\varphi]}(\sigma)$  be the set in $T(\sigma)$ consisting  of the marked stable curves 
with a given numerical  data ${\rm Num}(\varphi)$ of an automorphism.
Then a connected component of $T^{[\varphi]}(\sigma)$ is a complex manifold 
which is isomorphic to a direct product of pointed Teichm\"{u}ller spaces.
\par The image $M^{[\varphi]}(\sigma)$ of $T^{[\varphi]}(\sigma)$ in $\overline{M}_g \setminus M_g$ 
is an irreducible subvariety.
\vskip 3mm
$T^{[\varphi]}(\sigma)$ may be identified with the fixed point locus for a parabolic  (pseudo-periodic) 
element in $\Gamma_g$.
The direct factor of the structure of $T_g^{[\varphi]}$ appears as the space of pointed Riemann surfaces 
which come from each component of the base nodal Riemann surfaces of the cyclic covering associated with $\varphi$ 
(see the discussion in \S5.4 and the precise statement of Th.5.17).
Our discussion is similar to those of the moduli construction of the branched covering of Riemann 
surfaces (cf. \cite{V1994}, \cite{MN} etc.).
We are also inspired by Terasoma's argument \cite{Terasoma1998}.
Note that the equisymmetric theory on $M_g$ via various 
group actions was recently developed by 
Takamura \cite{Tak2019}, Hirakawa-Takamura \cite{HiTa}  et al. 

\bigskip

{\bf (C)} \ Kodaira \cite{Kodaira} constructed  the {\it basic member} of an elliptic surface (without multiple fibers) 
in a canonical way from 
the data of the monodromy and the $J$-invariant.
We will extend this result to any fibered complex surface of genus $g \geq 2$ (cf.~Remark 7.9). 
\par We start from the  local viewpoint.
Let $f: M \to \Delta$ be a degeneration of genus $g \geq 2$ with the degenerate fiber $F=f^{-1}(0)$ 
over a small disk $\Delta$, 
and  $\mu_f$ be the topological monodromy of $f$.
Then $\mu_f$ belongs to the conjugacy class 
$\hat{{\cal P}}_g^{(-)}(\sigma) \subset \hat{\Gamma}_g$ of pseudo-periodic maps 
of negative twists for some  $\sigma \in {\cal C}_g$ (cf. \cite{Nielsen2}, \cite{MM}).
Now as an extended notion of  local J-invariant, we define the {\it local orbifold moduli map} 
$$\{ \wJ: \wDelta \to B \subset  D_{\epsilon}(\sigma), 
\: J:\Delta \to M_{\epsilon}(\sigma)  \subset \overline{M}_g; G=\langle \mu \rangle \simeq \bZ/N\bZ \}$$
as follows. Here $J$ is the usual holomorphically extended moduli map for $f$ (cf. \cite{Im81}), 
$\wDelta \to \Delta$ is the totally ramified covering of disks with  covering degree N which is the 
pseudo-period of the monodromy $\mu_f$, 
$\wf: \wM \to \wDelta$ is the precise stable reduction of $f$ (\cite[\S2]{A2010}) and 
$\wJ$ is the natural map associated 
with the deformation $\wf$
to the Kuranishi space $B$ (as a chart of  $D_{\epsilon}(\sigma)$) 
of the stable curve $\wF=\wf^{-1}(0)$.
Then $\wJ$  can be considered as the lifting map of $J$ by the cyclic group $G$ whose generator 
$\mu$  essentially comes from the automorphism of $\wF$.
In order to describe the map $\wJ$ explicitly, we define the following;
\vskip 3mm
\noindent
{\bf Definition} (Defs. 6.8, 6.11)
A holomorphic function $\varphi(u)$ is called a {\it pseudo-periodic function} with  multiplicity $\gamma \in \bN$ 
and period $L \in \bN$ if the Taylor expansion at the origin is given by 
$\varphi(u)=\sum_{i=0}^{\infty} c_iu^{\gamma+iL}$ ($c_0 \not=0$).
We  call $\gamma/L \in \bQ$ the {\it analytic screw number} of $\varphi(u)$.
\vskip 3mm
\noindent
{\bf Main Theorem III} (Th.6.12) 
Let $(z_1, \cdots, z_k, z_{k+1}, \cdots z_{3g-3})$ be the Harris--Mumford 
coordinates at $\wJ(0)$ of $B \subset D_{\epsilon}(\sigma)$, 
where $z_i$ ($1 \leq i \leq k$) is a smoothing coordinate at a node of $\wF$, 
and $z_j$  ($k+1 \leq j \leq 3g-3$) a variable coordinate of a component of $\wF$ (see Def. 5.21).
Let 
$$\wJ: \wDelta \ni u \longmapsto (z_1, \cdots, z_{3g-3})=(\varphi_1(u), \cdots, \varphi_{3g-3}(u)) \in B \subset D_{\epsilon}(\sigma)$$
be the expression of the map $\wJ$ by these coordinates.
\par
 (i) If $1 \leq i \leq k$, then $\varphi_i(u)$ is a pseudo-periodic function 
whose multiplicity and  period are explicitly written by the numerical data of the topological monodormy $\mu_f$.
In particular, the analytic screw number of $\varphi_i(u)$ coincides with the absolute value of Nielsen's screw number of $\mu_f$ 
for a cut curve in $\sigma$.
\par
(ii) If $k+1 \leq j \leq 3g-3$ and $\varphi_j(u)$ is not identically zero, then $\varphi_j(u)$ is a pseudo-periodic function 
whose period and  fractional part of the analytic screw number are explicitly written by the numerical  data of $\mu_f$.
\vskip 3mm
Note that the fractional part of the analytic screw number of $\varphi_j(u)$ ($k+1 \leq j \leq 3g-3$) 
 is written by a logarithmic quadratic character induced from the total valency of $\mu_f$ by  
 the discussion in \S5.2.
But the integral part of the analytic screw number of $\varphi_j(u)$ is not a monodromy invariant, and  it is 
purely an analytic invariant of the fiber germ (see Remark 6.14).
\par Pseudo-periodicity has already appeared in Kodaira's local description \cite[\S8]{Kodaira}
of J-invariant 
by the action of the element of ${\rm SL}(2, \bZ)$ which is the homological monodromy of a degeneration 
of elliptic curves. 
\vskip 3mm
Now we discuss from the global point of view.
Let ${\bf f}:\bM \rightarrow \bB$ be a global fibered complex surface of genus $g \geq 2$ with 
descriminant locus ${\rm Disc}_{{\bf f}}(\bB)=\{ Q_1, \cdots, Q_s \} \subset \bB$.
We set $\bB^{(0)}=\bB \setminus {\rm Disc}_{{\bf f}}(\bB)$, and let 
$$\mu_{\bff}: \pi_1(\bB^{(0)}, b_0) \longrightarrow \Gamma_g, \:\:\: (b_0 \in \bB^{(0)})$$
be the monodromy representation of ${\bf f}$.
Let
$$\bB^{{\rm orb}}=
\left( \wbB^{(0)}, \pi_1(\bB^{(0)}, b_0) , \varphi_{\wbB^{(0)}}, \bB^{(0)} \right)
\bigcup_{1 \leq i \leq s} 
\left( \wDelta_i, \bZ/N_i\bZ , \varphi_{\wDelta_i}, \Delta_i \right)$$
be the natural orbifold structure over $\bB$, where $ \wbB^{(0)}$ is the universal 
covering of $\bB^{(0)}$
and $\wDelta_i \to  \Delta_i$ the  covering of disks around $Q_i$ whose degree coincides with the  
pseudo-period of the local monodromy at $Q_i$.
Then the orbifold moduli  map 
$$\bJ_{\bff}^{{\rm orb}}:\bB^{{\rm orb}} \longrightarrow \overline{M}_g^{orb}$$
is well-defined by patching the local orbifold moduli maps and the natural map from 
$\wbB^{(0)}$ to $T_g$ (Def. 7.7).\par
\vskip 3mm
\noindent
{\bf Main Theorem IV} (Th.7.8, Def.4.10) 
Any fibered complex surface ${\bf f}:\bM \rightarrow \bB$ can be pulled back in the orbifold sense 
 from the universal degenerating family of Riemann surfaces  
$\overline{\pi}: \overline{Y}_g^{orb} \longrightarrow  \overline{M}_g^{orb}$ 
via the orbifold muduli map $\bJ_{\bff}^{{\rm orb}}$.
\vskip 3mm
The above theorem may be considered as an accomplishment of the results of Imayoshi \cite{Im81} and 
\cite{M2012}.

\section{Mapping class groups of stable curves}
\label{MapClassStable}

In this section, we will define and study the mapping class group of  a stable curve $S$ of genus $g$.
We already described it in  \cite[Cor.4.7]{M2020} as a certain quotient group 
of a subgroup of the mapping class group of a Riemann surface, 
and called it {\it the Weyl group of genus $g$}.
However,  some parts of the argument 
were omitted there.
The aim of this section is to supply these points.

In \S \ref{LiftStableFN}, we define the  lifting map from $S$ to a 
Riemann surface $\Sigma_g$ via real blow-ups with some parameters.
For the study of the boundary of  Teichm\"{u}ller space, 
the contraction maps of simple closed curves on a Riemann surface to nodes   
are fundamental (cf.\cite{Bers1974}, \cite{Ab1980} etc.).
However, this contraction loses the information on twisting parameters of Fenchel--Nielsen coordinates 
along the geodesics isotopic to these curves.
By this lifting, we  
recover the above lost twisting parameters and make it possible 
to discuss the Fenchel--Nielsen deformation (or twist) 
at the boundary.
This argument  is essentially due 
to \cite{KN},\cite{U2001} and \cite{ACG}, and 
will be discussed in 
\S \ref{Subsec0302} in detail.

In \S \ref{Obstruction }, 
we discuss the liftability of a continuous map $f:S_1 \rightarrow S_2$ 
of stable curves to a map between two Riemann surfaces.
If $f$ is a holomorphic or real-analytic map, this is possible.
However, in the case where $f$ is a continuous map, 
there exists an obstruction to the lifting by the existence of 
a map which behaves as an infinite-angled rotation around a node.

In \S \ref{Weylgroup}, we show that this obstruction is 
cancelled in the isotopy class of  oriented homeomorphisms 
by the Alexander trick \cite{Al1923}, and 
we can discuss the mapping class group of a stable curve 
in the language of the mapping class group of a Riemann surface.
Then the 
geometric meaning of the Weyl groups will be clarified in Theorem \ref{Weylgpthm}.


\subsection{Lifting of stable curves and Fenchel--Nielsen twist at infinity}
\label{LiftStableFN}

In this subsection,  we will define and discuss  the  
{\it lifting map}\index{lifting map} 
from a stable curve to a 
Riemann surface via real blow-ups with some parameters.

Let $S$ be a {\it stable curve}\index{stable curve} over $\bC$, i.e. a complete algebraic curve with at most  nodes (ordinary double points) 
as singularities such that a nonsingular rational  component has at least three nodes of $S$.
Assume $S$ has genus $g >1$, and has $k$ nodes $P_1, \cdots, P_k$.
If $S$ is nonsingular (i.e. a compact Riemann surface), then $k=0$.
Let $S=\sum_{i=1}^r S_i$ be the irreducible decomposition 
and let 
$h:\whS=\coprod_{i=1}^r \whS_i \longrightarrow S$
be the normalization 
in the function field.
Topologically, the component $\whS_i$ is obtained 
from a connected component $S_i$ 
of the complement of nodes of $S$ 
by closing the punctures which are the pull back by $h$ of nodes 
to this component, say 
$\bP_i:=\sum_{j=1}^{r_i} P_{i,j}$.
If $\whS_i$ is a projective line, then $r_i \geq 3$.
We call the $r_i$-pointed Riemann surface 
\begin{equation}
\label{Eq0201no1}
(\whS_i, \bP_i)  
\end{equation}
the {\it normalized pointed component} of $S$.
Let
$$
\whpi_i:Bl_{\bP_i}(\whS_i) \longrightarrow \whS_i
$$
be the real blowing up\index{real blowing up} at $P_{i,1}, \cdots, P_{i,r_i}$.
By using a complex local coordinate ($U, z$) at $P_{i,j}= (z=0)$, 
$Bl_{\bP_i}(\whS_i)$ is defined locally by 
$\{ (z, \theta) \in U \times S^1 \:|\ z=|z|e^{\sqrt{-1}\theta} \}$ 
and $\whpi_i$ is induced from the first  projection.
Globally $\whpi_i$ is a real analytic isomorphism on the 
complement of the inverse images of $P_{i,j}$'s, and  each fiber
$
(\whpi_i)^{-1}(P_{i,j}) \:\:\:\:\:\:\: (1 \leq j \leq r_i)
$
is  the circle $S^1$, 
which is called {\it the exceptional circle}\index{exceptional circle} of $\whpi_i$.
The inverse image of the node $P_j$ under the composition map 
$\whh=(\coprod_{i=1}^r \whpi_i) \circ h: 
\coprod_{i=1}^r Bl_{\bP_i}(\whS_i) \longrightarrow S$ 
consists of  two exceptional circles, which we write
\begin{equation}
\label{twocircles}
(\whh)^{-1}(P_j)=C_j^{(1)} \coprod C_j^{(2)} \:\:\:\:\:\:\: (1 \leq j \leq k)
\end{equation}
where the order of two circles $C_j^{(1)}, C_j^{(2)}$ is  arbitrary.
We consider  a map
\begin{equation}
\label{twistparameter}
\Phi=\bigoplus_{j=1}^k \varphi_j \;:\; \bigoplus_{j=1}^k \left( C_j^{(1)} 
\longrightarrow  C_j^{(2)} \right) 
\end{equation}
where each $\varphi_j:C_j^{(1)} \longrightarrow C_j^{(2)}$ is a homeomorphism of circles.
Let $S(\Phi)$ be the topological surface obtained 
from Riemann surfaces with boundaries   
$Bl_{\bP_1}(\whS_1), \cdots, Bl_{\bP_r}(\whS_r)$ 
by pasting the corresponding boundaries (\ref{twocircles})
which are the exceptional circles 
via the identification map $\Phi$ of (\ref{twistparameter}). 
We obtain the natural continuous map
\begin{equation}
\label{liftedRiemann}
\pi(\Phi): S(\Phi) \longrightarrow S
\end{equation}
such that the fiber of the node $(\pi(\Phi))^{-1}(P_j)$ is a circle, say $C_j$, 
which is obtained from  $C_j^{(1)}$ and $C_j^{(2)}$ via the identification 
$\varphi_j$, and the restriction of $\pi(\Phi)$ to the complement 
of the inverse images of nodes 
$S(\Phi) \setminus \coprod_{1 \leq j \leq k} C_j 
\longrightarrow S \setminus \coprod_{1 \leq j \leq k} P_j$ is a homeomorphism. 
Figure I is a simple example in the case where $g=k=r=2$.


\begin{definition}
\label{defliftmap}
We call the map $\pi(\Phi)$ of {\rm (\ref{liftedRiemann})}
the {\rm{lifting}}  of $S$, 
and $S(\Phi)$ the {\rm{lifted Riemann surface}} of $S$ with twist $\Phi$.
We also call $\{ C_1, \cdots, C_k \}$ the {\rm{exceptional circles}} for $\pi(\Phi)$. 
\end{definition}

\begin{example} {\rm (a typical example of lifting)} \:
\label{rotaiontwist}
We identify $S^1$ with $\bR/2\pi \bZ$ via the map
\begin{equation}
\label{expid}
\bR/2\pi \bZ \in x \longmapsto \exp(\sqrt{-1}x) \in S^1.
\end{equation}
Let $\varphi_{\theta}: S^1 \longrightarrow S^1$ be the rotation map of angle 
$\theta$ defined by $x \mapsto x+\theta$.
By the definition of the real blowing up, we have a canonical identification 
$C_j^{(1)}=C_j^{(2)}=S^1$.
Therefore, by a $k$-ple of angles 
$\Theta=(\theta_1, \cdots, \theta_k) \in (S^1)^k$, 
the map
\begin{equation}
\label{thetalift}
\Phi(\Theta)=
\bigoplus_{j=1}^k \varphi_{\theta_j} \;:\; \bigoplus_{j=1}^k \left( C_j^{(1)} 
\longrightarrow  C_j^{(2)} \right)
\end{equation}
is defined.
By using $\Phi(\Theta)$ as the identification map, the lifting
$
\pi(\Phi(\Theta)): S(\Phi(\Theta)) \longrightarrow S
$
is constructed.
We call it
{\em the lifting of $S$ with the rotation angles $\Theta$}.
\end{example}

\begin{figure}[t]
\includegraphics[scale=0.6,bb=45 45 500 450]{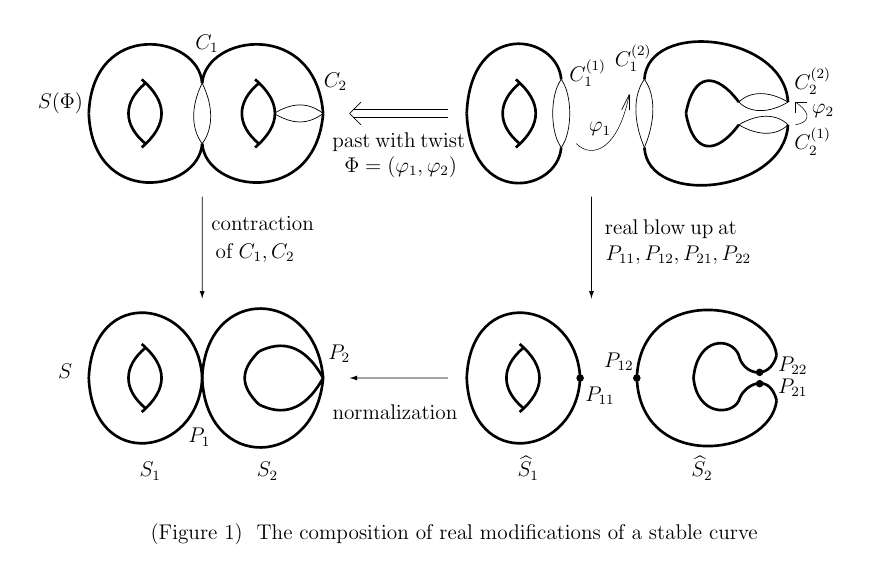}
\end{figure}

We choose the element $\bo=(0, \cdots, 0) \in (S^1)^k$.
By resetting $S(\Phi(\bo))=\Sigma_g$ and $\pi(\Phi(\bo))=\pi_S$, 
we write the lifting of $S$ with the rotation angle $\bo$ as
\begin{equation}
\label{canonicallift}
\pi_S: \Sigma_g \longrightarrow S,
\end{equation}
and call it 
{\em the canonical lifting of $S$}\index{canonical lifting}.
Since $S$ is a stable curve, 
the exceptional circles 
$\{ C_1, \cdots, C_k \}$ on $\Sigma_g$ for $\pi_S$ are not isotopic to 
each other. 
Therefore,  their isotopy classes determine a ($k-1$)-simplex of 
Harvey's curve complex\index{Harvey's curve complex} 
$\cC_g$ (see \S3.1, and \cite{Harvey}) of $\Sigma_g$, which we write
\begin{equation}
\label{excepsimp}
\sigma_{\pi_S}= \langle C_1, \cdots, C_k \rangle \in \cC_g,
\end{equation}
and call it {\it the exceptional simplex for $\pi_S$}.

The lifted Riemann surface $S(\Phi)$ of arbitrary twist $\Phi$ 
is reconstructed from $\Sigma_g$ by changing the pasting parameters along 
exceptional cycles by using (\ref{twistparameter}).
We naturally obtain a homeomorphism 
\begin{equation}
\label{excepsimp}
\tau_{\whPhi}: \Sigma_g \longrightarrow S(\Phi)
\end{equation}
which satisfies $\pi_S=\pi(\Phi) \circ \tau_{\whPhi}$.
If we admit a hyperbolic metric on $\Sigma_g$ 
such that the exceptional circles $C_1, \cdots, C_k$ are geodesics 
with respect to this metric, 
then the map (\ref{excepsimp}) is traditionally called 
the {\it Fenchel--Nielsen deformation}\index{Fenchel--Nielsen deformation}  
or the {\it Fenchel--Nielsen twist}\index{Fenchel--Nielsen twist} along $C_j$'s 
(cf. \cite[p.24]{W2010}, \cite[Chap.8]{IT}).
We also call $\tau_{\whPhi}$ the {\it Fenchel--Nielsen twist at infinity}, 
since $S$ itself sits on the boundary of moduli space 
as we will discuss afterwards.

Note that $\tau_{\whPhi}$ is not uniquely determined from $\Phi$ 
as a homeomorphism.
Since $\Phi$ of (\ref{twistparameter}) is a pasting map along circles 
in the construction of $S(\Phi)$, 
the integral rotation of each circle 
is ignored. 
However, if an identification of $S(\Phi)$ with $\Sigma_g$ is 
somehow fixed, 
the homeomorphism $\tau_{\whPhi}$ may be considered as an  
element in the mapping class group $\Gamma_g$ of $\Sigma_g$ as follows.
By identifying $S(\Phi)$ with $\Sigma_g$, 
the identification map $\varphi_j$ in (\ref{twistparameter}) naturally 
defines a self-homeomorphism 
$\varphi_j: \cA(C_j) \longrightarrow \cA(C_j)$  
of an annular neighborhood $\cA(C_j)$ of $C_j$.
In the mapping class group of the annulus $\cA(C_j)$, $\varphi_j$ 
is isotopic to a rotation map $\varphi_{\theta_j}$ 
with some real-valued angle $\theta_j \in \bR =( -\infty, +\infty)$. 
However , 
if the $\theta_j$'s are integral multiples of $2\pi$, 
the map $\tau_{\whPhi}: \Sigma_g \longrightarrow \Sigma_g=S(\Phi)$ 
is isotopic to a composition of integral Dehn twists along the exceptional circles.

Let $\Gamma(\sigma_{\pi_S})$ be the subgroup of $\Gamma_g$ generated by 
Dehn twists\index{Dehn twist} 
along the exceptional circles, and let 
$[\tau_{\whPhi}]$ be the isotopy  class of $\tau_{\whPhi}$.
By the above argument, 
if the $\theta_j$'s are integral multiples of $2\pi$, 
we have
\begin{equation}
\label{isotopyclassPhi}
[\tau_{\whPhi}] \in \Gamma(\sigma_{\pi_S}).
\end{equation}


\subsection{Obstruction to lifting homeomorphisms}
\label{Obstruction }

Here we discuss the liftability of a homeomorphism of stable curves 
to a homeomorphism of lifted Riemann surfaces with twists.

We start by sketching it locally.
Let  
$f:\Delta \longrightarrow \Delta$ be a map between 
the unit disks with $f(0)=0$ 
such that $f$ is a homeomorphism onto its image.
Let 
$L_{\theta}=\{ z \in \Delta \:|\: \arg z=\theta \}$ 
be a radial segment in $\Delta$ 
with the fixed angle $\theta \in \bR$.
If the limit
\begin{equation}
\label{anglelimit}
f|_{L_{\theta}}(0):=\lim_{z \in L_{\theta}, |z| \mapsto 0} \frac{f(z)}{|f(z)|}
\end{equation}
exists for each $\theta \in \bR$, then $f$ is said to be 
{\it finite-angled}\index{finite-angled}.
Otherwise $f$ is said to be 
{\it infinite-angled}\index{infinite-angled}.
The image of $L_{\theta}$ under an infinite-angled homeomorphism 
goes around the origin infinitely many times.
We first give one of this type of examples, and then state our claims.

\begin{example}
\label{infiniteangled}
Let $[r, \theta]$ be the polar coordinate of $\Delta$; $z=re^{\sqrt{-1}\theta}$.
We define 
$f[r,\theta]=[r, \theta+2\pi/r]$ for $[r, \theta] \not=[0,0]$, and 
$f[0,0]=[0,0]$.
Then $f$ is an infinite-angled homeomorphism of $\Delta$.
\end{example}

\begin{lemma}
\label{realblowuplift}
Let $f:\Delta \longrightarrow \Delta$ be a finite-angled homeomorphism onto 
its image between the unit disks with $f(0)=0$, and let 
$\pi_{0}: Bl_0(\Delta) \longrightarrow \Delta$ be the real blowing up at 
$0 \in \Delta$.
Then there exists uniquely a homeomorphism 
$\whf:Bl_0(\Delta) \longrightarrow Bl_0(\Delta)$ onto its image 
which is a lift of $f$, 
i.e. 
$\pi_{0} \circ \whf=f \circ \pi_{0}$ holds.
\end{lemma}

{\it Proof} \qquad 
We rename the sourse disk by $\Delta_1=\Delta$ and 
the target disk by $\Delta_2=\Delta$, and let $z_i$ be the complex coordinate 
of $\Delta_i$ ($i=1,2$). 
By definition, we have 
$Bl_0(\Delta_i)=
\{ (z_i, \theta_i) \in \Delta_i \times S^1 \:|\: z_i=|z_i|e^{\sqrt{-1}\theta_i} \}$ 
where the coordinate 
$\theta_i$ of $S^1$ is written via the identification (\ref{expid}).
By using (\ref{anglelimit}), we define the map 
$\whf:Bl_0(\Delta_1) \longrightarrow Bl_0(\Delta_2)$ by 
$$
(z_1, \theta_1) \longmapsto (z_2, \theta_2)=(f(z_1), {\rm arg}(f(z_1))
\:\:\:\:\: {\rm if} \:\: z_1 \not=0,
$$
$$
(0, \theta_1) \longmapsto (z_2, \theta_2)=(0,f|_{L_{\theta_1}}(0)).
\:\:\:\:\:\:\:\:\:\:\:\:\:\:\:\:\:\:\:\:\:\:\:\:\:\:\:\:\:\:\:\:\:
$$
Then $\whf$ satisfies the desired property.
\qed
\begin{definition}
\label{finiteangledmapstable}
A homeomorphism $f: S_1 \longrightarrow S_2$ 
of stable curves is said to be {\rm{finite-angled}} 
if the restrictions of $f$ to  small disk neighborhoods of both banks 
(local components) of each node of $S_1$ are finite-angled.
\end{definition}


\begin{prop}
\label{finiteangledstable}
Let  $f: S_1 \longrightarrow S_2$ be a finite-angled homeomorphism 
of stable curves.
We choose a lifting 
$\pi(\Phi_1): S_1(\Phi_1) \longrightarrow S_1$ 
with an arbitrary twist  $\Phi_1$.

Then there exist a unique lifting
$\pi(\Phi_2): S_2(\Phi_2) \longrightarrow S_2$ 
with some twist $\Phi_2$ and a unique homeomorphism 
$\whf:S_1(\Phi_1) \longrightarrow S_2(\Phi_2)$ 
which satisfy the following commutative diagram
$$
\xymatrix{
S_1(\Phi_1) \ar[r]^{\whf} \ar[d]^{\pi(\Phi_1)}
&S_2(\Phi_2) \ar[d]^{\pi(\Phi_2)}\\
S_1 \ar[r]^{f} &S_2.
}
$$

\end{prop}

{\it Proof} \quad 
Let $P_1$ be a node of $S_1$, and $P_2=f(P_1)$ be the image of $P_1$,  
which is a node of $S_2$.
Let $\Delta_i^{(j)}$ ($i=1,2, j=1,2$) be small disk neighborhoods 
in both banks of $P_i$ 
such that the restriction 
$f|_{\Delta_1^{(j)}}: \Delta_1^{(j)} \longrightarrow \Delta_2^{(j)}$ 
is a  homeomorphism onto its image, for $j=1,2$.
Let $\pi_{P_i}^{(j)}: Bl_{P_i}(\Delta_i^{(j)}) \longrightarrow \Delta_i^{(j)}$
be the real blowing up at the origin $P_i$ with exceptional circle 
$C_i^{(j)}=(\pi_{P_i}^{(j)})^{-1}(P_i)$.
Since $f|_{\Delta_1^{(j)}}$ is finite-angled by assumption,
it follows from Lemma \ref{realblowuplift} that 
there exists a unique homeomorphism 
$\whf|_{\Delta_1^{(j)}}:Bl_{P_1}(\Delta_1^{(j)}) \longrightarrow 
Bl_{P_2}(\Delta_2^{(j)})$
which satisfies 
\begin{equation}
\label{floclift}
\pi_{P_2}^{(j)} \circ \whf|_{\Delta_1^{(j)}} 
= f|_{\Delta_1^{(j)}} \circ \pi_{P_1}^{(j)}.
\end{equation}
Let $\varphi_1:C_1^{(1)} \longrightarrow C_1^{(2)}$ 
be the pasting homeomorphism associated with $\Phi_1$.
Then we define the homeomorphism 
$\varphi_2:C_2^{(1)} \longrightarrow C_2^{(2)}$ by
\begin{equation}
\label{newpatch}
\varphi_2=\whf|_{\Delta_1^{(2)}} \circ \varphi_1 
\circ (\whf|_{\Delta_1^{(1)}})^{-1}|_{C_2^{(1)}}.
\end{equation}
By pasting $C_2^{(1)}$ and $C_2^{(2)}$ via $\varphi_2$, 
the desired $S_2(\Phi_2)$ and $\whf$ are locally well-defined by 
(\ref{floclift}) and (\ref{newpatch}).
We do the same construction on disk neighborhoods of all the nodes of 
$S_1$ and $S_2$ and 
paste the resultants trivially 
to their complements  in $S_1$ and $S_2$.
Then we globally obtain the desired $S_2(\Phi_2)$ and $\whf$.
The uniqueness is clear from the construction.
\qed


\subsection{Weyl groups as mapping class groups}
\label{Weylgroup}

Here we characterize the  mapping class group of a stable curve 
as the Weyl group.

First we discuss the orientation of a stable curve $S$.
Since the local equation of a node of $S$ is $xy=0$ in $\bC^2$,
the link of the singularity is a positive Hopf link\index{Hopf link}.
Then $S$ has a natural orientation $\bnabla_S$ so that its  
normalized components  have orientations as Riemann surfaces 
which are compatible with the positivity of Hopf links at nodes.

On the other hand, we consider the 
canonical lifting  
$\pi_S: \Sigma_g \longrightarrow S$
described in (\ref{canonicallift}).
We compare the natural orientaion $\bnabla_{\Sigma_g}$ with $\bnabla_S$.
Then:

\begin{lemma}
\label{orientationdef}
The orientations $\bnabla_{\Sigma_g}$ and $\bnabla_S$ are 
compatible via the map $\pi_S$.
\end{lemma}

{\it Proof} \quad 
Since $\Sigma_g$ and $S$ are realized as fibers 
in an oriented manifold by Theorem \ref{ACGlift} (in the next section), 
the assertion is clear.
\qed

\bigskip

Let $S_i$ ($i=1,2$) be oriented stable curves with orientations $\bnabla_{S_i}$.
A homeomorphism $f:S_1 \longrightarrow S_2$ 
is said to be {\it{oriented}} if the pushed down  $f_*(\bnabla_{S_1})$ defines 
the orientation on $S_2$ which 
is equivalent to  $\bnabla_{S_2}$.
Two oriented homeomorphisms $f_i:S_1 \longrightarrow S_2$ ($i=0,1$) 
are said to be {\it{isotopic}}  to each other if there exists an isotopy 
$\{ f_t: S_1 \longrightarrow S_2 \}_{0 \leq t \leq 1}$ such that 
each $f_t$ is an oriented homeomorphism.
The isotopy class of an oriented homeomorphism 
$f:S_1 \longrightarrow S_2$ is said 
to be the {\it mapping class of $f$}, and is denoted by $[f]$.
For a fixed stable curve $S$, the set of mapping classes  
becomes a group by the composition of maps, which 
is called {\it the mapping class group of $S$}\index{mapping class group},  
and is denoted by $\Sigma(S)$. We have the following:


\begin{prop}
\label{isotopicfiniteangled}
Any oriented homeomorphism $f:S_1 \longrightarrow S_2$ of stable curves is 
isotopic to a finite-angled oriented homeomorphism $f':S_1 \longrightarrow S_2$.
\end{prop}

{\it Proof} \qquad 
Suppose that the restricted homeomorphism  to a disk neighborhood of a bank 
$f|_{\Delta_P^{(1)}}: \Delta_P^{(1)} \longrightarrow 
f(\Delta_P^{(1)}) \subset \Delta_{f(P)}^{(1)}$ 
of some node $P$ of $S_1$ is infinite-angled.
We choose another map 
$h:\Delta_P^{(1)} \longrightarrow \Delta_{f(P)}^{(1)}$ 
such that 

\noindent (i) \ $h(\Delta_P^{(1)})=f|_{\Delta_P^{(1)}}(\Delta_P^{(1)})$ 
as sets and the restricted maps to the boundary coincide: 
$h|_{\partial \Delta_P^{(1)}} \equiv f|_{\partial \Delta_P^{(1)}}$,

\noindent (ii) \ $h$ is a finite-angled homeomorphism onto its image 
so that $h$ and $f|_{\Delta_P^{(1)}}$ define the same orientation.

Note that we have infinitely many choices of $h$.
Since the composite homeomorphisms 
$h^{-1} \circ f|_{\Delta_P^{(1)}}: 
\Delta_P^{(1)} \longrightarrow \Delta_P^{(1)}$ 
coincide with each other at the boundary $\partial \Delta_P^{(1)}$, 
they are isotopic to the identity map ${\rm id}_{\Delta_P^{(1)}}$ 
on the whole disk $\Delta_P^{(1)}$ by the Alexander trick (\cite{Al1923}).
That is to say an $h$ is isotopic to $f|_{\Delta_P^{(1)}}$ 
without moving the points of the boundary.
We define an oriented homeomorphism $f':S_1 \longrightarrow S_2$ by 
$$
f'(x)=\begin{cases}
f(x)&{\rm if}\: x \in S \setminus \Delta_P^{(1)}\\
h(x)&{\rm if}\: x \in \Delta_P^{(1)}.
\end{cases}
$$
Then $f'$ is isotopic to $f$.
If $f'$ is infinite-angled at a certain node, 
then we repeat the same process as above.
After a finite number of steps, we reach a new $f'$ which satisfies the desired property.
\qed

\begin{coro}
\label{liftmappingclass}
Any mapping class $[f]$ of a stable curve $S$ 
has a lifting to a mapping class $[\whf]$ of a Riemann surface $\Sigma_g$ such that 
$(\pi_S)_*([\whf])=[f]$ holds.
\end{coro}

{\it Proof} \quad By identifying 
the Riemann surface $S(\Phi)$ 
with  $\Sigma_g=S(\Phi(\bo))$ as in \S \ref{LiftStableFN}, 
the assertion follows from  
Propositions \ref{isotopicfiniteangled} and \ref{finiteangledstable}.
\qed

\bigskip

The lifted class $[\whf]$ is not uniquely  determined by $[f]$.
This ambiguity comes from the ambiguity 
of the choice of the  map $h$ near a node $P$ as in the proof of 
Proposition \ref{isotopicfiniteangled}.
By the same argument as (\ref{isotopyclassPhi}), 
this ambiguity is cancelled modulo integral Dehn twists along the 
exceptional circles and the lifting $[\whf]$ is uniquely determined as 
an element of $\Gamma_g/\Gamma(\sigma_{\pi_S})$.
On the other hand, since $f$ permutes the nodes $P_1, \cdots, P_k$ of $S$, 
$[\whf]$ permutes the 
ambient isotopy classes $[C_1], \cdots, [C_k]$ of exceptional circles.
By \cite[Th.4.5]{M2020}, 
the subgroup of $\Gamma_g$ consisting the elements which permute 
$[C_1], \cdots, [C_k]$ is nothing but the normalizer\index{normalizer} 
$N\Gamma(\sigma_{\pi_S})$ 
of $\Gamma(\sigma_{\pi_S})$.
 We set 
$$
W(\sigma_{\pi_S})=N\Gamma(\sigma_{\pi_S})/\Gamma(\sigma_{\pi_S})
$$ 
and call it 
{\em the Weyl group of the exceptional simplex $\sigma_{\pi_S}$}\index{Weyl group}.
Note that our $W(\sigma_{\pi_S})$ is not necessarily a finite group.


\begin{theorem} {\rm (\cite[Cor.4.7]{M2020})}
\label{Weylgpthm}
\quad 
The mapping class group $\Gamma(S)$ of a stable curve $S$ 
is isomorphic to the Weyl group
$$
\Gamma(S) \cong W(\sigma_{\pi_S}).
$$

\end{theorem}


\section{Kuranshi families of stable curves and log structures}
\label{Kuranishisection}
Here we review  well-known facts about the  Kuranishi families of 
stable curves which will be used afterwards. 
We mainly refer to Arbarello--Cornalba--Griffiths \cite{ACG}.

In \S \ref{Subsec0301}, we first review the formal properties 
of the standard Kuranishi families 
of stable curves, and secondly the cohomological properties of them.
In particular, we discuss the parameter spaces of 
the smoothing (plumbing) deformations by 
deforming the  nodes to annuli, and those of variable deformations by 
deforming  the complex structures of irreducible components and  shifting the positions of nodes. 

In \S \ref{Subsec0302}, 
we first briefly review the notion of log geometry  
following Kato--Nakayama \cite{KN} and Usui \cite{U2001}, 
and then review in Theorem \ref{ACGlift} that 
the boundary of the log lifting of the Kuranishi family of 
stable curves parametrises the possible Fenchel-Nielsen twist at 
infinity which was discussed in \S \ref{LiftStableFN}.
The method of \cite{ACG} does not use the terminology of log geometry 
and shows this fact direcly by real blowing ups.


\subsection{Deformations of stable curves and Kuranishi families}
\label{Subsec0301}

In this subsection, we review  the standard results 
of the Kuranishi families\index{Kuranishi family} of stable curves.

We recall Harvey's curve complex $\cC_g$ of a Riemann surface $\Sigma_g$ (see \cite{Harvey}). 
By definition, $\cC_g$ is a ($3g-3-1$)-dimensional simplicial complex 
whose ($k-1$)-simplex $\sigma=\langle C_1, \cdots, C_k \rangle$ 
consists of mutually disjoint  and non-homotopic  isotopy classes of simple closed curves $C_j$ 
($1 \leq j \leq k$) on $\Sigma_g$.
An ($\ell-1$)-simplex 
$\rho=\langle C_{i_1}, \cdots, C_{i_{\ell}} \rangle$ 
is a face of $\sigma$, denoted by  $\rho < \sigma$, 
if $\{ C_{i_1}, \cdots, C_{i_{\ell}} \}$ is a subset of 
$\{ C_{1}, \cdots, C_{k} \}$ as isotopy classes.
We denote by
\begin{equation}
\label{cont}
{\rm cont}_{\sigma} : \Sigma_g \longrightarrow \Sigma_g(\sigma)
\end{equation}
the contraction map of $\sigma$, i.e. the continuous map which contracts  
each $C_j$ to a point and is homeomorphic on the complement of the $C_j$'s.  
Here $\Sigma_g(\sigma)$ is a nodal Riemann surface with $k$ nodes.
The map (\ref{cont}) is topologically identified with the map (\ref{canonicallift}).

Let $S$ be a stable curve whose topological type coincides with $\Sigma_g(\sigma)$.
By the {\it Kuranishi family} of $S$, we mean the fibration
\begin{equation}
\label{kuranshistablecurve}
\psi:X \longrightarrow B
\end{equation}
which has the following properties 
(cf. \cite[Chap.XI, \S4, \S6]{ACG}):

\noindent (i) \ $B$ (resp. $X$) is a ($3g-3$)-dimensional 
(resp. ($3g-2$)-dimensional) complex manifold 
with $S=\psi^{-1}(b_0)$  \ ($b_0 \in B$), 
$b_0$ being a fixed point, 

\noindent (ii) \ $\psi$ has the universal property with respect to 
local deformations. In other words, for any deformation 
$\varphi :V \longrightarrow Z$ with $\varphi^{-1}(e_0)=S$, 
the restricted family $\varphi_W: V_W \longrightarrow W$ 
to any sufficiently small connected neighborhood  $W$ ($\subset Z$) of $e_0$ 
has unique holomorphic maps $f:W \longrightarrow B$ and 
$\wf:V_W\longrightarrow X$ with $f(e_0)=b_0$ such that 
$\psi \circ \wf=f \circ \varphi_W$,

\noindent (iii) \ For any $b \in B$, 
the Kodaira-Spencer map\index{Kodaira-Spencer map} 
$$
T_{B,b} \longrightarrow 
{\rm Ext}_{\cO_{X_b}}^1({\rm \Omega}_{X_b}^1, \cO _{X_b})
$$
is an isomorphism, where $T_{B,b}$ is the tangent space of $B$ at $b$ 
and ${\rm \Omega}_{X_b}^1$ is the sheaf of K$\ddot{a}$hler differentials
of $X_b=\psi^{-1}(b)$,

\noindent (iv) \ The discriminant locus\index{discriminant locus}
 $D$ of $\psi$ is a normal 
crossing divisor on $B$ such that the fibers of $\psi$ over 
the irreducible component $D_{i_1 \cdots i_{\ell}}$ of 
$\ell$-codimensional open strata of $D$ are the stable curves whose topological 
types are $\Sigma_g(\rho)$ corresponding to the 
face $\rho=\langle C_{i_1}, \cdots, C_{i_{\ell}} \rangle < \sigma$.

\medskip

Moreover, if $\psi$ satisfies the additional condition (v), 
$\psi$ is said to be a 
{\it  standard Kuranishi family}\index{standard Kuranishi family}
 of $S$:

\noindent (v) \ The action of an analytic automorphism  of $S$ 
extends to a compatible action on $X$ and $B$ (the totality of which will be denoted by
${\rm Aut}(X/B)$).
Conversely, any isomorphism between fibers of $\psi$ is induced 
by an analytic automorphism of $S$ (the totality of which will be denoted by $G$).
That is to say, if there exists an isomorphism of fiberes 
$h_{b_1b_2}:\psi^{-1}(b_1) \rightarrow \psi^{-1}(b_2)$ (for two points $b_1, b_2 \in B$),
then there exists some $h_0 \in G$ and its extension 
$\overline{h}_0 \in {\rm Aut}(X/B)$ ($\overline{h}_0|_{\psi^{-1}(b_0)}=h_0$) 
such that the restriction map $\overline{h}_0|_{\psi^{-1}(b_1)}$ coincides with 
$h_{b_1b_2}$.

\medskip 

\noindent (vi) \ A Kuranishi family and a standard Kuranishi family of $S$ exist 
up to isomorphisms of families and up to shrinking the base $B$ near $b_0$ (\cite{ACG}).

Now we summarize  well-known facts 
about the characterization 
of smoothing (or plumbing) and non-smoothing deformations of $S$ 
for $\psi$
(cf.\cite[Chap.XI]{ACG}, \cite[\S1]{HM}).
Let 
$S=\sum_{i=1}^r S_i$ be a stable curve, and let ($\hat{S}_i, \bP_i$) and $P_j$ ($1 \leq j \leq k$) be as in \S \ref{LiftStableFN}.
Since the tangent space at $b_0$ of $B$ is isomorphic to 
${\rm Ext}^1_{\cO_{S}}(\Omega_{S}^1, \cO_{S})$, 
we can choose a local coordinate neighborhood $B_{loc}$ at $b_0$ of $B$ as an open neighborhood 
of the origin of this vector space ;
\begin{equation}
\label{locnbdbase}
B_{loc} \:\subset \:{\rm Ext^1}_{\cO_{S}}(\Omega_{S}^1, \cO_{S}) \:\cong \:\bC^{3g-3}.
\end{equation}
From the spectral sequence 
${\rm H}^q(S, Ext^p(\Omega_{S}^1, \cO_{S})) \Longrightarrow 
{\rm Ext}^{p+q}(\Omega_{S}^1, \cO_{S})$, we have the exact sequence
\begin{equation}
\label{fundexactseq}
0 \longrightarrow  {\rm H}^1(S, {\it Hom}_{\cO_{S}}(\Omega_{S}^1, \cO_{S})) 
 \longrightarrow {\rm Ext}^1_{\cO_{S}}(\Omega_{S}^1, \cO_{S}) 
 \longrightarrow {\rm H}^0(S, Ext^1(\Omega_{S}^1, \cO_{S})) 
 \longrightarrow 0.
\end{equation}
By using the Grothendieck--Serre duality, the dual vector spaces of each space in (\ref{fundexactseq}) are
\begin{equation}
\label{dualsp}
{\rm H}^1(S, {\it Hom}_{\cO_{S}}(\Omega_{S}^1, \cO_{S}))^* 
\cong \bigoplus_{i=1}^{r} {\rm H}^1(\hat{S}_i, T_{\hat{S}_i}(-\bP_i))^*
\cong \bigoplus_{i=1}^{r} {\rm H}^0 (\hat{S}_i, 2K_{\hat{S}_i}+ \bP_i),
\end{equation}
\vspace{-0.6cm}
$$
{\rm Ext}^1_{\cO_{S}}(\Omega_{S}^1, \cO_{S})^* 
\cong {\rm H}^0(S, \Omega_{S}^1 \otimes \omega_{S}), \:\:\:
{\rm H}^0(S, Ext^1(\Omega_{S}^1, \cO_{S}))^* 
\cong \bigoplus_{j=1}^{k} \tau_{P_j},\:\:\:\:\:\:\:\:\:\:\:\:\:
$$
where $\omega_{S}$ is the dualizing sheaf of $S$, 
$T_{\hat{S}_i}$ and $K_{\hat{S}_i}$ are the tangent sheaf and 
the canonical sheaf of $\hat{S}_i$ respectively, 
and $\tau_{P_j}$ is the 
torsion sheaf supported on $P_j$ which is described as follows (cf. \cite[p.33]{HM}): \ 
If  $S$ is locally defined by $xy=0$ near the node $P_j$, the differentials
$
\omega_1=dx^{\otimes 2}/x,  \:\:
\omega_2=dy^{\otimes 2}/y
$
have the relation $y\omega_1=x\omega_2$. 
Therefore 
\begin{equation}
\label{torsionsheaf}
y\omega_1=\frac{ydx^{\otimes 2}}{x}=\frac{xdy^{\otimes 2}}{y} 
\end{equation}
generates a one-dimensional submodule over $\bC$, which is nothing but the generator of  
$\tau_{P_j}$.
Hence the dual exact sequence of (\ref{fundexactseq}) is 
\begin{equation}
\label{dualfundamentalseq}
0 \longrightarrow  \bigoplus_{j=1}^{k} \tau_{P_j}
 \longrightarrow  {\rm H}^0(S, \Omega_{S}^1 \otimes \omega_{S}) 
 \longrightarrow 
 \bigoplus_{i=1}^{r} {\rm H}^0 (\hat{S}_i, 2K_{\hat{S}_i}+ \bP_i)
 \longrightarrow
0
\end{equation}
Then the deformation-theoretic meaning of (\ref{fundexactseq}) and 
(\ref{dualfundamentalseq}) 
is the following; 

\medskip

\noindent ${\rm (I)}$ \ The space ${\rm H}^1(\hat{S}_i, T_{\hat{S}_i}(-\bP_i))$ 
or  ${\rm H}^0 (\hat{S}_i, 2K_{\hat{S}_i}+ \bP_i)$ parametrizes 
{\it the variable deformations}\index{variable deformation}
 for 
varying the complex structures of $\hat{S}_i$ 
without smoothing the attaching nodes, 
 
\noindent ${\rm (II)}$ \ The torsion sheaf $\tau_{P_j}$ parametrizes 
{\it the smoothing deformations}\index{smoothing deformation}
 of the nodes ${P_j}$ to annuli.
In particular,  $\tau_{P_j}$ is generated by the 
{\it plumbing}\index{plumbing} 
at $P_j$ (Note that here the meaning of  ``plumbing'' is different from that in 
differential topology. For the meaning in the present context, 
see \cite[\S2, \S3]{Kra}, \cite[pp.184--186]{ACG}).
 
\medskip

The above facts ${\rm (I)}$ and ${\rm (II)}$ also follow
from the general theory of deformations of varieties of 
normal crossing (cf.\cite{F1983}).


\subsection{Log structure and real blowing up of Kuranishi families}
\label{Subsec0302}

In this subsection, we apply a part of the standard arguments of 
log geometry\cite{log geometry} 
to the Kuranishi family $\psi:X \rightarrow B$ of $S$.
For the teminology, see 
\cite[\S1]{KN}, \cite[\S1]{U2001}, \cite[Chap.2]{KU} etc.

Since the discriminant locus $D$ is a normal crossing divisor 
on $B$, the sheaf of finitely generated and saturated monoid 
$M=\{ f \in \cO_B\:|\ f\ {\rm is}\ {\rm invertible}\  {\rm outside}\: D \}$ 
embedded in $\cO_B$
defines the log structure\index{log structur} as follows. We set
$$
B^{{\rm log}}=\left\{ 
(b,h) \:|\: b \in B,h \in {\rm Hom}(M_{B,b}^{{\rm gp}}, S^1), 
h(f)=\frac{f(x)}{|f(x)|}, \:\forall f \in \cO_{B,b}^*
\right\}, 
$$
where $M$ is embedded in the abelian group 
$M^{{\rm gp}}=\{ a/b\:|\: a,b \in M \}$ as a sheaf on $B$,  
and $\cO_{B,b}^*$ is the stalk at $b$ of the non-vanishing 
holomorphic functions.
The first projection $(b,h) \mapsto b$ induces a map
\begin{equation}
\label{realblowupkuranshi}
\tau_B: B^{{\rm log}} \longrightarrow B.
\end{equation}
Then $B^{{\rm log}}$ has the  structure of 
a real analytic manifold\index{real analytic manifold} with corners 
such that the map (\ref{realblowupkuranshi}) 
may be considered as 
the real blowing up of $B$ along $D$.
Namely,  let ($U, z_1, \cdots, z_{3g-3}$) be a system of local coordinates 
of $B$ around $b_0$ so that $z_i=0$ ($1 \leq i \leq k$) defines 
locally an irreducible component $D_i$ of $D$.
A component $D_{i_1, \cdots, i_{\ell}}$ of the $\ell$-codimensional open 
strata of  $D$ on $U$ is defined by 
$z_{i_1}= \cdots z_{i_l}=0$,  
$z_j \not=0$ 
($j \in \{ 1, \cdots, 3g-3 \} \setminus \{ i_1, \cdots, i_{\ell} \}$).
Then $B^{{\rm log}}$ is defined near a point of $D_{i_1, \cdots, i_{\ell}} \cap U$ as
$$
\{ (z_1, \cdots, z_{3g-3}, \theta_1, \cdots, \theta_{\ell}) 
\in U \times (S^1)^{\ell} 
\:\:;\:\: z_{i_j}=|z_{i_j}|e^{\sqrt{-1}\theta_j}, \:1 \leq j \leq \ell \}
$$
by identifying $S^1=\bR/2\pi\bZ$.
The exceptional set $E$ of 
$B^{{\rm log}}$ has the stratification so that the component 
$E_{i_1, \cdots, i_{\ell}}=\tau_B^{-1}(D_{i_1, \cdots, i_{\ell}})$ 
is locally written by
\begin{equation}
\label{realblowupkuranshino2}
\{ z_{i_1}= \cdots z_{i_l}=0, \: 
(z_{j_1}, \cdots, z_{j_{3g-3-\ell}}, \theta_1, \cdots, \theta_{\ell})
\in (\bC^*)^{3g-3-\ell} \times (S^1)^{\ell}  \}
\end{equation}
where $j_1, \cdots,j_{3g-3-\ell}  \in  
\{ 1, \cdots, 3g-3 \} \setminus \{ i_1, \cdots, i_{\ell} \}$.
In particular, the restriction 
$E_{i_1, \cdots, i_{\ell}} \longrightarrow D_{i_1, \cdots, i_{\ell}}$ 
of  $\tau_B$ is an $(S^1)^{\ell}$-bundle.

On the other hand, since $\psi^{-1}(D)$ is a normal crossing divisor on $X$, 
we have the construction $\tau_X: X^{{\rm log}} \longrightarrow X$ similar to  (\ref{realblowupkuranshi}).
Then from \cite[Lemma 1.3]{KN}, we have the  log lifting
$\psi^{{\rm log}}: X^{{\rm log}} \rightarrow  B^{{\rm log}}$ of $\psi$ with which the diagram 
\begin{equation}
\label{Diagram0302no1}
\xymatrix{
X^{{\rm log}} \ar[r]^{\tau_X}\ar[d]^{\psi^{{\rm log}}}&X \ar[d]^{\psi}\\
B^{{\rm log}} \ar[r]^{\tau_B} &B
}
\end{equation}
is commutative.
The geometric meaning of the log lifting is clarified by \cite{U2001} 
in a more general setting. 
In the case of the Kuranishi family of a stable curve, it is  
a real analytic family of Riemann surfaces including 
the possible 
Fenchel--Nielsen twists at infinity\index{Fenchel--Nielsen twists at infinity} 
(\cite[p.149--156]{ACG}) :

\begin{theorem} {\rm (\cite{ACG}, \cite{U2001}, \cite{KN})}
\label{ACGlift}
Let $\psi:X \longrightarrow B$ be a Kuranishi family of a 
stable curve $S$. Then 
$\psi^{{\rm log}}: X^{{\rm log}} \rightarrow  B^{{\rm log}}$ 
is a real analytic family of 
Riemann surfaces 
such that

\noindent {\rm (i)} \ The restriction 
$(\tau_B \circ \psi^{{\rm log}})^{-1} (B \setminus D) \longrightarrow 
\tau_B^{-1}(B \setminus D)$  of $\psi^{{\rm log}}$ is isomorphic to 
the restriction 
$\psi^{-1}(B \setminus D) \longrightarrow B \setminus D$ of $\psi$.

\noindent {\rm (ii)} \ Let $Q$ be a point of the fiber $\tau_B^{-1}(P)$ 
of $(S^1)^{\ell}$-bundle 
$E_{i_1, \cdots, i_{\ell}} \longrightarrow D_{i_1, \cdots, i_{\ell}}$ 
$(P \in D_{i_1, \cdots, i_{\ell}})$ 
such that $\Theta=(\theta_1, \cdots, \theta_{\ell})$ is the fiber coordinates 
of $Q$ given in (\ref{realblowupkuranshino2}).
Then the fiber $(\psi^{{\rm log}})^{-1}(Q)$ is isomorphic to the Riemann surface 
which is the lifting of the stable curve $\psi^{-1}(P)$ by rotation angle $\Theta$.
That is to say, in the notation (\ref{thetalift}), we have
\begin{equation}
\label{fiberiso}
(\psi^{{\rm log}})^{-1}(Q) \cong \psi^{-1}(P)(\Phi(\Theta)).
\end{equation}
\end{theorem} 

\begin{remark}
The lifting map of a real analytic map between real analytic manifolds 
via real blowing  ups is also constructed by 
Hubbard--Papadopol--Veselov 
\cite[\S 5]{HPV} and is discussed for other purpose.
\end{remark}


\section{Construction of the universal degenerating family}
\label{ConstUniv}
Starting from Hubbard-Koch's result \cite{HK}, we constructed in 
\cite{M2020} a new orbifold structure 
on the Deligne-Mumford compactification\index{Deligne-Mumford compactification}
, which we will denote here by $\overline{M}_g^{{\rm orb}}$.
In this section, using this structure, we will construct an orbifold fiber space 
$\overline{\pi}^{{\rm orb}}:\overline{Y}_g^{{\rm orb}} \longrightarrow \overline{M}_g^{{\rm orb}}$ 
which is a family of stable curves with some universal properties.

In \S 4.1, we review the orbifold structure of 
$\overline{M}_g^{{\rm orb}}$ with some comments.
The ``biggest chart'' of $\overline{M}_g^{{\rm orb}}$ is the moduli space $M_g$ 
which is the quotient  of Teichmuller space by the mapping class group 
($M_g; T_g, \Gamma_g$).
The ``boundary charts'' are indexed by $\sigma \in \mathcal{C}_g$ (Harvey's curve complex), 
each of which is an  open set $M_\epsilon(\sigma)$ containing the locus $V(\sigma)$ of 
stable curves with  topological type $\Sigma_g(\sigma)$.  
$M_\epsilon(\sigma)$ is the quotient of the controlled deformation space $D_{\epsilon}(\sigma)$ 
by the Weyl group $W(\sigma)$.
Here $D_{\epsilon}(\sigma)$ is the space consisting of $\sigma$-marked 
stable curves whose generalized Fenchel--Nielsen coordinates 
are ``controlled'' so that $W(\sigma)$ acts on $D_{\epsilon}(\sigma)$ 
properly discontinuously.

In \S 4.2, we show that 
$D_{\epsilon}(\sigma)$ is expressed by patching the bases 
of the Kuranishi families of $\sigma$-marked stable curves.
From this, we construct a family of stable curves 
$\pi_{\sigma}: X_{\sigma} \longrightarrow D_{\epsilon}(\sigma)$ 
by patching  the Kuranishi families of stable curves.
This family is a refinement of Hubbard--Koch's 
family \cite[Th.10.1]{HK}, and also is 
an extension of 
Arbarello--Cornalba's theorem  \cite{AC2009} which says that 
the universal curve 
$C_g \longrightarrow T_g$ 
over Teichm\"uller space is expressed by patching 
the Kuranishi families of Riemann surfaces.
As the main tool of the proof of our theorem, we use the log lifting of 
the Kuranishi families discussed in \S \ref{Subsec0302}.

In \S 4.3, by patching the families $\pi_{\sigma}$ naturally, 
we construct the desired orbifold fiber space 
$\overline{\pi}^{{\rm orb}}$. 
We call it 
{\it the universal degenerating  family of Riemann surfaces 
for fibered complex surfaces}, because 
any fibered complex surface admitting unstable fibers 
can be pulled back from $\overline{\pi}^{{\rm orb}}$.
This point will be discussed in \S 7.


\subsection{Weyl marking and controlled deformation spaces}
\label{ContversusKur}

In this subsection,  we introduce the notions of 
the Weyl marking\index{Weyl marking} and 
controlled deformation space from 
\cite[Def.6.2]{M2020}, and review  related results proved in 
\cite{M2020} by adding some comments.

We fix a simplex $\sigma=\langle C_1, \cdots, C_k \rangle$ of 
$\cC_g$, and 
let $S$ be a stable curve whose topological type coincides with 
$\Sigma_g(\sigma)$ in (\ref{cont}).
By a {\it pre-marking} of $S$, we mean the isotopy class   $[w]$  of 
an oriented homeomorphism 
\begin{equation}
\label{fiberiso}
w:  \Sigma_g(\sigma) \longrightarrow S.
\end{equation}

Two such pre-marked stable curves 
($S_1, [w_1]$) and ($S_2, [w_2]$)
are said to be {\it equivalent} (and denoted by ($S_1, [w_1]$) $\sim$ ($S_2, [w_2]$) )
if there exists an analytic isomorphism 
$f:S_1 \longrightarrow S_2$ such that $f \circ w_1$ is isotopic to $w_2$. 
We denote the equivalence class by
$$
[S, w]=(S, [w])/\sim
$$
and call it  a 
{\it  stable curve with Weyl marking}\index{stable curve with Weyl marking} 
(or {\it $\sigma$-marked stable curve}).
This notion of the marking is different from the ones given in 
\cite[p.270]{HK} and \cite[p.490]{ACG}, 
since the action of $\Gamma(\sigma)$ 
is neglected in our case.
Then the 
{\it rigidity}\index{rigidity}, i.e. the automorphism-freeness of a  
stable curve with Weyl marking, is spelled out 
as follows:

\begin{lemma} 
\label{Lem0401no1}
Let $w:  \Sigma_g(\sigma) \rightarrow S$ be a Weyl marking, 
and $f:S \rightarrow S$ be an analytic automorphism such that 
$[S, f \circ w]=[S, w]$.
Then $f$ is the identity map.
\end{lemma}

{\it Proof} \quad If $f$ permutes a proper subset of nodes of $S$ non-trivially, 
then a lifting of $f$ to $\Gamma_g$ would permute a non-trivial subset of the
isotopy classes of $C_1, \cdots, C_k$, 
and $ f \circ w$ cannot be isotopic to $w$.
Therefore, $f$ stabilizes each irreducible component of $S$, and 
stabilizes each normalized component of it as a pointed Riemann surface.
Hence the assertion follows from the usual rigidity of 
the Teichm\"{u}ller marking\index{Teichm\"{u}ller marking}
 (cf. \cite[Prop.~6.8.1]{Hbook}).
\qed

\bigskip

The following is a criterion of the equivalence of the marking:

\begin{lemma} 
\label{Weylmarklemma}
Two oriented homeomorphisms 
$w_i: \Sigma_g(\sigma) \longrightarrow S_i$  $(i=1,2)$ give the same 
$\sigma$-Weyl marking, namely $[S_1, w_1]=[S_2, w_2]$, if and only if 
there exist an oriented homeomorphism 
$\whf: S_1(\Phi_1) \longrightarrow S_2(\Phi_2)$ which is a lift of 
an analytic isomorphism $f:S_1 \longrightarrow S_2$, 
oriented homeomorphisms 
$ \whw_i:  \Sigma_g \longrightarrow S_i(\Phi_i)$ $(i=1,2)$  
and an action $\alpha$ of the Dehn twists $\Sigma_g \longrightarrow \Sigma_g$ 
$(\alpha \in \Gamma(\sigma))$ 
such that the following diagram is homotopically commutative
$$
\xymatrix{
\Sigma_g \ar[r]^{\whw_1\:} \ar[d]^{\alpha}& S_1(\Phi_1) \ar[d]^{\whf}  \ar[r]^{\:\pi(\Phi_1)}&
S_1  \ar[d]^{f}\\
\Sigma_g  \ar[r]^{\whw_2\:} & S_2(\Phi_2)  \ar[r]^{\:\pi(\Phi_2)}&S_2.
}
$$
\end{lemma}

{\it Proof} \quad The assertion is clear from  (\ref{isotopyclassPhi}) 
and Corollary \ref{liftmappingclass}. 
\qed

\bigskip

We denote by 
$T(\sigma)= \bigcup_{w: \sigma\text{-}Weyl\; marking} [S, w]$  the  set of $\sigma$-marked
stable curves.
For a  face $\rho < \sigma$, we denote by 
$T(\rho)= \bigcup_{w:\rho\text{-}Weyl\;marking} [S, w]$  the  set of $\rho$-marked
stable curves similarly.
If $\rho=\emptyset$, then $T(\emptyset)=T_g=\bigcup_{w: \emptyset\text{-}Weyl\; marking} [S, w]$ 
is the Teichm\"{u}ller space, since a $\emptyset$-$Weyl\;marking\; [S, w]$ 
is nothing but a Riemann surface $S$ with a usual Teichm\"{u}ller marking $w$.

We have the following real structure on $T_g\bigcup_{\rho < \sigma}T(\rho)$ 
as a subspace of 
the augmented Teichm\"{u}ller space\index{augmented Teichm\"{u}ller space}
 $\whT_g$ (cf. \cite[Chap.4]{W2010}).
We fix a maximal simplex 
$$\tilde{\sigma}=\langle C_1, \cdots, C_k, C_{k+1}, \cdots, C_{3g-3} \rangle$$
containing $\sigma=\langle C_1, \cdots, C_k\rangle$, and let 
($\ell_1, \cdots, \ell_{3g-3},\tau_1, \cdots, \tau_{3g-3}$) 
be the associated Fenchel-Nielsen coordinates.
Wolpert \cite{W2008} proves that these coordinates, 
except for the first $k$ twist coordinates $\tau_1, \cdots, \tau_{k}$, 
extend to $T_g \cup T(\sigma)$ continuously as 
\begin{equation} 
\label{extendedFN}
((\ell_j, \tau_j), \ell_i):
T_g\bigcup_{\rho < \sigma}T(\rho) \longrightarrow
\prod _{j=k+1}^{3g-3} (\bR_{>0} \times \bR) \times 
\prod _{j=1}^{k} (\bR_{\geq 0}).
\end{equation} 

Moreover, 
if we consider  a point $p=[S, w] \in T(\rho)$ for a face 
$\rho=\langle C_{i_1}, \cdots, C_{i_s} \rangle$, 
then 
the nodes of $S$ correspond to 
 $\ell_{i_j}(p)=0$ for $1 \leq j \leq s$. 

Let $M$ be a universal constant such that two distinct 
simple closed geodesics\index{closed geodesic} 
on any Riemann surface of genus $g$ 
are disjoint if their lengths are smaller than 
$M$ (\cite{Keen1974}, \cite[\S3.3]{Ab1980}), 
which is sometimes called a 2-dimensional 
{\it Margulis constant}\index{Margulis constant}.
Take a number $0 < \epsilon \leq M$ and fix it throughout 
the discussion.
By using the 
extended Fenchel--Nielsen coordinates\index{extended Fenchel--Nielsen coordinates}
 (\ref{extendedFN}), 
we define the subspace $\cU_{\epsilon}(\sigma)$ 
of $T_g\bigcup_{\rho < \sigma}T(\rho)$ by
\begin{equation}
\label{lengthcond}
\cU_{\epsilon}(\sigma)
=\left\{\: 0 \leq \ell_i < \epsilon \:\:(1 \leq i \leq k), \:\:
{\rm max} \{ \ell_1. \cdots, \ell_k \} < \ell_j \:\:\:
(k+1 \leq j \leq 3g-3) \: \right\}.
\end{equation}
Since $T(\sigma)$ is defined by $\ell_j=0$ ($1 \leq j \leq k$), 
we have a natural inclusion $T(\sigma) \subset \cU_{\epsilon}(\sigma)$ 
and $\Gamma(\sigma)$ naturally acts on $\cU_{\epsilon}(\sigma)$.

\begin{definition} {\rm (\cite[Def.6.2]{M2020})} \quad 
\label{ontrolledDS}
The quotient space 
$$
D_{\epsilon}(\sigma)=\cU_{\epsilon}(\sigma)/\Gamma(\sigma)
$$
is called the 
{\rm{controlled deformation space}}\index{controlled deformation space} 
with respect to the simplex $\sigma$.
\end{definition} 

These spaces are considered as refinements of Bers' deformation spaces 
(\cite{Bers1974}), and have the following properties:

\noindent (I) (topological properties) $D_{\epsilon}(\sigma)$ is a 
Hausdorff topological space and 
the Weyl group $W(\sigma)$ acts on  $D_{\epsilon}(\sigma)$ 
as the change of the marking
\begin{equation} 
\label{changeWmark}
[S, w] \longrightarrow [S, w \circ \varphi^{-1}] 
\:\:\:\:\:\: {\rm for} \:\: \varphi \in W(\sigma).
\end{equation} 
This action is properly discontinuous 
(\cite[Lemmata 6.1 $\sim$ 6.4]{M2020}).

\noindent (II) (analytic properties) $D_{\epsilon}(\sigma)$ has a 
complex structure on which $W(\sigma)$ acts holomorphically, and  
the quotient space 
$M_{\epsilon}(\sigma)=D_{\epsilon}(\sigma)/W(\sigma)$ has an analytic open 
embedding into $\overline {M}_g$ (\cite[Lemmata 6.5 $\sim$ 6.8]{M2020}).

\medskip

The property (I) depends on the real analytic augumented  Teichm\"{u}ller 
theory (cf.\cite{Ab1980}) and some 
Weil--Petersson geometry (cf. \cite{W2008}, \cite{Yamada}).
Note that the delicate condition (\ref{lengthcond}) 
for the geodesic lengths $\ell_i$ comes from 
the determination of the region on which $W(\sigma)$ 
acts properly discontinuously (see \cite[Remark 6.3]{M2020}).
The property (II) is essentially comes   
from Hubbard--Koch's theorem  \cite[Th.10.1]{HK}.
In the next subsection, we add a  stronger analytic property.


\subsection{Kuranishi families over controlled deformation spaces}
\label{ProofThmstronganalytic}

In this subsection, we first intoduce  the complex structure on $D_{\epsilon}(\sigma)$ 
by a method which is independent of Hubbard-Koch's theorem 
\cite[Th.10.1]{HK}, and 
construct an analytic family over $D_{\epsilon}(\sigma)$ 
by patching Kuranishi families of stable curves. 
The basic idea of the method here is due to 
Arbarello--Cornalba \cite[\S\S3,4]{AC2009} and 
Arbarello--Cornalba--Griffiths \cite[\S8]{ACG}, 
but our proof of the following theorem 
will be rather direct using the 
topological properties (I) of $D_{\epsilon}(\sigma)$ in \S \ref{ContversusKur}.

\begin{theorem} 
\label{KuranishioverCont}
There exist a $(3g-2)$-dimensional complex manifold $X_{\epsilon}(\sigma)$, 
a complex structure $D_{\epsilon}(\sigma)= \bigcup_{i \in \bI} B_i$
and a holomorphic map 
$$
\pi_{\sigma}: X_{\epsilon}(\sigma) \longrightarrow D_{\epsilon}(\sigma)
$$
such that 
each restricted family
$
\pi_{\sigma}|_{X_i} : X_i=(\pi_{\sigma})^{-1}(B_i) 
\longrightarrow B_i
$
coincides with 
a standard Kuranishi family 
of a stable curve\index{standard Kuranishi family of stable curve} 
$(\pi_{\sigma})^{-1}(b_i)$ for some $b_i \in B_i$.

The Weyl group\index{Weyl group} $W(\sigma)$ acts on
$\pi_{\sigma}$  holomorphically and 
properly discontinuously\index{properly discontinuous} so that 

\begin{picture}(200,71)(-8,-3)
\put(80,47){\makebox(6,6){$X_{\epsilon}(\sigma)$}}
\put(80,5){\makebox(6,6){$D_{\epsilon}(\sigma)$}} 
\put(200,48){\makebox(12,6){$Y_{\epsilon}(\sigma) \cong X_{\epsilon}(\sigma)/W(\sigma)$}} 
\put(204,4){\makebox(12,6){$M_{\epsilon}(\sigma) \cong D_{\epsilon}(\sigma)/W(\sigma)$.}} 
\put(81,43){\vector(0,-1){24}} 
\put(170,43){\vector(0,-1){24}} 
\put(99,51){\vector(1,0){51}}
\put(99,7){\vector(1,0){51}}
\put(68,26){\makebox(6,6){$\pi_{\sigma}$}} 
\put(182,26){\makebox(6,6){$\overline{\pi}_{\sigma}$}} 
\put(121,-3){\makebox(6,6){$\varphi_{\sigma}$}} 
\put(121,39){\makebox(6,6){$\psi_{\sigma}$}} 
\end{picture}

\noindent is a commutative diagram, 
where $\varphi_{\sigma}$ and $\psi_{\sigma}$ are 
the quotient holomorphic maps  to the normal analytic spaces 
$Y_{\epsilon}(\sigma)$ and $M_{\epsilon}(\sigma)$, 
and $\overline{\pi}_{\sigma}$ is the induced holomorphic map.

Moreover $M_{\epsilon}(\sigma)$ has an holomorphic open 
embedding\index{holomorphic open 
embedding} into $\overline {M}_g$.
\end{theorem} 

We need the following:

\begin{lemma} 
\label{imbKuranishitoCDS}
Let $p=[S,w] \in D_{\epsilon}(\sigma) \cap T(\sigma)$ 
be a point, and let 
$\psi: X \longrightarrow B$ be a standard Kuranishi family 
of $\psi^{-1}(b_0)=S$ with a sufficiently small base $B$.
Then there exists a natural topological open embedding
$$
\iota_B: B \longrightarrow D_{\epsilon}(\sigma).
$$
\end{lemma}

{\it Proof} \quad {\it Step 1} \quad 
We set $S=\psi^{-1}(b_0)$.
Let $\whw:\Sigma_g \longrightarrow S(\Phi_0)$ be the lift of $w$ with 
a certain twist $\Phi_0$.
By Theorem \ref{ACGlift}, there exists a point $e_0 \in \tau_B^{-1}(b_0)$ 
such that $S(\Phi_0)$ is canonically identified with the fiber 
\begin{equation} 
\label{firstidfiber}
(\psi^{{\rm log}})^{-1}(e_0)=S(\Phi_0).
\end{equation} 
We consider the stable curve $X_P=\psi^{-1}(P)$ for any $P \in B$, and 
choose a point 
$e_P \in \tau_B^{-1}(P$). Then we  have
\begin{equation} 
\label{secondidfiber}
(\psi^{{\rm log}})^{-1}(e_P)=X_P(\Phi_P) 
\end{equation} 
with a certain  twist $\Phi_P$.
By the connectivity of $B^{{\rm log}}$, 
we can choose a smooth path $L_{e_0e_P}$ 
connecting the points $e_0$ and $e_P$ on $B^{{\rm log}}$.
The family $\psi^{{\rm log}}: X^{{\rm log}}  \longrightarrow B^{{\rm log}}$ is 
a differentiable family of Riemann surfaces over a manifold with corners. 
Usui's theorem \cite{U2001} which is an extension of Ehresmann's theorem (\cite{E1947}) to 
this situation states that there exists a diffeomorphism 
\begin{equation} 
\label{fiberdiffeo}
\varphi_{e_0e_P}: (\psi^{{\rm log}})^{-1}(e_0) 
\longrightarrow (\psi^{{\rm log}})^{-1}(e_P)
\end{equation}
by connecting diffeomorphisms of fibers along $L_{e_0e_P}$.
From (\ref{firstidfiber}), (\ref{secondidfiber}) 
and (\ref{fiberdiffeo}), we have a diffeomorphism 
\begin{equation} 
\label{liftmark}
\whw_P:=\varphi_{e_0e_P} \circ \whw \:: \:
\Sigma_g  \longrightarrow X_P(\Phi_P).
\end{equation}
Note that $\whw_P$ is uniquely determined modulo the action of $\Gamma(\sigma)$,  
independently of 
the choices of $\Phi_0$, $e_P$ and $L_{e_0e_P}$.
Therefore, from Lemma \ref{Weylmarklemma} and (\ref{liftmark}), 
we have a $\rho$-Weyl marking 
\begin{equation} 
\label{aimmark}
w_P \:: \:
\Sigma_g(\rho)  \longrightarrow X_P
\end{equation}
where the face $\rho$ (possibly, $\rho=\emptyset, \:\sigma$) is 
determined by the stratum of $B$ which contains $P$.
Since the extended Fenchel--Nielsen coordinates 
(\ref{extendedFN}) moves continuously (\cite{W2008}), 
the length condition (\ref{lengthcond}) is  satisfied for 
[$X_P, w_P$].
Hence the continuous map 
\begin{equation} 
\label{aimincmap}
\iota_B: B \ni P \longmapsto [X_P, w_P] \in D_{\epsilon}(\sigma)
\end{equation}
is well-defined.

\medskip

{\it Step 2} \quad 
We will prove that the map $\iota_B$ is injective.
We assume $\iota_B(P_1)=\iota_B(P_2)$ for some $P_1, P_2 \in B$.
By definition, the stable curves $ X_{P_i}=\psi^{-1}(P_i)$  ($i=1,2$)  have 
$\rho$-Weyl markings 
$w_i: \Sigma_g(\rho) \longrightarrow X_{P_i}$ for some $\rho$ such that 
there exists an isomorphism 
$g_{P_1P_2}: X_{P_1} \longrightarrow X_{P_2}$ 
satisfying 
$w_2 \simeq g_{P_1P_2} \circ w_1$ (homotopic).

It follows from the property (v) 
of a standard Kuranishi family, stated  
in \S\ref{ACGlift},  
that there exist an automorphism 
$g_{b_0} :X_{b_0} \longrightarrow X_{b_0}$ ($X_{b_0}=S$) 
and a relative automorphism of the family 
\begin{equation}
\label{DiagLem4-5First}
\xymatrix{
X \ar[r]^{\:g_X\:} \ar[d]^{\psi}&X \ar[d]^{\psi}\\
B\ar[r]^{\:g_B\:} &B
}
\end{equation}
such that $g_{b_0}$ and $g_{P_1P_2}$ are induced by the restriction to the fibers 
of the automorphism $g_X$.
The stable curves $X_{P_1}$ and $X_{P_2}$ have the  same topological type 
$\Sigma_g(\rho)$.
Let $B(\rho)$ be the stratum in $B$ such that 
the fibers of $\psi$  over $B(\rho)$ are of type 
$\Sigma_g(\rho)$.

Let ${\bf L}_1$ be a real smooth path on $B$ connecting the points 
$P_1$ and $b_0$ such that 
${\bf L}_1^*={\bf L}_1 \setminus \{ b_0 \}$  is contained in $B(\rho)$.
Then {\it the real proper transform}\index{real proper transform}
 $\widetilde{\bf L}_1$  of 
${\bf L}_1$ on $B^{{\rm log}}$ via the map 
$\tau_B:B^{{\rm log}} \rightarrow B$ is defined as follows.

Assume that ${\rm dim}\:\rho=k_0-1$ and ${\rm dim}\:\sigma=k-1$ ($k_0 \leq k$).
Let ($U;z_1, \cdots, z_{3g-3}$) be local coordinates at $b_0$ of $B$ 
such that 
$B(\rho) \cap U=\{ z_1= \cdots =z_{k_0}=0, 
z_{k_0+1} \not= 0, \cdots, z_{3g-3} \not= 0 \}$
and 
$b_0=\{ (z_1, \cdots, z_k, z_{k+1}, \cdots, z_{3g-3})=(0, \cdots, 0, c_{k+1}, 
\cdots, c_{3g-3}) \}$  
for some $c_{k_1+1},\cdots, c_{3g-3} \in \bC^*$.
Assume that ${\bf L}_1 \cap U$ is defined by
$$
z_{i}=0, \:\:\:
z_{j}=(1-t)r_{j}(t){\rm exp} (\sqrt{-1} \alpha_{j}(t)), \:\:\:
z_{\ell}=c_{\ell} 
$$
for $1 \leq i \leq k_0, \: k_0+1 \leq j  \leq k, \: k+1 \leq \ell \leq 3g-3$, 
where 
$r_{j}(t)$'s  are non-vanishing smooth functions
and $\alpha_{j}(t)$'s  are
smooth functions 
with respect to a variable $t \in [\delta, 1]$ ($0 \leq \exists \delta < 1$).
On the other hand, as stated in \S3.2, 
$\tau_B^{-1}(U)$ is defined by 
$$
\{ (z_1, \cdots, z_{3g-3}, \theta_1, \cdots, \theta_{k}) 
\in U \times (S^1)^{k} 
\:\:;\:\: z_{j}=|z_{j}|e^{\sqrt{-1}\theta_j}, \:1 \leq j \leq k \}
$$
by identifying $S^1=\bR/2\pi\bZ$.
We define the smooth section 
$(\bs_1)|_{{\rm loc}}:{\bf L}_1 \cap U \rightarrow \tau_B^{-1}(U)$ by
$$
\theta_j=0 \:\:\:(1 \leq j \leq k_0), \:\:\:
\theta_j=\alpha_j(t) \:\:\:(k_0+1 \leq j \leq k)
$$
using the fiber coordinates of $U \times (S^1)^{k} \rightarrow U$.
By the similar arguments for other charts on $B$ 
and patching them, 
we obtain the smooth section 
$\bs_1:{\bf L}_1 \rightarrow B^{{\rm log}}$. Then we set  
$\widetilde{\bf L}_1=\bs_1({\bf L}_1)$.
The restriction map $(\tau_B)|_{\widetilde{\bf L}_1}: \widetilde{\bf L}_1 
\rightarrow {\bf L}_1$ is a diffeomorphism by definition.

Now we set 
\begin{equation}
\label{EqLe4-5NewNo2}
{\bf L}_2=g_B({\bf L}_1).
\end{equation}
Then ${\bf L}_2$ is  a real smooth path on $B$ 
connecting  $P_2$ and $b_0$ such that 
${\bf L}_2^*={\bf L}_2 \setminus \{ b_0 \}$ is contained in $B(\rho)$.
Let $\varphi_i:[0,1] \rightarrow {\bf L}_i$ ($i=1,2$) 
be the parametrization such that 
$\varphi_i(0)=P_i, \:\varphi_i(1)=b_0$ and 
$\varphi_2(t)=g_B \circ \varphi_1(t)$ ($ t \in [0,1]$).
By the restriction of $g_X:X \rightarrow X$ 
in (\ref{DiagLem4-5First}) to the fibers, 
we have the family of isomorphisms
\begin{equation}
\label{EqLe4-5NewNo3}
\{ g_t=(g_X)|_{X_{\varphi_1(t)}}: X_{\varphi_1(t)} 
\longrightarrow X_{\varphi_2(t)} \}_{0 \leq t \leq 1}
\end{equation}
such that $g_0=g_{P_1P_2}$ and $g_1=g_{b_0}$.

Let $\widetilde{\bf L}_2=\bs_2({\bf L}_2)$ be the real proper transform of ${\bf L}_2$ 
via $\tau_B$ given by the smooth section 
$\bs_2:{\bf L}_2 \rightarrow B^{{\rm log}}$.
Let $\wvarphi_i=\bs_i \circ \varphi_i:[0,1] \stackrel{\varphi_i}{\longrightarrow} 
\bL_i
\stackrel{\bs_i}{\longrightarrow} \widetilde{\bL}_i$ 
be the parametrization.
By the definition of the real proper transform, or 
the  argument in Proposition  \ref{finiteangledstable} and Theorem \ref{ACGlift},
there exists the homeomorphism 
$\wg_t: X_{\wvarphi_1(t)}^{{\rm log}}=(\psi^{{\rm log}})^{-1}(\wvarphi_1(t)) 
\longrightarrow 
X_{\wvarphi_2(t)}^{{\rm log}}=(\psi^{{\rm log}})^{-1}(\wvarphi_2(t)) $
with  the  commutative diagram
\begin{equation}
\label{DiagLem4-5Second}
\xymatrix{
X_{\wvarphi_1(t)}^{{\rm log}} \ar[r]^{\:\wg_t\:} \ar[d]^{\tau_X}&X_{\wvarphi_2(t)}^{{\rm log}} \ar[d]^{\tau_X}\\
X_{\varphi_1(t)} \ar[r]^{\:g_t\:} &X_{\varphi_2(t)}
}
\end{equation}
where we use the same symbol $\tau_X$ as the restriction to the fibers of 
$\tau_X:X^{{\rm log}} \rightarrow X$.

First we consider the case $t=0$ in (\ref{DiagLem4-5Second}).
By the same reasoning as above, the $\rho$-marking 
$w_i: \Sigma_g(\rho) \rightarrow X_{P_i}$ and the homotopy relation 
$w_2 \simeq g_0 \circ w_1$ are lifted to 
the homeomorphism 
$\ww_i:\Sigma_g \rightarrow X_{\wvarphi_i(0)}^{{\rm log}}$  and 
the homotopy relation $\ww_2 \simeq \wg_0 \circ \ww_1$.

Secondly we consider (\ref{DiagLem4-5Second}) for any $t \in [0,1]$.
By applying \cite{E1947} to the family $\psi^{{\rm log}}$, 
we have the diffeomorphism  $f_{\wvarphi_i(0)\wvarphi_i(t)}:
X_{\wvarphi_i(0)}^{{\rm log}} \rightarrow 
X_{\wvarphi_i(t)}^{{\rm log}}$ 
along the path $\widetilde{\bL}_i$ 
connecting $\wvarphi_i(0)$ and $\wvarphi_i(t)$.
Let $\ww_i(t)=f_{\wvarphi_i(0)\wvarphi_i(t)} \circ \ww_i:
\Sigma_g \rightarrow X_{\wvarphi_i(t)}^{{\rm log}}$ be the composite 
homeomorphism.
By the commutativity of (\ref{DiagLem4-5Second}), we obtain the family 
of the homotopy relations 
\begin{equation}
\label{EqLe4-5NewNo4}
\{ \ww_2(t) \simeq \wg_t \circ \ww_1(t) \}_{0 \leq t \leq 1}.
\end{equation}

Lastly we consider the case $t=1$.
The homeomorphism 
$\ww_i(1):\Sigma_g \rightarrow X_{\wvarphi_i(1)}^{{\rm log}}$ 
descends to the homeomorphism 
$w_i(1):\Sigma_g(\sigma) \rightarrow X_{b_0}$ via the contraction maps 
${\rm cont}_{\sigma}:\Sigma_g \rightarrow \Sigma_g(\sigma)$ 
and 
$(\tau_{X})|_{X_{\wvarphi_i(1)}^{{\rm log}}}:X_{\wvarphi_i(1)}^{{\rm log}} 
\rightarrow X_{b_0}$. 
Then the relation  (\ref{EqLe4-5NewNo4}) for $t=1$ 
descends to 
\begin{equation}
\label{EqLe4-5NewNo5}
w_2(1) \simeq g_{b_0} \circ w_1(1) \:\:\:\: ({\rm homotopic}).
\end{equation}
From  Lemma \ref{Lem0401no1} and (\ref{EqLe4-5NewNo5}), the automorphism 
$g_{b_0}$ is the identity map of $S$.
Hence $P_1=P_2$, i.e. the injectivity of the map $\iota_B$ is 
proved.

\medskip

{\it Step 3} \quad 
By shrinking $B' \subset B$ preserving 
the properties (i) $\sim$(v) in \S\ref{ACGlift},
we may assume that the map 
$\iota_{\overline{B'}}: \overline{B'} \longrightarrow D_{\epsilon}(\sigma)$ 
for the closure $\overline{B'}$ in $B$ 
is injective.
Since $\overline{B'}$ is compact and $D_{\epsilon}(\sigma)$ is Hausdorff, 
$\iota_{\overline{B'}}$ is homeomorphic onto its image.
Therefore the assertion of Lemma \ref{imbKuranishitoCDS} 
follows for  $B'$.
\qed

\bigskip

{\it Proof of Theorem \ref{KuranishioverCont}} :

{\it Step 1} \quad 
Let $p=[S, w] \in D_{\epsilon}(\sigma)$  
 be a point, and let 
$\psi: X \longrightarrow B$ be a standard Kuranishi family of $S$.
By the same argument as in Lemma \ref{imbKuranishitoCDS}, 
there exist an open neighborhood $U_p$ of $p$ in  $D_{\epsilon}(\sigma)$ 
and a homeomorphism $\iota_B: B \longrightarrow U_p$. 
Let
\begin{equation} 
\label{kuranishibasepatch}
D_{\epsilon}(\sigma)=\bigcup_{p \in D_{\epsilon}(\sigma)} U_p
\end{equation}
be the open covering by $U_p$'s of these types. 
By identifying $U_p$ with $B$, the complex structures $\{ B \}$ induce a 
complex structure on $D_{\epsilon}(\sigma)$.

In fact, the coordinate transformations are biholomorphic as follows.
Let $U_q=B'$ be another base of a standard Kuranishi family of 
$S'$ for $q=[S', w'] \in D_{\epsilon}(\sigma)$ such that $U_ p \cap U_q \not= \emptyset$.
Since the Kodaira-Spencer map at any point $r \in U_ p \cap U_q=B \cap B'$ 
is an isomorphism by the property (iii) of the Kuranishi family, stated in \S \ref{ACGlift}, 
$B \cap B'$ is a base of a Kuranishi family of the stable curve 
$\psi^{-1}(r)$ (cf. \cite[Chap.XI, Cor.(4.6)]{ACG}).
By the uniqueness of the Kuranishi family modulo isomorphism, 
the coordinate transformations are nothing but these isomorphisms.
Note that $ D_{\epsilon}(\sigma)$ is a Hausdorff space.
Therefore, (\ref{kuranishibasepatch}) defines the desired complex structure on 
$ D_{\epsilon}(\sigma)$.

These standard Kuranishi families 
$\{ \psi_p:V_p=X^{(p)} \longrightarrow U_p=B^{(p)} \}_{p \in  D_{\epsilon}(\sigma)}$
are patched together, and define a complex manifold 
$X_{\epsilon}(\sigma)$ and a holomorphic map
$$
\pi_{\sigma}= \cup  \psi_p \: : \:
X_{\epsilon}(\sigma)=\bigcup V_p 
 \longrightarrow 
\bigcup U_p 
= D_{\epsilon}(\sigma).
$$
These are also well-patched by the uniqueness of  Kuranishi families 
modulo isomorphisms. 
Therefore, we have constructed the desired family.

Since the action of $ \varphi \in W(\sigma)$ on $ D_{\epsilon}(\sigma)$ is 
given by (\ref{changeWmark}), 
$ \varphi$ sends the open neighborhood $U_p=B^{(p)}$ of $p=[S,w]$ which  is the base 
of the standard Kuranishi family of $S$ to 
an open neighborhood of $U_{\varphi(p)}=B^{(\varphi(p))}$ 
of $\varphi(p)=[S,w \circ \varphi^{-1}]$.
This is nothing but the isomorphism of 
the bases of Kuranishi families, and is holomorphic.
Hence $ \varphi$ acts on $D_{\epsilon}(\sigma)$ holomorphically.

\medskip

{\it Step 2} \quad 
It suffices to prove that $W(\sigma)$ acts on $X_{\epsilon}(\sigma)$ 
properly discontinuously, i.e. 

\noindent (i) let $x', x'' \in X_{\epsilon}(\sigma)$ be points such that 
$g(x') \not= x''$ for an element $g \in W(\sigma)$.
Then there exist open neighborhoods $U'$ and $U''$ of $x'$ and $x''$ respectively 
such that $g(U') \cap U''=\emptyset$,

\noindent (ii) the isotropy sub-group $G_x$ ($ \subset W(\sigma))$ of a point  
$x \in X_{\epsilon}(\sigma)$ is finite, 

\noindent (iii) there exists a $G_x$-invariant neighborhood $U$ of $x$ 
such that, if $g(U) \cap U \neq \emptyset$ for some $g \in W(\sigma)$, 
then $g \in G_x$.

Since $W(\sigma)$ acts on $D_{\epsilon}(\sigma)$ properly discontinuously 
by \cite[Lemma 6.4]{M2020}, 
the assertion (i) is clear in the case where $x'$ and $x''$ belong to 
distinct fibers of $\pi_\sigma$.
Assume that $x'$ and $x''$ belong to the same fiber of $\pi_\sigma$, say 
$x'=[[S], w']$ and $x''=[[S], w'']$ with an isomorphism class 
$[S] \in \overline {M}_g$ and distinct markings.
Then, after forgetting the markings, $W(\sigma)$ acts on $S$ as the analytic automorphism group of $S$ which is finite.
Therefore the assertion (i) holds.

The assertion (ii) is clear, since  $G_x$ is a subgroup of 
the isotropy sub-group of  $\pi_{\sigma}(x)$ by the action of $W(\sigma)$ 
on $D_{\epsilon}(\sigma)$, which is finite 
by \cite[Lemma 6.4]{M2020}.

By the first half of this theorem,  
there exists an open neighborhood $V$ of $\pi_{\sigma}(x)$ in  
$D_{\epsilon}(\sigma)$ such that 
$\pi_{\sigma}^{-1}(V) \longrightarrow V$ is a standard Kuranishi family 
of the stable curve $\pi_{\sigma}^{-1}(\pi_{\sigma}(x))$. 
Therefore the assertion (iii) follows from 
the property (v) stated in \S \ref{ACGlift}.
Hence $W(\sigma)$ acts on $X_{\epsilon}(\sigma)$ 
properly discontinuously, and 
the quotient space $X_{\epsilon}(\sigma)/W(\sigma)$ is a normal 
analytic space by Cartan's theorem (\cite{Cartan}).

The holomorphic embedding of $D_{\epsilon}(\sigma)/W(\sigma)$ into 
$\overline {M}_g$ follows from the above holomorphy  and the 
property (II)  stated in  \S \ref{ContversusKur}. 
\qed


\subsection{Orbifold fiber space over the Deligne-Mumford compactification}
\label{ConstunivdegRiemann}

In this subsection, as a globalization of Theorem \ref{KuranishioverCont}, 
we construct an orbifold fiber space (or orbifold fibration for short) 
$\overline{\pi}_g:
\overline{X}_g^{orb} \longrightarrow \overline{M}_g^{orb}$ 
so that the Bers fiber space (\cite{Bers1973}) over Teichm\"{u}ller space is contained as an open chart 
of $\overline{\pi}_g$.

By a {\it complex orbifold}\index{complex orbifold} $M$, 
we mean here that $M$ is a normal complex space such that 
$M$ is covered by (an atlas of) orbifold charts 
$\{ ( \wU_i, G_i, \varphi_i, U_i) \}_{i \in \bI}$ which satisfy the following conditions:

\noindent {\rm (i)} \ Each chart consists of a complex manifold $\wU_i$, 
a (not necessary finite) group $G_i$ acting on $\wU_i$ holomorphically 
and properly discontinuously 
(admitting non-effective action), 
an open set $U_i$ of $M$ and a folding map 
$\varphi_i: \wU_i \longrightarrow U_i$ 
which induces a natural analytic isomorphism $\wU_i/G_i \longrightarrow U_i$.

\noindent {\rm (ii)}  (compatibility condition) \ 
For $x \in \wU_i $ and $y \in \wU_j$ 
such that $\varphi_i(x)=\varphi_j(y) \in U_i \cap U_j$, 
there exists a biholomorphic map $\varphi_{x,y}:\wU_x^\prime \longrightarrow \wU_y^\prime$ 
from an open neighborhood of $x$ in $\wU_i$ to an open 
neighborhood of $y$ in $\wU_j$ such that 
$\varphi_{x,y}(x)=y$ and $\varphi_j \circ \varphi_{x,y}=\varphi_i$.

For simplicity, we sometimes write this orbifold structure by 
$M= \{ ( U_i; \wU_i, G_i) \}_{i \in \bI}$ or $M=\bigcup_{i \in \bI} \wU_i/G_i$ 
if there is no fear of confusion.

Let $M'$ be another complex orbifold with its orbifold charts 
$\{ ( \wV_k , H_k,  \phi_k, V_k) \}_{k \in \bK}$.
By an {\it orbifold map} $h:M' \rightarrow M$, we mean that 
$h$ is a holomorphic map of  normal complex spaces which satisfies the following conditions:

\noindent {\rm (i)} \
For any points $p \in M'$ and  $h(p) \in M$, let 
$( \wV_k , H_k,  \phi_k, V_k)$ be a small orbifold chart of $M'$ with $p \in V_k$ and 
$( \wU_i, G_i, \varphi_i, U_i)$ be that of $M$ with $h(p) \in U_i$.
Then there exists a holomorphic map $\wh_{ki}:\wV_k \rightarrow  \wU_i$ 
such that 
$\varphi_i \circ \wh_{ki}=h|_{V_k} \circ \phi_k$.

\noindent {\rm (ii)}  (compatibility condition) \ 
Assume that $x \in \wV_k, y \in \wV_{\ell}, 
\wh_{ki}(x) \in \wU_i, \wh_{\ell j}(y) \in \wU_j$
such that $\phi_k(x)=\phi_{\ell}(y)$ and 
$\varphi_i(\wh_{ki}(x))=\varphi_j(\wh_{\ell j}(y))$.
Then there exist open neighborhoods 
$\wV_x, \wV_y, \wU_{\wh_{ki}(x)}$,  $\wU_{\wh_{\ell j}(y)}$ of 
$x, y, \wh_{ki}(x), \wh_{\ell j}(y)$ in  
$\wV_k, \wV_{\ell}, \wU_i, \wU_j$ respectively 
and biholomorphic maps 
$\phi_{x,y}:\wV_x \rightarrow \wV_y$ and 
$\varphi_{h(x),h(y)}:\wU_{\wh_{ki}(x)} \rightarrow \wU_{\wh_{\ell j}(y)}$ 
such that $\varphi_{h(x),h(y)} \circ \wh_{ki}|_{\wV_x}
=\wh_{\ell, j}|_{\wV_y} \circ \phi_{x,y}$.

Moreover, we define the following:

\begin{definition}
\label{deforbfib}
Let $f:X \longrightarrow M$ be an orbifold map\index{orbifold map} of 
complex orbifolds 
of relative dimension $\geq 1$ (i.e. $dim_{\bC}X \geq dim_{\bC}M + 1$).
We say that $f$ has a structure of an {\em orbifold fibration} 
if the following conditions hold.

\noindent {\rm (i)} \ 
Let $\{ ( \wU_i, G_i, \varphi_i, U_i) \}_{i \in \bI}$ be the orbifold charts  
of $M$.
For each $i \in \bI$, there exists a normal complex space $\wW_i$ 
and a holomorphic map $\wf_i:  \wW_i \longrightarrow  \wU_i$, 
and $G_i$ acts holomorphically and relatively on $\wf_i$ 
such that the diagram

\begin{picture}(200,71)(-8,-3)
\put(80,50){\makebox(6,6){$\wW_i$}}
\put(80,5){\makebox(6,6){$\wU_i$}} 
\put(200,51){\makebox(12,6){$f^{-1}(U_i) \cong \wW_i/G_i$}} 
\put(195,6){\makebox(12,6){$U_i  \cong \wU_i/G_i$.}} 
\put(81,43){\vector(0,-1){24}} 
\put(181,43){\vector(0,-1){24}} 
\put(93,51){\vector(1,0){64}}
\put(93,7){\vector(1,0){70}}
\put(68,26){\makebox(6,6){$\wf_i$}} 
\put(186,26){\makebox(6,6){$f$}} 
\put(121,-3){\makebox(6,6){$\varphi_i$}} 
\put(121,39){\makebox(6,6){$\psi_i$}} 
\end{picture}

\noindent 
is commutaive, where $\psi_i$ is the projection map by $G_i$.

\noindent {\rm (ii)} \ 
Each $\wW_i$ is an open complex sub-orbifold of $X$ in the following sense: 
$\wW_i$ has the orbifold charts 
$\{ ( \widehat{W}_{i,j}, H_{i,j}, \varphi_{i,j}, \wW_{i,j}) \}_{j \in \bJ(i)}$ 
(the set $\bJ(i)$ of suffixes  depends on $i$), 
and there exist a group $G_{i,j}$ containing $H_{i,j}$ as a normal subgroup  
and an exact sequence of groups 
$$
1 \longrightarrow  H_{i,j}  \longrightarrow G_{i,j} \longrightarrow G_i \longrightarrow 1, 
$$
such that the  orbifold charts of $X$ are given by
$\{ ( \widehat{W}_{i,j}, G_{i,j}, \psi_i \circ \varphi_{i,j}, 
 \psi_i \circ \varphi_{i,j}( \widehat{W}_{i,j})) \}_{i \in \bI, j \in \bJ(i)}$.
\end{definition}

The typical example of an orbifold fibration will be given in \S6.2.
Note that the set
\begin{equation}
\label{Eq0403pseudochart}
\{ ( \wW_i, G_i, \psi_i, f^{-1}(U_i) )\}_{i \in \bI}
\end{equation}
 is not an atlas of orbifold charts of $X$ in general, because 
$\wW_i$ may have singularities.
Nevertheless it is important in our discussion in \S6.2. 
See Definition \ref{defspecialorb} (i).

\begin{definition}
\label{defspecialorb}
{\rm (i)} \ We call (\ref{Eq0403pseudochart}) 
{\em the pseudo-orbifold charts} of $X$ 
with respect to the orbifold fibration $f: X \rightarrow M$.

\noindent {\rm (ii)} \ If $\wW_i$ is nonsingular for each $i$, 
i.e. if the set (\ref{Eq0403pseudochart}) gives an atlas of orbifold charts of $X$, 
we call $f: X \rightarrow M$ {\em a strong orbifold fibration}.
\end{definition}

Now the theorem (\cite[Th.6.11]{M2020}) says that the Deligne-Mumford 
compactification $\overline{M}_g$ has the orbifold charts
\begin{equation}
\label{mgbarchars}
\{ ( D_{\epsilon}(\sigma), W(\sigma), \varphi_{\sigma}, M_{\epsilon}(\sigma)) \}_{\sigma \in \cC_g/ \Gamma_g}, 
\end{equation}
where $\sigma$ moves in the curve complex modulo the action of the mapping class group 
$\Gamma_g$.
From now on, since this orbifold structure is different 
from the usual one (cf. \cite[Chap.XII]{ACG}) of $\overline{M}_g$, 
we use the symbol $\overline{M}_g^{orb}$ in order to specify it.
Over this base structure,
we have the following strong orbifold fiber space.

\begin{theorem} 
\label{orbunivdeg}
There exists a strong  orbifold fibration\index{strong  orbifold fibration}
\begin{equation}
\label{uinvdegfam}
\overline{\pi}: \overline{Y}_g^{orb} \longrightarrow  \overline{M}_g^{orb}
\end{equation}
such that the orbifold charts of $\overline{M}_g^{orb}$ 
are given by (\ref{mgbarchars}), and those of 
$\overline{Y}_g^{orb}$ are  
\begin{equation}
\label{totalunivchars}
\{ ( X_{\epsilon}(\sigma), W(\sigma), \psi_{\sigma}, 
Y_{\epsilon}(\sigma)) \}_{\sigma \in \cC_g/ \Gamma_g}
\end{equation}
given in Theorem \ref{KuranishioverCont}.
\end{theorem} 

{\it Proof} \quad 
It suffices to prove the compatibility condition (ii) of 
Def.\ref{defspecialorb}.
Let $[S] \in \overline{M}_g$ be an isomorphism class, and let 
$p_i=[S, w_i], \;w_i: \Sigma_g(\sigma) \longrightarrow S$ ($i=1,2$) be 
two marked stable curves such that 
the point $p_i$ belongs to $D_{\epsilon}(\rho_i)$ for $\rho_i \leq \sigma$.
By Theorem \ref{KuranishioverCont}, 
there exists an open neighborhood $U_i \subset D_{\epsilon}(\rho_i)$ of $p_i$ 
such that the restricted family $\pi_i:X_i \longrightarrow U_i$ 
of $\pi_{\sigma}: X_{\epsilon}(\rho_i) \longrightarrow D_{\epsilon}(\rho_i)$ 
over $U_i$  coincides with the standard Kuranishi family of $S$.
By the uniqueness of the standard Kuranishi family modulo isomorphism, 
there exist biholomorphic maps $\overline{h}:X_1 \longrightarrow X_2$ 
and $h: U_1 \longrightarrow U_2$ which satisfy
$h \circ \pi_1=\pi_2 \circ \overline{h}$ 
after a  suitable shrinking  of $U_i$.
Therefore, the desired compatibility condition is satisfied.

The other required properties for the strong  orbifold fibration
are given in Theorem \ref{KuranishioverCont}.
In particular, the normal analytic structure of 
$\overline{Y}_g^{orb}$ is given by patching 
those of $Y_{\epsilon}(\sigma)$'s obtained by Cartan's theorem 
via the above local biholomorphic maps.
\qed

\begin{definition}
\label{defunivdeg}
We call 
$\overline{\pi}: \overline{Y}_g^{orb} \longrightarrow  \overline{M}_g^{orb}$
the universal degenerating family of Riemann surfaces\index{universal degenerating family of Riemann surfaces}
 of genus $g$.
\end{definition}

We believe that the family $\overline{\pi}$ is {\it universal} in the sense that 
every orbifold fibration with a Riemann surface of genus $g$ as a general fiber 
can be pulled back from this universal orbifold fibration. Our present 
achievement is, however, rather modest, and we have proved the universality 
of $\overline{\pi}$ only for fibered complex surfaces, namely, for orbifold fibrations  
whose base spaces are  
of dimension $1$. See \S7.
For a precise definition of \lq\lq orbifold pull-back\rq\rq, see the following.

\begin{definition}
\label{definition0403no3}
Let $f:X \rightarrow M$ 
and 
$f':X' \rightarrow M'$  be orbifold fibrations\index{orbifold fibration} 
and $h:M' \rightarrow M$ be an orbifold map\index{orbifold map}.
We say that  $f'$ is {\em the orbifold pull back\index{orbifold pull back}
 from $f$ via $h$}  
if the following conditions are satisfied:

\noindent {\rm (i)} \ 
Let 
$\{ ( \wU_i, G_i, \varphi_i, U_i) \}_{i \in \bI}$ and 
$\{ ( \wV_j , H_j,  \phi_j, V_j) \}_{j \in \bJ}$ 
be  the orbifold charts of $M$ and 
$M'$ respectively.
For any $j \in \bJ$, there exist some $i \in \bI$ depending on $j$ 
and a holomorphic map $h_{ji}=h|_{V_j}:V_j \rightarrow U_i$.
Let $\wh_{ji}:\wV_j \rightarrow  \wU_i$ be the lifting of 
$h_{ji}$ and 
$\varphi_i^{(j)}=\varphi_i |_{\wU_i^{(j)}}: 
\wU_i^{(j)}=\wh_{ji}(\wV_j) \longrightarrow U_i^{(j)}=h_{ji}(V_j)$ 
be the restriction map to the image, i.e.

\begin{picture}(230,71)(-35,-6)
\put(80,50){\makebox(6,6){$\wV_j$}}
\put(52,5){\makebox(12,6){$\wV_j/H_j \cong V_j$}} 
\put(190,51){\makebox(12,6){$\wU_i^{(j)} \subset \wU_i$}} 
\put(215,6){\makebox(12,6){$U_i^{(j)} \subset U_i  \cong \wU_i/G_i$}} 
\put(81,43){\vector(0,-1){24}} 
\put(181,43){\vector(0,-1){24}} 
\put(93,51){\vector(1,0){72}}
\put(93,7){\vector(1,0){72}}
\put(68,27){\makebox(6,6){$\phi_j$}} 
\put(190,31){\makebox(6,6){$\varphi_i^{(j)}$}} 
\put(125,-5){\makebox(6,6){$h_{ji}$}} 
\put(125,38){\makebox(6,6){$\wh_{ji}$}} 
\end{picture}

\noindent 
is a commutative diagram.
Then there exists an injective group homomorphism
\begin{equation}
\label{equation0403nox2}
H_j \hookrightarrow G_i
\end{equation}
such that the subgroup $H_j$ of $G_i$ acts on 
$\wU_i^{(j)}$ holomorphically.

\noindent {\rm (ii)} \ 
There exists a lifted orbifold map $k:X' \rightarrow X$ of $h$, i.e. 
$f \circ k=h \circ f'$, which has the following property: 
Let $\{ ( \wW_i, G_i, \psi_i, f^{-1}(U_i)) \}_{i \in \bI}$ and 
$\{ ( \wT_j , H_j,  \omega_j, (f')^{-1}(V_j)) \}_{j \in \bJ}$ 
be  the pseudo-orbifold charts of $X$ and $X'$ respectively.
Let $\wk_{ji}:\wT_j \rightarrow  \wW_i$ be the lifting of 
$k_{ji}=k|_{(f')^{-1}(V_j)}:(f')^{-1}(V_j) \longrightarrow f^{-1}( U_i)$,
and 
$\psi_i^{(j)}=\psi_i |_{\wW_i^{(j)}}: 
\wW_i^{(j)}=\wk_{ji}(\wT_j)= \psi_i^{-1}(\wU_i^{(j)}) 
\longrightarrow \wU_i^{(j)}$
be the restriction map.
Then the group $H_j$ relatively acts on 
$\psi_i^{(j)}$ 
such that the diagram

\begin{picture}(200,71)(-25,-3)
\put(80,50){\makebox(6,6){$\wT_j$}}
\put(80,5){\makebox(6,6){$\wV_j$}} 
\put(191,50){\makebox(12,6){$\wW_i^{(j)} \subset \wW_i$}} 
\put(191,5){\makebox(12,6){$\wU_i^{(j)} \subset \wU_i$}} 
\put(81,43){\vector(0,-1){24}} 
\put(181,43){\vector(0,-1){24}} 
\put(93,51){\vector(1,0){71}}
\put(93,7){\vector(1,0){71}}
\put(68,28){\makebox(6,6){$\omega_j$}} 
\put(191,30){\makebox(6,6){$\psi_i^{(j)}$}} 
\put(123,-5){\makebox(6,6){$k_{ji}$}} 
\put(123,38){\makebox(6,6){$\wk_{ji}$}} 
\end{picture}

\noindent 
expresses the fiber product compatible with the action of $H_j$, i.e.
$\wT_j$ is isomorphic to 
$\wW_i^{(j)} \times_{\wU_i^{(j)}} \wV_j$ such that 
$\alpha \circ \wk_{ji}=k_{ji} \circ \alpha$ for any $\alpha \in H_j$.
\end{definition}


\section{Automorphisms of stable curves and cyclic equisymmetric strata  on $\overline{M}_g^{orb}$}
\label{Pullbacksurface}

We study an automorphism $\varphi$ of a stable curve $S$ 
and the associated cyclic branched covering from several points of view.
In particular, we define the numerical data ${\rm Num}(\varphi)$  which consist  of two 
kinds of data, i.e.
the  action on the dual graph of $S$ and 
the system of branch data of cyclic branched coverings 
which every irreducible component of $S$ naturally has.
We study the locus $T_{\sigma}^{[\varphi]}$ on the boundary charts of 
$\overline{M}_g^{orb}$ consisting of marked stable curves with the automorphisms 
of this type of numerical data, and describe its structure  
 in terms of pointed Teichm\"{u}ller spaces of lower genera 
(Theorem \ref{stableequithm}).
This result is an extension of Harvey--Broughton's  theorem about 
the equisymmetric strata on $T_g$ and $M_g$, 
to the boundaries for cyclic groups.
 
In \S5.1, we review some known results about automorphisms of Riemann surfaces.
First, we  review the notion of total valency essentially due to Nielsen \cite{Nielsen1}  
and Harvey \cite{H1966}, whic provides precise information on the branch data 
of  the cyclic covering associated with the automorphism.
Secondly, we review Harvey--Broughton's equisymmetric strata on $T_g$ and $M_g$ 
(\cite{H1971}, \cite{Broughton1990}) in the case of cyclic groups.

In \S5.2, as a modified discussion of Eichler's trace formula, 
we propose a method to determine the characters of the representation 
of  automorphisms of pointed Riemann surfaces into the space of 
logarithmic quadratic differential forms.

In the terminology of augmented Teichm\"{u}ller theory, 
the true boundary $T(\sigma)$  on the boundary chart of  $\overline{M}_g^{orb}$ 
should be  called the {\it little Teichm\"{u}ller space}.
In \S5.3, we interpret this notion by the cohomological terminologies of 
Kuranishi spaces of stable curves.

The aim of  \S5.4 is to define  ${\rm Num}(\varphi)$ naturally. 
The  point is to analyze the cyclic branched covering 
$\pi_{\varphi}: S \rightarrow W$ associated with $\varphi$.
Although the base $W$ is a nodal Riemann surface, the dual graph of $W$ is not a simple 
quotient graph of the dual graph of $S$.
In fact, we extend the notion of graphs to those admitting open edges, and then 
define the notion of compact quotient graphs as the desired one.
By combining these data of the graphs and the system of the total valencies 
which every component of $S$ naturally has, we have the definition of ${\rm Num}(\varphi)$.

In \S5.5, we prove the structure theorem of the equisymmetric strata 
$T_{\sigma}^{[\varphi]}$ on $T({\sigma})$.
The space $T_{\sigma}^{[\varphi]}$ is described locally by the Kuranishi space 
of each normalized component of $W$, and globally by the pointed 
Teichm\"{u}ller spaces. 
Our basic method is, 
in the case of Riemann surfaces, closely related to  the one in \cite{V1994} or \cite[\S\S 4,5]{MN}.
In the case of stable curves, 
the pioneering work of Terasoma 
\cite{Terasoma1998} already 
shows the connectivity of the moduli space  $\overline{M}_g^{[\varphi]}$ in our teminology 
by using the level structures. 

In \S5.6, we define a special system of local coordinates 
around a point of $T_{\sigma}^{[\varphi]}$. 
These coordinates consist of the eigenvectors 
for the action of $\varphi$
on the standard chart of the Kuranishi space of $S$.  
Harris--Mumford \cite[\S1]{HM} intrinsically used 
this type of coordinates  
systematically for a certain local analysis of $\overline{M}_g$.


\subsection{Automorphisms of Riemann surfaces and equisymmetric strata}
\label{autoriemannequi}

In this subsection, we review some results about the automorphisms of Riemann surfaces,  
the equisymmetric stratification on  Teichm\"{u}ller space $T_g$ and 
the moduli space $M_g$ (\cite{H1971}, \cite{Broughton1990}).

We consider a Fuchsian group\index{Fuchsian group} $F$, i.e. a discrete subgroup of ${\rm Aut}(\cH) \cong {\rm PSL}(2, \bR)$ 
where $\cH$ is the upper half plane, such that $\cH/F$ is compact.
As an abstract group, $F$ is generated by 
$a_1, b_1, \cdots, a_{\overline{g}}, b_{\overline{g}}, x_1, \cdots, x_s$ with the relations 
$$
x_1 \cdots x_s \prod_{i=1}^{\overline{g}} [a_i, b_i] =1, \:\:\: x_i^{\lambda_i}=1 \:\:(1 \leq i \leq s).
$$
The ordered set ($\overline{g}; \lambda_1, \cdots, \lambda_s$) is called  {\it the signature} of $F$.
We assume that there exists a Fuchsian group $K$ with signature $(g, -)$ (where $-$ means 
that $\{ x_1, \cdots \}$ is empty) and exact sequence
\begin{equation}
\label{Fuchianex}
1 \rightarrow  K  \xrightarrow{i_*}  F  \xrightarrow{j_*}   
G_{\varphi} \rightarrow 1
\end{equation}
such that $G_{\varphi}=\langle \varphi \rangle \cong \bZ/n\bZ$ is  the 
cyclic group of order $n$ generated by $\varphi$ and 
$i_*(K)$ is a normal subgroup of $F$.
In this case, $G_{\varphi}$ is geometrically characterized as follows. 
The Riemann surface $S \cong \cH/K$ of genus $g$ has an analytic automorphism 
$\varphi:S \longrightarrow S$ of order $n$ such that the associated branched covering 
$\pi_{\varphi}:S \longrightarrow W=S/ G_{\varphi}$ 
has the following properties: 
The genus of $W$ coincides with $\overline{g}$. 
Let $P_1, \cdots. P_s \in W$ be the branch points for $\pi_{\varphi}$.
Then $\lambda_i$  ($1 \leq i \leq s$) coincides with 
the ramification index  at the point  
$\widetilde{P}_i \in \pi_{\varphi}^{-1}(P_i)$.
The map $\varphi^{n/\lambda_i}$ fixes $\widetilde{P}_i$ and rotates 
the disc neighborhood of $\widetilde{P}_i$ by the angle 
$2\pi\delta_i/\lambda_i$.
Let $1 \leq \sigma_i \leq \lambda_i-1$ 
be the natural number with $\sigma_i\delta_i \equiv 1$ (${\rm mod}\: \lambda_i$).
The triple ($m_i, \lambda_i, \sigma_i$) is called the {\it valency} 
of $\varphi$ at $\widetilde{P}_i$ 
(\cite{Nielsen1}, \cite[Def.~1.5]{MM}), and sometimes is written by $ \sigma_i/\lambda_i$.
Then:

\noindent {\rm (a1)} (Hurwitz formula) \:
$2(g-1)/n=2(\overline{g}-1)+\sum_{i=1}^{s} (1-1/\lambda_i)$, 

\noindent {\rm (a2)} (Nielsen \cite[(4.6)]{Nielsen1}) \: 
$\sum_{i=1}^s \sigma_i/\lambda_i$ 
is an integer,

\noindent {\rm (a3)} (Wiman \cite{Wiman}) \: $n \leq 4g+2$,

\noindent {\rm (a4)} (Harvey \cite[Th.4]{H1966}) \: 
We set $M={\rm lcm} (\lambda_1, \cdots, \lambda_s)$. Then

{\rm (a4-1)} \: 
${\rm lcm} (\lambda_1, \cdots, \widehat{\lambda_i}, \cdots, \lambda_s)=M$ 
for all $i$, where $\widehat{\lambda_i}$ denotes the omission of $\lambda_i$.

 {\rm (a4-2)}\: $M$ divides $n$, and if $\overline{g}=0$, then $M=n$.

{\rm (a4-3)}\: $s \not=1$, and, if $\overline{g}=0$, then $s \geq 3$.

 {\rm (a4-4)}\: If $2|M$, the number of $\lambda_1, \cdots, \lambda_s$ 
which are divisible by the maximal power of 
$2$ dividing $M$ is even.

\medskip

We symbolically write these data by 
\begin{equation}
\label{totalvalency}
{\rm TV}(\varphi)=
\left( g, \overline{g}, n: \frac{\sigma_1}{\lambda_1} + \cdots +\frac{\sigma_s}{\lambda_s} \right), \:\:
{\rm TV}^c(\varphi)=
\left( g, \overline{g}, n: \left\{ \frac{\delta_1}{\lambda_1}, \cdots, 
\frac{\delta_s}{\lambda_s} \right\} \right)
\end{equation}
and call ${\rm TV}(\varphi)$ (resp.${\rm TV}^c(\varphi)$)
 {\it the total valency}\index{total valency}
  (resp.~{\it the total co-valency}) 
 of $\varphi$. Note that the total (co-)valency is determined by $j_*$ of (\ref{Fuchianex}) from \cite[Th.7, Lem.6]{H1971}.

\begin{remark}
\label{lowgenusrem}
{\rm (i)}
The valency was introduced by Nielsen \cite{Nielsen1} as the unique essential conjugacy invariant 
of periodic maps (i.e. finite automorphisms) in $\Gamma_g$. A periodic map is realized as an 
analytic automorphism of a certain complex structure on $\Sigma_g$ by a standard augument 
or as a corollary of Kerchoff \cite{Ke}. Therefore, the total valency may be considered as an 
invariant of both periodic maps in $\Gamma_g$ and anlytic automorphisms.\par
\noindent{\rm (ii)} 
Even if $g=0,1$, a finite automorphism $\varphi$ satisfies the conditions 
{\rm (a1)} $\sim$ {\rm (a3)} and we   use the same terminology 
(\ref{totalvalency}) in these cases.\par
\noindent {\rm (iii)} 
For the classification of ${\rm TV}(\varphi)$ 
in the case where $g=2, 3$, see e.g. \cite[p.199]{AI}.
\end{remark}

As discussed in  \cite{H1971},
the choice of the inclusion map $i_*$ 
in (\ref{Fuchianex}) determines the  
Teichm\"{u}ller marking, and 
the choice of the surjective map $j_*$ determines  the generators of 
$G_{\varphi}$, which determine the data {\rm (a2)} in turn.

The conditions {\rm (a1)} $\sim$ {\rm (a4)} are not only 
necessary conditions but also  sufficient conditions 
for the existence of automorphisms.
This point  is widely discussed from 
the moduli theoretic viewpoint by Harvey and Broughton as follows.
We fix  the numerical data  (\ref{totalvalency}) for $g \geq 2$, 
and  consider the locus in $T_g$  (or $M_g$) of 
Riemann surfaces which have this type of automorphisms.
Note that, in the following results  (I) and (II), 
the case where $g=2$ and  $\varphi$ is the hyperelliptic involution is excluded 
as an exceptional case:
 
\medskip

\noindent
{\rm (I)}  (\cite[Th.2 and (6) in p.392]{H1971}, \cite[Prop.2.5]{Broughton1990}) 
Let $\varphi : S \to S$ be a periodic map in $\Gamma_g$. Let $T_g^\varphi$ be the subset of $T_g$ consisting of the 
points $p=[S,f] \in T_g$ which are fixed by the action of $\varphi$.  We define two periodic maps $\varphi$ and $\psi$ 
to be {\it equivalent} if their total valencies are the same: ${\rm TV}(\varphi)={\rm TV}(\psi)$. As we remarked in 
Remark \ref{lowgenusrem} (i), 
Nielsen \cite{Nielsen1} proved that this is equivalent to saying that $\varphi$ and $\psi$ are conjugate. Let $[\varphi]$ denote 
the equivalence class to which $\varphi$ belongs, and we define $T_g^{[\varphi]}$ as
$$T_g^{[\varphi]} =\bigcup_{\psi \in [\varphi]}T_g^{\psi}.$$
Recall that $\pi_{\varphi}: S \to W=S/G_{\varphi}$ is the branched covering associated with the analytic action $\varphi: S \to S$ 
of order $n$, and $\overline{g}$ is the genus of $W$. $P_1, \cdots, P_s \in W$ are the branch points for $\pi_{\varphi}$. 
Then $T_g^{\varphi}$ is real analytically isomorphic to the Teichm\"uller space $T_{\overline{g}, s}$ of $s$-pointed Riemann surfaces of 
genus $\overline{g}$.  Since $T_g^\varphi$ is a connected component $(T_g^{[\varphi]})^{(0)}$  of $T_g^{[\varphi]}$, we have 
$$T_g^{\varphi}\cong (T_g^{[\varphi]})^{(0)} \cong T_{\overline{g},s}\;\;\;\;\;(\text{real analytically}).$$
Each connected component of $T_g^{[\varphi]}$ is $T_g^\psi$ for some $\psi \in [\varphi]$, and we can say the
same thing for $T_g^\psi$. Thus the space $T_g^{[\varphi]}$ is a countable union of $(3\overline{g}-3+s)$-dimensional complex manifolds 
so that each of them is real analytically isomorphic to $(T_g^{[\varphi]})^{(0)}$.\par

 \noindent {\rm (II)} (\cite[Th.2.2]{Broughton1990}) \:
 The locus $M_g^{[\varphi]}$ of the isomorphism classes of Riemann surfaces which have automorphisms  
 whose total valencies coincide with $\rm{TV}(\varphi)$ 
 is a closed irreducible algebraic subvariety of $M_g$.

 \medskip
 
 The loci $T_g^{[\varphi]}$ and $M_g^{[\varphi]}$ are called 
 {\it the equisymmetric strata}\index{equisymmetric strata}
  for $[\varphi]$ of $T_g$ and $M_g$ respectively.
 Note that, in the excluded case where $g=2$ and $\varphi$ is the hyperelliptic involution, we have 
 $T_2^{[\varphi]}=T_2$ and $M_2^{[\varphi]}=M_2$.
 
 \begin{remark}
 \label{Rem0501no2b}
 For the study of the equisymmetric strata on $T_g$, Kuribayashi's method  (\cite{Kuribayashi1966})
using invariant quadratic differentials is important.
 Recently, Takamura (\cite{Tak2019}) and Hirakawa and Takamura (\cite{HiTa}) 
 developed this type of argument and gave a method for 
 the detailed analysis of the stratification for various group actions 
 on Riemann surfaces.
 \end{remark}


\subsection{Logarithmic quadratic representation of automorphisms}
\label{Meromorphicplurirep}

Here we discuss a representation 
of an automorphism of a pointed Riemann surface to the logarithmic quadratic differentials\index{logarithmic quadratic differential} 
for later use.

Let $\varphi:S \longrightarrow S$ be an automorphism of a Riemann surface of 
genus $g \geq 0$ with the given total (co-)valency  (\ref{totalvalency}), and 
$\pi_{\varphi}:S \longrightarrow W=S/G_{\varphi}$ be the associated 
cyclic branched covering (see also Remark \ref{lowgenusrem} (ii)).
We consider the vector space 
\begin{equation}
\label{eigenspaceV}
V=H^0(S, 2K_S+\sum_{i=1}^{s'} \pi_{\varphi}^{-1}(P_i) 
+\sum_{i=1}^{t} \pi_{\varphi}^{-1}(P_{s+i}))
\end{equation}
where $\{ P_1, \cdots, P_{s'} \}$ is a 
subset (it may be empty) of the set of branch points 
$\{ P_1, \cdots, P_{s'},$ $  P_{s'+1}, \cdots, P_{s} \}$ ($ s' \leq s$) 
for $\pi_{\varphi}$, and 
$\{ P_{s+1}, \cdots, P_{s+t} \}$
is a subset of the set of un-branched points for $\pi_{\varphi}$ 
(which  is intended to be the set of possible poles, and may be empty).  We assume
\begin{equation}
\label{assumpcardpoints}
2\overline{g}-2+s+t >0.
\end{equation}
An element $v \in V$ is considered as a logarithmic quadratic 
differenticial on $S$ whose poles are at most in 
$\sum_{i=1}^{s'} \pi_{\varphi}^{-1}(P_i) 
+\sum_{i=1}^{t} \pi_{\varphi}^{-1}(P_{s+i})$.
Since $\varphi$ acts on each fiber of $\pi_{\varphi}: S \to W$ as a 
permutation of points, $\varphi$ naturally acts on the differential forms $v$ by 
\begin{equation}
\label{actsonform}
\varphi(v)=v \circ \varphi^{-1} \in V
\end{equation}
(similarly to \cite[p.269]{FK}),  
and  $\varphi$ induces a linear automorphism  of $V$.
We call it {\it the representation of $\varphi$ to the logarithmic quadratic 
forms $V$}.
Fundamental facts of the representation to 
holomorphic differential forms discussed in \cite[\S V2]{FK} are easily extended to 
this type of representation to logarithmic forms.
For example, its eigenvalues are $n$-th roots of unity.

The proof of the following Proposition is 
a slight modification of I. Guerrero's argument written in 
\cite[p.274--277]{FK}.
For a rational number $x$, we write $[x]$ the maximal integer not exceeding $x$ and 
$\{ x \}=x- [x]$ its fractional part.
We also use the symbol $\be(x)=e^{2\pi i x}$.

\begin{prop}
\label{eigenspaceformula}
Let $\varphi: S \longrightarrow S$ be an automorphism of order $n$ of a Riemann surface $S$ 
with the total (co-)valency (\ref{totalvalency}).
Then, under the assumption  (\ref{assumpcardpoints}), 
the dimension of the eigenspace  
of eigenvalue $\be(\alpha/n)$ $(0 \leq \alpha \leq n-1)$ for the 
the action $\varphi$ on $V$ defined by  (\ref{actsonform})
is given as follows. 
$$
h^0\left(S, 2K_S+\sum_{i=1}^{s'} \pi_{\varphi}^{-1}(P_i) +\sum_{i=1}^{t} \pi_{\varphi}^{-1}(P_{s+i}) \right)_{\be(\alpha/n)}
\:\:\:\:\:\:\:\:\:\:\:\:\:\:\:\:\:\:\:\:\:\:\:\:\:\:\:\:\:\:\:\:\:\:\:\:
\:\:\:\:\:\:\:\:\:\:\:\:\:\:\:\:\:\:\:\:\:\:\:
$$
\vspace{-0.3cm}
$$
=3\overline{g}-3+2s +t
-\sum_{i=1}^{s'} \left( 
\left\{ \frac{ -\alpha \sigma_i-1}{\lambda_i} \right\} 
+\frac{1}{\lambda_i}
\right)
-\sum_{i=s'+1}^{s} \left( 
\left\{ \frac{ -\alpha \sigma_i-2}{\lambda_i} \right\} 
+\frac{2}{\lambda_i}
\right).
$$
\end{prop}

{\it Proof} \quad 
If $1 \leq i \leq s'$ or $s+1 \leq i \leq s+t$, then we set $r_i=1$.
If $s'+1 \leq i \leq s$, then we set $r_i=0$.
We denote the eigenspace of eigenvalue $\be(\alpha/n)$ by  
$V_{\alpha}=H^0(S, 2K_S+\sum_{i=1}^{s+t} r_i \pi_{\varphi}^{-1}(P_i))_{\be(\alpha/n)}$ 
for simplicity.
Moreover, we set 
\begin{equation}
\label{inequality}
\lambda_i=1, \delta_i=\sigma_i=0, \;\;\;
\text{for $s+1 \leq i \leq s+t$.}
\end{equation}
For $1 \leq i \leq s+t$, 
we denote the fiber by
$$
\pi_{\varphi}^{-1}(P_i)=\coprod_{1 \leq j \leq n/\lambda_i} P_j^{(i)}.
$$
Assuming  $V_{\alpha} \not= \emptyset$, we fix an element $v_0 \in V_{\alpha}$.
Then for any $v \in V_{\alpha}$, 
the element $v/v_0$ is $\varphi$-invariant, and it is a 
meromorphic function on $W=S/G_{\varphi}$.
Therefore, there exits a divisor $D_{\alpha}$ on $W$ such that 
\begin{equation}
\label{wtohzero}
V_{\alpha} \cong H^0(W, D_{\alpha}).
\end{equation}

We determine $D_{\alpha}$ explicitly.
By using a local coordinate $z$ around 
$P_j^{(i)}$, 
the Laurant expansion of  $v$ is written as  
$v=\sum_{k \geq 0} A_kz^{k-r_i}dz^2$ for $A_k \in \bC$.
Since the map $\varphi^{-n/\lambda_i}$ is written here by  $z \mapsto \be(-\delta_i/\lambda_i)z$, 
 we have
$$
 (\varphi^{n/\lambda_i})^*v=
 \sum_{k \geq 0} A_k\be \left( \frac{-\delta_i(k+2-r_i)}{\lambda_i} \right) z^{k-r_i} 
 dz^2
$$
which coincides with 
$\sum_{k \geq 0}  \be(\alpha/\lambda_i)A_k z^{k-r_i} dz^2$ 
by definition.
Therefore, $A_k=0$ for $-\delta_i(k+2-r_i) \not\equiv \alpha$ (${\rm mod}\:\lambda_i$), i.e. 
$k \not\equiv -\alpha \sigma_i-2+r_i$  (${\rm mod}\:\lambda_i$).
It follows that 
$$v=\sum_{k \in \bZ_{K}} A_kz^{k-r_i}dz^2, \;\;\text{ where $\bZ_{K}=\left\{ k=\lambda_ia-\alpha \sigma_i-2 +r_i \geq 0\; | \;\exists a \in \bZ \;\right\}. $}$$

Let $b_i$ be the vanishing order of 
$v_0$ at $P_j^{(i)}$, 
which is automatically independent of $j$.
Then there exists some $\widetilde{b}_i \in \bZ$ such that 
\begin{equation}
\label{alphaexpress}
b_i= \lambda_i \widetilde{b}_i-\alpha \sigma_i-2 +r_i \geq 0.
\end{equation}
Let $Q_1, \cdots, Q_u$ be the images by $\pi_{\varphi}$ of zeros of $v_0$ 
except for $P_1, \cdots, P_{s+t}$. Put 
$\pi_{\varphi}^{-1}(Q_i)=\coprod_{1 \leq j \leq n} Q_j^{(i)}$.
Let $\beta_i$ be the vanishing order of $v_0$ at $Q_j^{(i)}$.
For a meromorphic function $\widetilde{h}$ on $W$, 
the element $(\widetilde{h} \circ \pi_{\varphi}) \cdot v_0$ belongs to 
$V_{\alpha}$ if and only if   $\widetilde{h}$ satisfies 
${\rm ord}_{P_i} \widetilde{h} \geq 
- b_i/\lambda_i
, \: {\rm ord}_{Q_i} \widetilde{h} \geq -\beta_i$ 
and is holomorphic outside $P_i$'s and $Q_i$'s.
Since $-\alpha \sigma_i-2 +r_i < 0$,  the  integral divisor 
which satisfies these conditions should be $D_\alpha = \sum_{i=1}^{s+t}(\widetilde{b}_i+\left[(-\alpha\sigma_i-2+r_i)/\lambda_i\right])P_i
+\sum_{i=1}^u\beta_iQ_i$. In particular,
\begin{equation}
\label{divisorexpress}
{\rm deg} \:D_{\alpha}=\sum_{i=1}^{s+t} 
\left(\widetilde{b}_i+\left[ \frac{-\alpha \sigma_i-2 +r_i }{\lambda_i} \right]\right)
+\sum_{i=1}^u \beta_i. 
\end{equation}

We rewrite the expression (\ref{divisorexpress}) of ${\rm deg}D_\alpha$  so that it is independent of the choice of $v_0$. \\
First we clearly have ${\rm deg}\; v_0=\sum_{i=1}^{s+t}(n/\lambda_i)b_i+n\sum_{i=1}^{u}\beta_i$. On the other hand, 
since $v_0$ belongs to $V_\alpha$, we have ${\rm deg} v_0 = {\rm deg}\left( 2K_S+ \sum_{i=1}^{s+t}r_i  \pi_{\varphi}^{-1}(P_i)\right) 
= 4(g-1)+\sum_{i=1}^{s+t}(n/\lambda_i)r_i$.    Therefore it follows from (\ref{alphaexpress}) that

\begin{equation}
\label{middleformula}
\frac{{\rm\; deg}  v_0}{n}=\sum_{i=1}^{s+t}\left(\widetilde{b_i} +\frac{-\alpha\sigma_i-2+r_i}{\lambda_i}\right) +\sum_{i=1}^t \beta_i
=\frac{4(g-1)}{n}+\sum_{i=1}^{s+t}\frac{r_i}{\lambda_i}.
\end{equation}                                                                                    
The Riemann--Hurwitz formula for the covering $\pi_{\varphi}$ says that 
\begin{equation}
\label{rhformula}
\frac{2g-2}{n}=2\overline{g}-2+\sum_{i=1}^{s} \left( 1-\frac{1}{\lambda_i} \right).
\end{equation}
Therefore,  
from (\ref{inequality}),  (\ref{divisorexpress}),  (\ref{middleformula}) and (\ref{rhformula}), we obtain 
$$
{\rm deg}\:D_{\alpha}
=\sum_{i=1}^{s+t} 
\left( \widetilde{b}_i+\frac{-\alpha \sigma_i-2 +r_i }{\lambda_i}
 \right)  
+\sum_{i=1}^t \beta_i 
- \sum_{i=1}^{s+t} \left\{ \frac{-\alpha \sigma_i-2 +r_i }{\lambda_i} \right\}
\:\:\:\:\:\:\:\:\:\:\:\:\:
$$
$$
=4\overline{g}-4+2 \sum_{i=1}^{s} \left( 1-\frac{1}{\lambda_i} \right)
+ \sum_{i=1}^{s+t} \left( \frac{r_i}{\lambda_i}-\left\{ \frac{-\alpha \sigma_i-2 +r_i }{\lambda_i} \right\} \right)
$$
\begin{equation}
\label{degreefinal}
=4\overline{g}-4+2s+t
-\sum_{i=1}^{s} \left( 
\left\{ \frac{ -\alpha \sigma_i-2+r_i}{\lambda_i} \right\} 
+\frac{2-r_i}{\lambda_i}
\right).
\:\:\:\:\:\:\:\:\:
\end{equation}
By (\ref{assumpcardpoints}) and (\ref{degreefinal}), we have
${\rm deg}(K_W \otimes D_{\alpha}^{-1})
\leq 2-2\overline{g}-s-t<0$. Hence

\begin{equation}
\label{h1vanishing}
H^0(W, K_W \otimes D_{\alpha}^{-1})=0.
\end{equation}
From the Riemann--Roch  formula, the Serre duality and (\ref{wtohzero}), (\ref{degreefinal}), (\ref{h1vanishing}), we have 
$${\rm dim} V_\alpha = {\rm dim} H^0(W, D_\alpha)={\rm deg} D_\alpha -\overline{g}+1=
3\overline{g}-3+2s+t-\sum_{i=1}^s\left(\left\{\frac{-\alpha\sigma_i-2+r_i}{\lambda_i}\right\}+\frac{2-r_i}{\lambda_i}\right).$$
This equality coincides with the desired one. \qed

\bigskip

We consider ($S, \bP$) as a $k'$-pointed Riemann surface, where 
$\bP=\sum_{i=1}^{s'} \pi_{\varphi}^{-1}(P_i) 
+\sum_{i=1}^{t} \pi_{\varphi}^{-1}(P_{s+i})$  in 
(\ref{eigenspaceV}) 
and $k'=\sharp(\bP)$ (cardinality), and write the set of the eigenvalues for the action of $\varphi$ on  
the ($3g-3+k'$)-dimensional vector space $V=H^0(S, 2K_S+\bP)$ by 
\begin{equation}
\label{equation0502newno1}
\left\{ \be \left( \frac{\theta_1}{n} \right), \be \left( \frac{\theta_2}{n} \right), \cdots, 
\be \left( \frac{\theta_{3g-3+k'}}{n} \right) \right\} \:\:\:\:\: 
(0 \leq \theta_1 \leq \theta_2 \leq \cdots \leq \theta_{3g-3+k'} \leq n-1).
\end{equation}
Here each eigenvalue is counted as many times 
as the dimension of its eigenspace.
Then:
\begin{definition}
\label{definition0502newno1}
By using (\ref{equation0502newno1}), 
we define the ordered set of the 
{\em log-quadratic characters}\index{log-quadratic character} for  
the automorphism $\varphi$ of 
the pointed Riemann surface $(S, \bP)$ as 
\begin{equation}
\label{equation0502newno2}
{\rm Ch}_{\varphi}H^0(S,2K_S+\bP)=
\left\{  \frac{\theta_1}{n},  \frac{\theta_2}{n}, \cdots, 
\frac{\theta_{3g-3+k'}}{n} \right\}.
\end{equation}
We also define the ordered set of the {\em eigenbasis} 
$\{ v_1. \cdots, v_{3g-3+k'} \}$ of 
$V$ as the basis consisting of eigenvectors of each of 
${\rm Ch}_{\varphi}H^0(S,2K_S+\bP)$, i.e. 
\begin{equation}
\label{equation0502newno3}
\varphi^*(v_i)= \be \left( \frac{\theta_i}{n} \right)v_i \:\:\:\:\:\:\: 
(1 \leq i \leq 3g-3+k').
\end{equation}
\end{definition}

\begin{remark}
\label{remark0502newno1}
Since the dimension of the eigenspace for some eigenvalue 
might be strictly greater than $1$, 
the eigenbasis of $V$ is not unique in general.
\end{remark}

\begin{example}
\label{example0502no1}
Let $\varphi:(S, P) \rightarrow (S, P)$ 
be the automorphism of order $7$ of the one-pointed genus $3$ Riemann surface 
with the total valency 
$({\bf 6/7}+6/7+ 2/7, \bar{g}=0)$ 
where ${\bf 6/7}$ is attached to $P$.
Since $s=3$, $s'=1$ and $t=0$, it follows from Prop.~\ref{eigenspaceformula} 
that 
$$
h^0\left(C, 2K_C+P \right)_{\be(\nu/7)}
=3-
\left( \left\{ \frac{-6\nu-1}{7} \right\}+\frac{1}{7} \right)
-\left( \left\{ \frac{-2\nu-2}{7} \right\}+\frac{2}{7} \right)
-\left( \left\{ \frac{-6\nu-2}{7} \right\}+\frac{2}{7} \right)
$$
$$
\:\:\:\:\:\:\:\:\:\:\:\:\:\:\:
=0,1,2,1,1,1,1 \:\:\:\:\:
(\nu=0,1,2,3,4,5,6, \:\: {\rm respectively}).
$$
Hence ${\rm Ch}_{\varphi}H^0(S,2K_S+P)=
\{ 1/7, 2/7, 2/7, 3/7, 4/7, 5/7, 6/7 \}$.
\end{example}

For other examples, see Example  \ref{example0604no2} or \cite[\S 1]{HM}.


\subsection{The little Teichm\"{u}ller space 
in an orbifold chart of $\overline{M}_g^{orb}$}
\label{litteTeich}
In this subsection, we interpret the 
little Teichm\"{u}ller spaces in the augumented Teichm\"{u}ller theory 
(cf.\cite[\S7]{HK})
by the language of Kuranishi families,  
using the facts in \S \ref{Subsec0301}.

We follow the notation of \S \ref{ContversusKur}.
The maximal codimensional strata 
$T(\sigma) \subset 
D_{\epsilon}(\sigma) \setminus T_g$ 
is written by the extended Fenchel--Nielsen coordinates (\ref{extendedFN}) as 
$$
\ell_1= \cdots =\ell_k=0, \:\: \ell_j > 0 \:\:\:(k+1 \leq j \leq 3g-3).
$$
The space $T(\sigma)$ is called 
{\it the little Teichm\"{u}ller space\index{little Teichm\"{u}ller space}  
for $\sigma$}, which is isomorphic to the product of 
lower-dimensional  Teichm\"{u}ller spaces of pointed Riemann surfaces 
complex analytically  (cf.\cite[p.289]{HK}).
Here we first review the real analytic structure of $T(\sigma)$, and then 
explain its complex analytic structure  from  the viewpoint 
of \S \ref{Subsec0301}.

Let $\Sigma_g(\sigma)=\sum_{i=1}^r R_i$ be the irreducible decomposition 
of the source stable curve and ($\hat{R}_i, \hat{\bP}_i$) ($1 \leq i \leq r$)
be the pointed Riemann surfaces obtained from the normalization of $\Sigma_g(\sigma)$ 
as in \S\ref{LiftStableFN}.
Let $g_i$ be the genus of $\hat{R}_i$ and $k_i=\sharp(\hat{\bP}_i)$ 
the number of points.
From the count of the dimensions of the vector spaces 
in the exact sequence (\ref{dualfundamentalseq}), we have
$$
\sum_{i=1}^r (3g_i-3+k_i)=3g-3-k.
$$
Now each member of the curve system 
$\tilde{\sigma} \setminus \sigma=\langle C_{k+1}, \cdots, C_{3g-3} \rangle$ 
may be considered as a simple closed curve on a unique member of 
these pointed Riemann surfaces via the normalization map. Let 
$$
\{ k+1, \cdots, 3g-3 \}=\coprod_{i=1}^r \bI_i, \:\:\: \sharp(\bI_i)=3g_i-3+k_i
$$
be the decomposition of suffixes so that 
$\{ C_{i_j} \}_{i_j \in \bI_i}$ is a system of maximal simple closed curves 
on ($\hat{R}_i, \hat{\bP}_i$)
which induces a pants decomposition of this surface.
Let $[S, w]$ be a $\sigma$-marked stable curve with the marking 
$w:\Sigma_g(\sigma) \longrightarrow S$ and let  
$S=\sum_{i=1}^r S_i$ be the irreducible decomposition, 
$(\hat{S_i}, \bP_i)  \;(1 \leq i \leq r)$ being the pointed Riemann surfaces obtained by 
the normalization. 
Then the restriction of $w$ to each component is lifted via the normalization 
to a Teichm\"{u}ller marking
\begin{equation}
\label{partialmark}
w_i: (\hat{R}_i, \hat{\bP}_i) \longrightarrow (\hat{S_i}, \bP_i).
\end{equation}
Then the Fenchel--Nielsen coordinates 
$(l_{i_j}(p), \tau_{i_j}(p)) \in (\bR_{>0})^{3g_i-3+k_i} 
\times \bR^{3g_i-3+k_i}$ are defined 
for the point $p=[(\hat{S_i}, \bP_i), w_i]$  of the 
pointed Teichm\"{u}ller space $T_{g_i, k_i}$.  They are the 
hyperbolic length $l_{i_j}(p)$ of the unique geodesic in the isotopy class of 
$w_i(C_{i_j})$ ($i_j \in \bI_i$) and the twisting parameter $\tau_{i_j}(p)$ along $C_{i_j}$. 
Then the real analytic isomorphism
$T({\sigma}) \cong \prod_{1 \leq i \leq r} T_{g_i, k_i}$ is induced via these 
Fenchel--Nielsen coordinates:
$$
\{ \ell_i([S,w]), \tau_i([S,w]) \}_{k+1 \leq i \leq 3g-3} 
\longmapsto
\prod_{1 \leq i \leq r}\{  (l_{i_j}([(\hat{S_i}, \bP_i), w_i]), \tau_{i_j}([(\hat{S_i}, \bP_i), w_i])  \}_{i_j \in \bI_i}.
$$

Over this real structure of $T({\sigma})$, the complex structure is described as 
follows.
By shrinking the base spaces of standard Kuranishi families,  
it follows from Theorem \ref{KuranishioverCont} and (\ref{locnbdbase}) 
that there exists complex coordinate patching 
$D_{\epsilon}(\sigma)=\cup_{\alpha \in \bI} B_{\alpha},$ 
where each 
$B_{\alpha}:=B_{S_\alpha}$ is the base space of a standard Kuranishi 
family for a  $\sigma$-marked stable curve $[S_\alpha, w_{\alpha}]$ 
which is embedded in 
${\rm Ext^1}_{\cO_{S_\alpha}}(\Omega_{S_\alpha}^1, \cO_{S_\alpha})$.
Then $T({\sigma}) \cap B_{\alpha}$ is isomorphic to 
${\rm H}^1(S_\alpha, {\it Hom}_{\cO_{S_\alpha}}
(\Omega_{S_\alpha}^1, \cO_{S_\alpha}))
\cap  B_{\alpha}$ 
by (\ref{fundexactseq}),(\ref{dualsp}) and the property (I)
in \S \ref{Subsec0301}.

Let $S_\alpha=\sum_{i=1}^r S_{\alpha, i}$ be the irreducible decomposition,  
and ($\hat{S}_{\alpha, i}, \bP_{\alpha, i}$)
be the pointed Riemann surface of $S_{\alpha, i}$ via the normalization.
Then the direct factor $T_{g_i, k_i}$ of $T({\sigma})$ comes from 
the direct factor 
${\rm H}^1(\hat{S}_{\alpha, i}, T_{\hat{S}_{\alpha, i}}(-\bP_{\alpha, i}))$ 
of ${\rm H}^1(S_\alpha, {\it Hom}_{\cO_{S_\alpha}}
(\Omega_{S_\alpha}^1, \cO_{S_\alpha}))$ in 
(\ref{dualsp}), i.e. 
$T_{g_i, k_i} \cap B_{\alpha}$ is isomorphic to 
${\rm H}^1(\hat{S}_{\alpha, i}, T_{\hat{S}_{\alpha, i}}(-\bP_{\alpha, i})) \cap B_{\alpha}$ 
$\cong {\rm H}^0 (\hat{S}_{\alpha, i}, 2K_{\hat{S}_{\alpha, i}}+ \bP_{\alpha, i})^* \cap B_{\alpha}$ 
by the property (I)  in \S  \ref{Subsec0301}.
These spaces are globally well-patched by the arguments 
in \cite[Chap.XV,\S2]{ACG}. Thus:

\begin{lemma} 
\label{littlepatch}
In the above notation,
the complex analytic structures of the little Teichm\"{u}ller space $T({\sigma})$ 
is described  by the coordinate patchings
$$
T({\sigma})=\bigcup_{\alpha \in \bI} 
\left( {\rm H}^1(S_\alpha, {\it Hom}_{\cO_{S_\alpha}}
(\Omega_{S_\alpha}^1, \cO_{S_\alpha}))
\cap  B_{\alpha} \right),
\:\:\:T({\sigma}) \cong \prod_{1 \leq i \leq r} T_{g_i, k_i} \:({\rm analytically}), 
\:\:\:\:\:
$$
\vspace{-0.5cm}
$$
T_{g_i, k_i}=\bigcup_{\alpha \in \bI}
\left( {\rm H}^1(\hat{S}_{\alpha, i}, T_{\hat{S}_{\alpha, i}}(-\bP_{\alpha, i})) \cap B_{\alpha} \right)
=\bigcup_{\alpha \in \bI}
\left( {\rm H}^0 (\hat{S}_{\alpha, i}, 2K_{\hat{S}_{\alpha, i}}+ \bP_{\alpha, i})^* \cap B_{\alpha} \right).
$$
\end{lemma}


\subsection{Automorphisms and cyclic branched coverings  of stable curves}
\label{Subsec0504}

In this subsection, we study automorphisms\index{automorphism} 
of stable curves and the associated 
cyclic branched coverings from stable curves to nodal Riemann surfaces.

We start from the discussion of graphs\index{graph} and their automorphisms.
A graph 
$\cG=\{ v_i, \vec{e}_j \}_{1 \leq i \leq r, 1\leq j \leq k}$ 
is a 1-dimensional  finite oriented   ``open" cell complex 
embedded in Euclidian 3-space $\bE^3$ in the following sense.
A 0-cell $v_i$ is a vertex.
A 1-cell $\vec{e}_j$ is  an oriented edge with one of the following two types.
The first type is an {\it ordinary edge} 
$\vec{e}_j=\vec{e}_j(v_{h_1(j)}, v_{h_2(j)})$ 
which connects a vertex $v_{h_1(j)}$ to a vertex $v_{h_2(j)}$ in this direction.
(Note that  $v_{h_1(j)}=v_{h_2(j)}$ may occur. This case corresponds to a self-intersection of an irreducible component.)
The second type is an {\it open edge} 
which emanates from a vertex $v_{h_1(j)}$ such that the end point is a point 
in $\bE^3$ which is not contained 
in the set of vertices, and we  symbolically write it as $\vec{e}_j=\vec{e}_j(v_{h_1(j)}, *)$. 
The {\it absolute edge} $|\vec{e}_j|$ is the usual edge obtained 
by neglecting the orientation of $\vec{e}_j$.
If all the oriented edges of $\cG$ are ordinary, we call $\cG$ a 
{\it compact graph}\index{compact graph}.

We define the contraction map of $\cG$ 
$$
{\rm cont}: \cG \longrightarrow \cG^{c}
$$
as the identity map on the complement of  open  edges 
such that any open edge $\vec{e}_j(v_{h_1(j)}, *)$ is contracted to  
the vertex $v_{h_1(j)}$.
Then $\cG^{c}$ is a compact graph.

An automorphism $\bar{\varphi}: \cG \rightarrow \cG$ 
of a compact graph $\cG$ means that 
$\bar{\varphi}$ is a homeomorphism in the Euclidian topology 
which preserves the set of vertices and 
the set of absolute edges
such that $\bar{\varphi}^n$ is the identity map for some $n$.
The least natural number $n$ which enjoys the above property is called the {\it order} of $\bar{\varphi}$.

For a vertex $v_i$ (resp.~an edge $\vec{e}_j$), 
there exists a minimal natural number 
$m(v_i)$ (resp.~$m(\vec{e}_j)$) which satisfies 
$\bar{\varphi}^{m(v_i)}(v_i)=v_i$ 
(resp.~$\bar{\varphi}^{m(\vec{e}_j)}(\vec{e}_j)=\vec{e}_j$). 
If $m(\vec{e}_j)$ is even and 
$\bar{\varphi}^{m(\vec{e}_j)/2}$ stabilizes 
$|\vec{e}_j|$ by reversing its orientation, i.e. 
$\bar{\varphi}^{m(\vec{e}_j)/2}(\vec{e}_j)=-\vec{e}_j$, 
then $\vec{e}_j$ is said to be an {\it amphidrome edge} for $\bar{\varphi}$.
Otherwise, $\vec{e}_j$ is said to be a {\it non-amphidrome edge}. 
Let ${\bf V}, \:{\bf NE}, \:{\bf AE}$ be the set of vertices, 
non-amphidrome edges
and amphidrome edges
of $\cG$ for $\bar{\varphi}$, respectively. Let 
$$
{\bf V}=\coprod_{1 \leq i \leq \overline{r}} \:\coprod_{0 \leq \ell \leq m(v_i)-1} \bar{\varphi}^{\ell}(v_i),
\:\:\:\:\:\:\:\:\:\:\:\:\:\:\:\:\:\:\:\:\:\:\:\:\:\:\:\:\:\:\:\:\:\:\:\:\:\:\:\:\:\:\:\:\:\:\:\:\:\:\:\:\:\:\:\:\:\:\:\:\:\:\:\:\:\:\:\:\:\:\:\:\:\:\:\:\:\:\:\:\:\:\:
$$
\vspace{-0.4cm}
$$
{\bf NE}=\coprod_{1 \leq j \leq \overline{k}_1} \:\coprod_{0 \leq \ell \leq m(\vec{e}_j)-1} \bar{\varphi}^{\ell}(\vec{e}_j),\:\:\:
{\bf AE}=\coprod_{ \overline{k}_1+1 \leq j \leq \overline{k}_1+\overline{k}_2} \:
\coprod_{0 \leq \ell \leq m(\vec{e}_j)/2-1} \bar{\varphi}^{\ell}(\vec{e}_j)
$$
be the orbit decompositions for $\bar{\varphi}$.
Here  $\{ v_1, \cdots, v_{ \overline{r}} \}$, 
$\{ \vec{e}_1, \cdots, \vec{e}_{\overline{k}_1} \}$ and 
$\{ \vec{e}_{\overline{k}_1+1}, \cdots, \vec{e}_{\overline{k}_1+\overline{k}_2} \}$
are assumed to belong to mutually distinct orbits in 
${\bf V}$, ${\bf NE}$ and ${\bf AE}$ respectively.
We have  $r=\sum_{1 \leq i \leq \overline{r}} m(v_i)$, 
$k=\sum_{1 \leq j \leq \overline{k}_1} m(\vec{e}_j) 
+\sum_{ \overline{k}_1+1 \leq j \leq \overline{k}_1+\overline{k}_2} m(\vec{e}_j)/2$.
 
Let $G_{\bar{\varphi}} \cong \bZ/n\bZ$ be the cyclic group 
generated by $\bar{\varphi}$.

\begin{definition}
\label{Def0504no1new}
{\rm (i)}
The quotient map  and the  quotient graph\index{quotient graph} of $\cG$ by 
$G_{\bar{\varphi}}$
$$
\pi_{\bar{\varphi}}:  \cG \longrightarrow \cW
:=\cG/G_{\bar{\varphi}}
$$
are defined by the following two conditions;

\noindent {\rm (ia)} 
The vertices and the edges are given by 
$\cW=\{ v_i^{\sharp}, (\vec{e}_j)^{\sharp} \}
_{1 \leq i \leq \overline{r}, 1\leq j \leq \overline{k}_1+\overline{k}_2}$ 
so that 
$\pi_{\bar{\varphi}}(\bar{\varphi}^{\ell}(v_i))$ $=v_i^{\sharp}$ 
and 
$\pi_{\bar{\varphi}}(\bar{\varphi}^{\ell}(\vec{e}_j))
=(\vec{e}_j)^{\sharp}$ for any $i, j, \ell$. 

\noindent {\rm (ib)} 
Suppose $\vec{e}_j=\vec{e}_j(v_{h_1(j)}, v_{h_2(j)})$. 
If $1 \leq j \leq \overline{k}_1$, then $(\vec{e}_j)^{\sharp}$ 
is an ordinary edge given by 
$(\vec{e}_j)^{\sharp}
(\pi_{\bar{\varphi}}(v_{h_1(j)}), \pi_{\bar{\varphi}}(v_{h_2(j)}))$.
If $\overline{k}_1+1 \leq j \leq \overline{k}_1+\overline{k}_2$, 
then $(\vec{e}_j)^{\sharp}$ 
is an open  edge given by
\begin{equation}
\label{Eq0504no3new}
(\vec{e}_j)^{\sharp}
(\pi_{\bar{\varphi}}(v_{h_1(j)}),*).
\end{equation}

\noindent {\rm (ii)}
The compact quotient map and the  compact quotient graph of $\cG$ by 
$G_{\bar{\varphi}}$ are defined by 
$$
\pi_{\bar{\varphi}}^c={\rm cont} \circ \pi_{\bar{\varphi}}: \cG \longrightarrow 
(\cW)^c= (\cG/G_{\bar{\varphi}})^c.
$$
\end{definition}

With respect to (\ref{Eq0504no3new}), 
we may write $(\vec{e}_j)^{\sharp}
(\pi_{\bar{\varphi}}^{\prime}(v_{h_2(j)}),*)$ 
because $\vec{e}_j$ is amphidrome in this case and 
$v_{h_1(j)}$ and $v_{h_2(j)}$ are on the same orbit.

\begin{example}
\label{Eq0504no1new}
We consider the compact graph 
$\cG=\{ v_i, \vec{e}_j(v_1, v_2) \}_{1 \leq i, j \leq 2}$.
Let $\bar{\varphi}, \bar{\psi}: \cG \rightarrow \cG$ be automorphisms 
of order $2$ defined as follows.
The first $\bar{\varphi}$ interchanges $v_1$ and $v_2$, and stablizes 
$|\vec{e}_j|$ by reversing their orientations.
By using the usual symbols of permutations, $\bar{\varphi}$ is written by 
$(v_1, v_2), (\vec{e}_1, -\vec{e}_1), (\vec{e}_2, -\vec{e}_2)$.
The second $\bar{\psi}$ interchanges $v_1$ and $v_2$, and also $\vec{e}_1$ 
and $\vec{e}_2$, i.e. $\bar{\psi}$ is written by 
$(v_1, v_2), (\vec{e}_1, -\vec{e}_2)$.

Then $\vec{e}_1$ and $\vec{e}_2$ are amphidrome (resp. non-amphidrome) 
for $\bar{\varphi}$ (resp. for $\bar{\psi}$), and 
we have
$\cG/G_{\bar{\varphi}}
=\{ v_1^{\sharp}, \:(\vec{e}_1)^{\sharp}(v_1^{\sharp}, *)$,  
$(\vec{e}_2)^{\sharp}(v_1^{\sharp}, *) \}$, 
$(\cG/G_{\bar{\varphi}})^c=\{ v_1^{\sharp}\}$ (empty edge), 
$\cG/G_{\bar{\psi}}=(\cG/G_{\bar{\psi}})^c
=\{ v_1^{\sharp}, \:(\vec{e}_1)^{\sharp}(v_1^{\sharp}, v_1^{\sharp})\}$
as shown in Figure II.
\end{example}

\begin{figure}[h]

\setlength{\unitlength}{0.65mm}
\begin{picture}(200,30)(-25,0)

\put(5, 20){\circle*{4}}
\put(25, 20){\circle*{4}}
\qbezier(6.5,21)(14,28)(23.5,21)
\put(20.5,21){\line(1,0){3}}
\put(23.5,21){\line(-1,2){1.5}}
\put(4,26){\makebox(2,2){$v_1$}}
\put(24,26){\makebox(2,2){$v_2$}}
\put(14,28){\makebox(2,2){$\vec{e}_1$}}
\qbezier(6.5,19)(14,12)(23.5,19)
\put(20.5,19){\line(1,0){3}}
\put(23.5,19){\line(-1,-2){1.5}}
\put(14,10){\makebox(2,2){$\vec{e}_2$}}
\put(14,2){\makebox(2,2){$\cG$}}

\put(60, 20){\circle*{4}}
\put(61.5,21){\line(2,1){8}}
\put(69,25){\line(-1,0){2}}
\put(69,25){\line(-1,-2){1}}
\put(61.5,19){\line(2,-1){8}}
\put(69,15){\line(-1,0){2}}
\put(69,15){\line(-1,2){1}}
\put(59,26){\makebox(3,3){$v_1^{\sharp}$}}
\put(74,27){\makebox(3,3){$(\vec{e}_1)^{\sharp}$}}
\put(74,11){\makebox(3,3){$(\vec{e}_2)^{\sharp}$}}
\put(69,0){\makebox(2,5){$\cG/G_{\bar{\varphi}}$}}

\put(115, 20){\circle*{4}}
\put(115,0){\makebox(2,5){$(\cG/G_{\bar{\varphi}})^c$}}
\put(114.5,27){\makebox(3,3){$v_1^{\sharp}$}}

\put(170, 20){\circle*{4}}
\cbezier(172,21)(176,28)(181,28)(182,20)
\cbezier(172,19)(176,12)(181,12)(182,20)
\put(172,19){\line(1,0){2}}
\put(172,19){\line(0,-1){2}}
\put(177,-1){\makebox(2,8){$\cG/G_{\bar{\psi}}=(\cG/G_{\bar{\psi}})^c$}}
\put(169,27){\makebox(3,3){$v_1^{\sharp}$}}
\put(187,19){\makebox(3,3){$(\vec{e}_1)^{\sharp}$}}

\end{picture}

\begin{center}
(Figure II) The quotient and the compact quotient graphs in Example \ref{Eq0504no1new}
\end{center}

\end{figure}

Now let $S$ be a stable curve of genus $g \geq 2$ with its irreducible decomposition 
$S=\sum_{i=1}^{r}S_i$,  
and let $\bP=\{ P_1, \cdots, P_k \}$ be its set of  nodes.
By this configuration, a compact graph 
${\rm rg}(S)=\{ v_{S_i}, \vec{e}_{P_j} \}_{1 \leq i \leq r, 1 \leq j \leq k}$ 
is defined as follows:
An irreducible component $S_i$ corresponds to a vertex $v_{S_i}$.
If a node $P_j$ is an intersection point of the irreducible components 
$S_{h_1(j)}$ and $S_{h_2(j)}$, 
then it is expressed by an ordinary edge 
$\vec{e}_{P_j}=\vec{e}_{P_j}(v_{S_{h_1(j)}},v_{S_{h_2(j)}})$ (the orientation is arbitrary).
We call ${\rm rg}(S)$ 
{\it the reduced dual graph}\index{reduced dual graph} of $S$.

Let $\varphi: S \rightarrow S$ be  an analytic automorphism 
of order $N$, and $G=\langle \varphi \rangle \cong \bZ/N\bZ$ the subgroup of 
${\rm Aut}(S)$ generated by $\varphi$.
Then $\varphi$ clearly induces the automorphism 
$\varphi_{{\rm rg}(S)}: {\rm rg}(S) \rightarrow {\rm rg}(S)$.
We translate the terminologies defined for 
(${\rm rg}(S), \varphi_{{\rm rg}(S)}$) into 
($S, \varphi$).
That is to say, $m(S_i)$ and $m(P_j)$ are the minimal natural numbers  
which satisfy 
$\varphi^{m(S_i)}(S_i)=S_i$ and 
$\varphi^{m(P_j)}(P_j)=P_j$, 
and $P_j$ is an {\it amphidrome node}\index{amphidrome}
 (resp.~{\it a non-amphidrome node}\index{non-amphidrome}) 
for $\varphi$ 
if $\vec{e}_{P_j}$ is an amphidrome edge (resp.~a non-amphidrome edge) 
for $\varphi_{{\rm rg}(S)}$.
{\it The orbit-irreducible decomposition} of $S$ for $\varphi$ is defined by
\begin{equation}
\label{Eq0504no1}
S=\sum_{i=1}^{\overline{r}}\sum_{j=0}^{m(S_i)-1} \varphi^{j}(S_i), 
\end{equation}
where each $S_i$ for $1 \leq i \leq \overline{r}$ has mutually distinct orbits and 
$r=\sum_{i=1}^{\overline{r}}m(S_i)$.

Let 
$h:\whS=\coprod_{i=1}^r \coprod_{j=1}^{m(S_i)-1} \widehat{\varphi^{j}(S_i)} 
\longrightarrow S$ be the normalization map, and 
$\varphi_{i,j}: \widehat{\varphi^{j}(S_i)} 
\longrightarrow \widehat{\varphi^{j}(S_i)}$
the lifting of  
$\varphi^{m(S_i)}|_{\varphi^{j}(S_i)}: {\varphi^{j}(S_i)} 
\longrightarrow {\varphi^{j}(S_i)}$.
Since $\widehat{\varphi^{j}(S_i)} \cong \whS_i$ and 
$\varphi_{i,j}$ is congruent to $\varphi_{i,0}$ for  
$0 \leq j \leq m(S_i)-1$, 
we may consider 
\begin{equation}
\label{Eq0504no2}
\varphi_i:=\varphi_{i,0}: \whS_i \longrightarrow \whS_i \:\:\:\:\:\:\: 
(1 \leq i \leq \overline{r})
\end{equation} 
as the representatives of the $\varphi_{i,j}$'s.
Let $n_i$ be the order of $\varphi_i$.
We have an $n_i$-fold cyclic branched covering of Riemann surfaces
\begin{equation}
\label{smallcover}
\pi_{\varphi_i}: \whS_i \longrightarrow \whW_i \cong \whS_i/G_{\varphi_i},  
\:\:\:\:G_{\varphi_i}=\langle \varphi_i \rangle \cong \bZ/n_i\bZ.
\end{equation}
Note that  $\varphi_i$ also defines an automorphism 
of the pointed Riemann surface  
\begin{equation}
\label{Eq0504no3}
\varphi_i: (\whS_i, \bP_i) \longrightarrow (\whS_i, \bP_i).
\end{equation}
We divide the set $\bP_i$ and the set $\bQ_i=\pi_{\varphi_i}(\bP_i)$ on $\whW_i$ into
\begin{equation}
\label{Eq0504no4b}
\bP_i=\bP_i^{\bN} \coprod \bP_i^{\bA}, \:\:\:\:\: 
\bQ_i=\bQ_i^{\bN} \coprod \bQ_i^{\bA} 
\end{equation}
where $\bP_i^{\bN}$ (resp.~$\bP_i^{\bA}$) consists of the points $P$ such that 
the nodes $h(P)$ in $S$ are non-amphidrome (resp.~amphidrome) for $\varphi$, 
and $\bQ_i^{\bN}=\pi_{\varphi_i}(\bP_i^{\bN})$, $\bQ_i^{\bA}=\pi_{\varphi_i}(\bP_i^{\bA})$.

The  {\it dual graph}\index{dual graph}  
${\rm dg}(S)=\{ (v_{S_i},g(\whS_i )), \vec{e}_{P_j} \}_{1 \leq i \leq r, 1 \leq j \leq k}$ 
is the weighted graph 
obtained from the reduced dual graph ${\rm rg}(S)$ by attaching the weight 
$g(\whS_i)$ to each vertex $v_{S_i}$, 
which is the genus of $\whS_i$.
(If $g(\whS_i )=0$, it is sometimes omitted.)
An automorphism $\varphi:S \rightarrow S$ also induces an automorphism
$\varphi_{{\rm dg}(S)}: {\rm dg}(S) \rightarrow {\rm dg}(S)$, since 
$\varphi$ preserves $g(\whS_i)$.

The following two lemmata guarantee the 
existence of a natural quotient of $S$ by $G$ as a nodal Riemann surface.
\begin{lemma}
\label{Lem0504no1}
There exists a finite holomorphic map 
\begin{equation}
\label{quotnodal}
\pi_{\varphi}:S \longrightarrow W
\end{equation} 
to a nodal Riemann surface $W$ which has the following properties:

\noindent {\rm (i)} 
$W$ has an irreducible decomposition 
$\sum_{i=1}^{\overline{r}} W_i$ such that
$W_i=\pi_{\varphi}(\varphi^{j}(S_i))$ for any $j$.

\noindent {\rm (ii)} 
The normalizations of $S$ and $W$ naturally induce 
cyclic branched coverings
$$ 
\pi_{\varphi_{i,j}}: \widehat{\varphi^j(S_i)} \longrightarrow 
\widehat{W_i} \:\:\:\:\:\:\:
(1 \leq i \leq \overline{r}, \: 0 \leq j \leq m(S_i)-1)
$$ 
such that $\pi_{\varphi_{i,j}}$ is isomorphic to $\pi_{\varphi_i}$ 
in (\ref{smallcover}) for any $j$.

\noindent {\rm (iii)} The reduced dual graph ${\rm rg}(W)$ coincides with 
the compact quotient graph 
$({\rm rg}(S)/G_{\varphi_{{\rm rg}(S)}})^c$  
so that the dual graph is written as ${\rm dg}(W)=
\{ (v_{W_i}, g(\whW_i)), \vec{e}_{Q_j} \}
_{1 \leq i  \leq \overline{r}, 1 \leq j \leq \overline{k}_1}$.
In particular, a non-amphidrome node in $S$ is sent by $\pi_{\varphi}$ 
to a node of $W$, while an amphidrome node is sent to a non-singular point. 
\end{lemma}

{\it Proof} \quad
By using $\{ (\whW_i, \bQ_i^{{\bf N}}) \}_{1 \leq i \leq \overline{r}}$ in  
(\ref{smallcover}),  (\ref{Eq0504no4b})
as the building blocks for patchings, one can construct the desired 
nodal surface $W$ so that 
$\bQ_i^{{\bf N}}$ are patched as the nodes of $W$.
This patching  process is constructed 
by identifying  the 
two local components of $xy=0$ at the origin of $\bC^2$ 
with the disk neighborhoods of $\whW_i$'s at the suitable points in $\bQ_i^{{\bf N}}$.
\qed

\begin{lemma}
\label{Lem0504no2}
$W$ (in Lemma \ref{Lem0504no1}) is isomorphic to the quotient analytic space $S/G$.
\end{lemma}

{\it Proof} \quad 
We consider the local ring $\cO_{S,P}$ at the node $P$ of $S$, and 
let $\widehat{\cO}_{S,P} \cong \bC[[x, y]]/(xy)$ be its completion by 
the maximal ideal of $\cO_{S,P}$.
If $P$ is non-amphidrome whose 
covalencies at both banks are $\delta^{(1)}/\lambda^{(1)}$ and 
$\delta^{(2)}/\lambda^{(2)}$, then the action is written  by 
\begin{equation}
\label{Eq0504no2b}
\varphi^{m(P)}: (x,y) \longmapsto
(\be(\delta^{(1)}/\lambda^{(1)})x, (\be(\delta^{(2)}/\lambda^{(2)})y).
\end{equation}
We may assume $\lambda^{(1)} \geq \lambda^{(2)}$.
Then the invariant subring of $\widehat{\cO}_{S,P}$ for 
$\varphi^{m(P)}$ is given by
$$
(\widehat{\cO}_{S,P})^{\varphi^{m(P)}} 
\cong \bC[[z, w]]/(zw),  
$$
where $z=x^{\lambda^{(1)}}$, $w=y^{\lambda^{(2)}}$ and 
$zw=x^{\lambda^{(1)}-\lambda^{(2)}}(xy)^{\lambda^{(2)}}=0$.
If $P$ is amphidrome for $\varphi^{m(P)}$ whose 
covalency is $\delta/\lambda$, 
then 
\begin{equation}
\label{Eq0504no2c}
\varphi^{m(P)/2}: (x,y) \longmapsto
(\be(\delta/2\lambda)y, (\be(\delta/2\lambda)x).
\end{equation}
Hence
$$
(\widehat{\cO}_{S,P})^{\varphi^{m(P)}} 
\cong \bC[[z]], 
$$
where $z=x^{m(P)}=y^{m(P)}$. 
If $P$ is a nonsingular point of $S$,
the similar invariant subring is obviously regular.

It follows that $S/G$ exists as  a complex curve 
with at most nodes.
The natural morphism $S \rightarrow S/G$ has the same 
properties as those in Lemma \ref{Lem0504no1}.
Therefore,  $W$ is isomorphic to $S/G$.
\qed

\begin{definition}
\label{Def0504no1}
We write (\ref{quotnodal}) as $\pi_{\varphi}:S \longrightarrow W=S/G$ 
and call it the cyclic branched covering\index{cyclic branched covering}
 associated with $\varphi$.
\end{definition}

\begin{example}
\label{EX0504no1}
Let $S$ be a stable curve of genus $3$ with two irreducible components and two nodes; 
$S=\sum_{i=1}^2 S_i, \:S_1 \cap S_2=\{ P_1, P_2 \}$.
Assume that each  $S_i$ is isomorphic to an elliptic curve with the period $\sqrt{-1}$; 
$S_i \simeq \bC/(\bZ+\sqrt{-1}\bZ)$.
We consider the following two automorphisms $\varphi, \psi:S \rightarrow S$ 
of order 8 so that  the automorphisms $\varphi_{{\rm rg}(S)}$ and  $\psi_{{\rm rg}(S)}$ 
are given in Example \ref{Eq0504no1new}.
First, $\varphi$ fixes $P_1$ and $P_2$ and interchanges $S_1$ and $S_2$ such that 
the total valencies of $\varphi^2|_{S_i}$ $(i=1,2)$ 
are ${\bf 3/4}+{\bf 3/4}+1/2$ (where ${\bf 3/4}$ are attached 
to $P_1$ and $P_2$).
Secondly, 
$\psi$ interchanges $P_1$ and $P_2$ and interchanges $S_1$ and $S_2$ such that 
the total valencies of $\psi^2|_{S_i}$ are ${\bf 3/4}+{\bf 3/4}+1/2$.

Then $P_1$ and $P_2$ are amphidrome nodes (resp. non-amphidrome nodes) 
for $\varphi$ (resp.  $\psi$).
Let $\pi_{\varphi}:S \rightarrow W_{\varphi}$, 
$\pi_{\psi}:S \rightarrow W_{\psi}$ be the cyclic branched covering associated with 
$\varphi$ and $\psi$.
Then $W_{\varphi}$ is isomorphic to $\bP^1$, while $W_{\psi}$ is a stable curve of genus 1 
with one node.
Their reduced dual graphs  are given in Example \ref{Eq0504no1new}.
\end{example}

\begin{remark}
\label{nonreducedquot}
It may be natural to re-write $\pi_{\varphi}: S \longrightarrow W$ as 
$\pi_{\varphi}:S \longrightarrow \overline{W}=\sum_{i=1}^{\overline{r}} m_iW_i$
where $\overline{W}$ is a non-reduced scheme such that 
$\overline{W}^{{\rm red}}=W$.
This point will be also discussed in \S \ref{ppsubsec}
(see also \cite[Chap.3]{MM}). 
\end{remark}

Let $\bB_i \subset \whW_i$ be the set of branch points of 
$\pi_{\varphi_i}:\whS_i \rightarrow \whW_i$
in (\ref{smallcover}).
By comparing with (\ref{Eq0504no4b}), 
we define the cardinalities of the sets by 
\vspace{-0.1cm}
$$
\sharp(\bB_i)=s_i, \;\;
\sharp(\bB_i \cap \bQ_i)=s_i^{\prime},\;\;
\sharp(\bB_i \cap \bQ_i^{\bN})=(s_i^{\prime})^{(1)},
\sharp(\bB_i \cap \bQ_i^{\bA})=(s_i^{\prime})^{(2)}, 
\:s_i^{\prime}=(s_i^{\prime})^{(1)}+(s_i^{\prime})^{(2)},
$$
\begin{equation}
\label{Eq0504no4c}
\sharp(\bQ_i \setminus \bB_i)=t_i,\;\;
\sharp(\bQ_i^{\bN} \setminus \bB_i)=t_i^{(1)},\:\:
\sharp(\bQ_i^{\bA} \setminus \bB_i)=t_i^{(2)}, \:\:
\:t_i=t_i^{(1)}+t_i^{(2)}.
\:\:\:\:\:\:\:\:\:\:\:
\end{equation}
We add the valency $1$ for each non-branch point in 
$\bQ_i \setminus \bB_i$ to the total valency (\ref{totalvalency}) 
for $\varphi_i$, and define the 
{\it dressed total valency}\index{dressed total valency} for  
$\varphi_i$ by
$$
{\rm DTV}(\varphi_i)= 
\left( 
g_i, \overline{g}_i, n_i: {\bf \frac{\sigma_1^{(i)}}{\lambda_1^{(i)}}} 
+ \cdots +{\bf \frac{\sigma_{(s'_i)^{(1)}}^{(i)}}{\lambda_{(s'_i)^{(1)}}^{(i)}}}
+\bigl( \bigl( \frac{\sigma_{(s'_i)^{(1)}+1}^{(i)}}
{\lambda_{(s'_i)^{(1)}+1}^{(i)}} \bigr) \bigr) + \cdots 
+\bigl( \bigl( \frac{\sigma_{s'_i}^{(i)}}{\lambda_{s'_i}^{(i)}} \bigr) \bigr)
\right.
\:\:\:\:\:\:\:\:
$$
\vspace{-0.3cm}
\begin{equation}
\label{Eq0504no5b}
\left.
+\frac{\sigma_{(s'_i)+1}^{(i)}}{\lambda_{s'_i+1}^{(i)}}  + \cdots 
\frac{\sigma_{s_i}^{(i)}}{\lambda_{s_i}^{(i)}} 
+\underbrace{{\bf 1}+ \cdots+{\bf 1}}_{t_i^{(1)}}
+\underbrace{((\bf 1))+ \cdots+((1))}_{t_i^{(2)}}
\right)
\end{equation}
where the valencies are ordered for 
$\bB_i \cap \bQ_i^{\bN}$, 
$\bB_i \cap \bQ_i^{\bA}$,
$\bB_i \setminus \bQ_i$, 
$\bQ_i^{\bN} \setminus \bB_i$ and 
$\bQ_i^{\bA} \setminus \bB_i$.
These symbols (bold faced valency for a non-amphidrome 
branch point, 
soft double bracket valency for an amphidrome branch point, etc.)
are borrowed from \cite[\S2]{AI}.
Then:

\begin{definition}
\label{numdetaauto}
We define the 
{\rm numerical data}\index{numerical data} of an automorphism 
$\varphi:S \longrightarrow S$ of a stable curve (of order $N$) by  
\begin{equation}
\label{defnumdata}
{\rm Num}(\varphi)=
\left(N, \varphi_{{\rm dg}(S)}, \coprod_{i=1}^{\overline{r}}{\rm DTV}(\varphi_i) \right).
\end{equation}
\end{definition}

\begin{remark}
\label{exEtype}
These numerical data are essentially the same as 
the numerical data of pseudo-periodic maps of negative twist 
(see \S \ref{ppsubsec}).
From this viewpoint,  
${\rm Num}(\varphi)$  for $g=3$ are  
classified  in 
\cite[Table 2, Lemma 3.4, Prop.~3.8]{AI}.
\end {remark}


\subsection{Equisymmetric strata at the boundary charts of $\overline{M}_g^{orb}$}
\label{Subsec0505}

We prove a  structure theorem  for equisymmetric strata\index{equisymmetric strata} 
consisting of marked stable curves 
which have the same numerical invariants of automorphisms.

Let $w:\Sigma_g(\sigma) \longrightarrow S$  be a $\sigma$-marking 
such that $S$ has an automorphism $\varphi$ with the preassigned numerical data ${\rm Num}(\varphi)$ 
 as in (\ref{defnumdata}).
From the orbit-irreducible decomposition 
$S=\sum_{i=1}^{\overline{r}}\sum_{j=0}^{m_i-1} \varphi^{j}(S_i)$ ($m_i= m(S_i)$) 
as in (\ref{Eq0504no1}), 
the little Teichm\"{u}ller space $T({\sigma})$ which contains  the point 
$[S, w]$ is isomorphic to 
\begin{equation}
\label{speciallittleTeih}
T({\sigma}) \cong \prod_{i=1}^{\overline{r}} 
\underbrace{\left(T_{g_i, k_i} \times \cdots \times T_{g_i, k_i}\right)}_{m_i},
\end{equation}
where $g_i=g(\whS_i), k_i=\sharp(\bP_i)$. See Lemma \ref{littlepatch}.

In \S5.1, we defined, with the explanation (I), the subspaces $T_g^{\varphi}$ and $T_g^{[\varphi]}$ for a 
periodic map $\varphi: S \to S$. There $S$ was a Riemann suface. We consider here a similar definition 
for a stable curve $S$.\par
As in the previous subsection \S5.4, let $\varphi: S \to S$ be an analytic automorphism of order $N$ of a stable curve $S$.
 Let $T_{\sigma}^{\varphi}$ be the set of points $p=[S, w]$ in $T(\sigma)$ which is fixed 
by the action of $\varphi_{*}: T(\sigma) \to T(\sigma)$. (Through the $\sigma$-marking $w: \Sigma_g(\sigma) \to S$, the 
analytic automorphism $\varphi$ is considered to be  an element $w^{-1}\circ\varphi \circ w$ of the Weyl group $W(\sigma)$. The action $\varphi_*$ 
is the action as an element of $W(\sigma)$. See (\ref{changeWmark}).) As in the case of Riemann surfaces, 
we define the equivalence class $[\varphi]$ to be the set of those analytic automorphisms $\psi: S \to S$ 
which have the same numerical data (\ref{defnumdata}) as $\varphi$ : ${\rm Num}(\psi)={\rm Num}(\varphi)$. 
Then $T_\sigma^{[\varphi]}$ is defined as follows:
$$T_\sigma^{[\varphi]} = \bigcup_{\psi\in [\varphi]}T_\sigma^{\psi}.$$
The subspace $\overline{M}_g^{[\varphi]}$ of $\overline{M}_g$ is defined 
to be the quotient space of $T_\sigma^{[\varphi]}$ by the Weyl group:  $\overline{M}_g^{[\varphi]} =T_\sigma^{[\varphi]}/W(\sigma)$.

\begin{remark}
\label{Remnumconj}
As we will prove in \S6.1, an analytic automorphism $\varphi$ of a stable  
curve $S$ 
is lifted to a pseudo-periodic map of negative twist 
$\tilde{\varphi}$ of a Riemann surface 
such that ${\rm Num}(\varphi)$ may be identified 
with the invariants {\rm (a)}, {\rm (c)} of $\tilde{\varphi}$ in Th.~6.1.
Hence ${\rm Num}(\varphi)$ determines the  conjugacy class of 
the 
analytic automorphism $\varphi$ of the stable curve $S$ by the results of 
\cite{Nielsen2}, \cite{MM}.
\end{remark}

\begin{theorem}
\label{stableequithm} 
{\rm (i)} There exist  analytic embeddings of 
pointed Teichm\"{u}ller spaces\index{pointed Teichm\"{u}ller space}
$$
\psi_{i, j}: T_{\overline{g}_i, s_i+t_i} \longrightarrow  T_{g_i, k_i}\:\:
(1 \leq i \leq \overline{r}, \:0 \leq j \leq m_i-1), 
$$
where $t_i$ is given in (\ref{Eq0504no4c}), and also an analytic  embedding 
\begin{equation}
\label{embPhi}
\Phi:\prod_{i=1}^{\overline{r}}T_{\overline{g}_i, s_i+t_i} 
\hookrightarrow 
\prod_{i=1}^{\overline{r}} 
\underbrace{\left(T_{g_i, k_i} \times \cdots \times T_{g_i, k_i}\right)}_{m_i} 
\cong T({\sigma}),
\end{equation}
\vspace{-0.4cm}
\begin{equation}
\label{embmapexplicit}
\Phi:(x_1, \cdots, x_{\overline{r}}) \longmapsto 
(\underbrace{\psi_{1,0}(x_1), \cdots, \psi_{1, m_i-1}(x_1)}_{m_1},\cdots, 
\underbrace{\psi_{\overline{r},0}(x_{\overline{r}}), \cdots, 
\psi_{\overline{r},m_{\overline{r}-1}}(x_{\overline{r}})}_{m_{\overline{r}}}),
\end{equation}
such that the connected component $(T_{\sigma}^{[\varphi]})^{(0)}$ 
of $T_{\sigma}^{[\varphi]}$ 
containing the point $[S, w]$ is analytically isomorphic to 
the image 
$\Phi(\prod_{i=1}^{\overline{r}}T_{\overline{g}_i, s_i+t_i})$. 
In particular,  $(T_{\sigma}^{[\varphi]})^{(0)}$ is a 
$\sum_{i=1}^{\overline{r}}(3\overline{g}_i-3+s_i+t_i)$-dimensional 
complex submanifold of $T({\sigma})$.
The space $T_{\sigma}^{[\varphi]}$ itself is a countable union of 
submanifolds\index{countable union of 
submanifolds}
 of these types of  connected components. 

\noindent 
{\rm (ii)} $\overline{M}_g^{[\varphi]}$ is a locally closed  irreducible subvariety 
of $\overline{M}_g$ which is isomorphic to 
$\prod_{i=1}^{\overline{r}}M_{\overline{g}_i, s_i+t_i}$.
\end{theorem}

{\it Proof} \quad {\it Step 1} \quad 
Let $S$ be a stable curve with the irreducible decomposition
$S=\sum_{i=1}^r S_i$ such that $S$ has 
an automorphism $\varphi$ with the preassigned 
numerical data  (\ref{defnumdata}).
Let $\psi:X \rightarrow B$ be a standard Kuransihi family of 
$S= \psi^{-1}(b_0)$ ($b_0 \in B$).
We may assume that $B$ is a small open ball around the origin ($=b_0$) of the vector space 
${\rm Ext}_{\cO_S}^1(\Omega_S^1, \cO_S)$, which we write 
$B={\rm Ext}_{\cO_S}^1(\Omega_S^1, \cO_S) \cap B$ if we want to emphasize this 
ambient space.
From (\ref{dualsp}) and (I) of \S \ref{Subsec0301},
we restrict $B$ to the subspace which parematrizes the variable deformation
\begin{equation}
\label{defnonsmoothhere}
{\rm H}^1(S, {\it Hom}_{\cO_{S}}(\Omega_{S}^1, \cO_{S})) \cap B
\cong 
\bigoplus_{i=1}^{r} V_i, \:\:{\rm where}\:\:
V_i:={\rm H}^0 (\hat{S}_i, 2K_{\hat{S}_i}+ \bP_i)^* \cap B.
\end{equation}
Now $\varphi$ relatively acts on $\psi$, and let 
$B^{\varphi}={\rm Ext}_{\cO_S}^1(\Omega_S^1, \cO_S)^{\varphi} \cap B$ 
be the invariant subspace.
Since the action of $\varphi$ on $S$ preserves the set of nodes, 
the action of  $\varphi$ on $B$ preserves the subspace 
$\bigoplus_{i=1}^{r} V_i$.
Set $m_i=m(S_i)$.
Since the iterations of the action of $\varphi$ on the direct factor $V_i$ 
isomorphically map  
$V_i \rightarrow \varphi(V_i) \rightarrow \varphi^2(V_i) \rightarrow \cdots$
and stabilize $\varphi^{m_i}(V_i)=V_i$, 
the $\varphi$-invariant subspace 
$(\bigoplus_{i=1}^r V_i)^{\varphi}$ 
is isomorphic to the direct sum of $\varphi^{m_i}$-invariant subspaces 
\begin{equation}
\label{directsuminvsubsp}
\left( \bigoplus_{i=1}^r V_i \right)^{\varphi}  \cong
\bigoplus_{i=1}^{\overline{r}} 
\bigoplus_{j=0}^{m_i-1} \left( \varphi^j(V_i) \right)^{\varphi_i}, 
\:\: {\rm where} \:\:\varphi_i:=\varphi^{m_i}.
\end{equation}
The direct factor $V_i^{\varphi_i}$ of (\ref{directsuminvsubsp}) 
is described as follows.
We consider the covering 
$\pi_{\varphi_i}: \whS_i \longrightarrow \whW_i$ 
in (\ref{smallcover}), and set 
$\bB_i=\{ Q_1, \cdots, Q_{s_i} \}$,
$\bQ_i \setminus \bB_i=\{ Q_{s_i+1}, \cdots, Q_{s_i+t_i} \}$ 
from (\ref{Eq0504no4c}).
By appling Proposition \ref{eigenspaceformula} to the $0$-eigenspace, 
the dimension of $V_i^{\varphi_i}$ is equal to 
$3\overline{g}_i-3+s_i+t_i$.
More precisely, the discussion in the proof of 
Proposition \ref{eigenspaceformula} says that 
\begin{equation}
\label{invsublogquaddiff}
V_i^{\varphi_i} \cong
H^0(\whW_i, 2K_{\whW_i}+\sum_{j=1}^{s_i+t_i} Q_j)^* \cap {\overline B},
\end{equation}
where ${\overline B}$ is a small open ball around the origin 
($=\overline{b_0}$) of 
$H^0(\whW_i, 2K_{\whW_i}+\sum_{j=1}^{s_i+t_i} Q_j)^*$.

\medskip

{\it Step 2} \quad We show the existence of a natural family of 
pointed Riemann surfaces over $V_i^{\varphi_i}$.
The connected component of the normalization of the restricted family 
$\psi^{-1}(V_i) \rightarrow V_i$ 
naturally induces a family of pointed Riemann surfaces
\begin{equation}
\label{isubkuranishi}
\psi_i:(\cS, \cP) \longrightarrow V_i, 
\end{equation}
which should be a standard Kuranishi family of the pointed 
Riemann surface $\psi_i^{-1}(b_0)=(\widehat{S}_i, \bP_i)$ 
by $H^1(\widehat{S}_i, T_{\widehat{S}_i}(-\bP_i)) \cong
H^0(\widehat{S}_i, 2K_{\widehat{S}_i}+\bP_i)^*$ (cf. \cite[p.177]{ACG}).
Therefore, the natural action of the stabilizer 
$G_{\varphi_i}={\rm Stab}(S_i)=\langle \varphi_i \rangle$ 
on $\widehat{S}_i$ is extended relatively to the family $\psi_i$.
We consider the relative quotient map
\begin{equation}
\label{quotfamimportant}
\overline{\psi}_i: \cW=\cS/G_{\varphi_i} \longrightarrow 
\overline{V_i}=V_i/G_{\varphi_i}  \cong V_i^{\varphi_i}.
\end{equation}
The central fiber $\overline{\psi}_i^{-1}(\overline{b_0})$ 
of (\ref{quotfamimportant}) is isomorphic to 
$\widehat{S}_i/G_{\varphi_i}=\widehat{W_i}$.
For any $\overline{b} \in \overline{V_i}$, 
the fiber  $\overline{\psi}_i^{-1}(\overline{b})$ is described as follows.
Since (\ref{isubkuranishi}) is a standard Kuranishi family, its fundamental property 
(cf. \cite[Chap.XI (6.13)]{ACG}) says that 
the fiber $\psi_i^{-1}(b)$ ($b \in V_i$) admits a subgroup 
$G_{\varphi_i, b} \subset {\rm Aut}(\psi_i^{-1}(b))$ which is isomorphic to $G_{\varphi_i}$.
Then
$$
\overline{\psi}_i^{-1}(\overline{b}) 
\cong \psi_i^{-1}(b)/G_{\varphi_i, b}:=\whW_{i, b}.
$$
In other words, the family (\ref{quotfamimportant}) is a deformation of 
$\whW_i=\overline{\psi}_i^{-1}(\overline{b_0})$ 
so that each fiber of $\overline{\psi}_i$ is a quotient surface 
of each fiber of $\psi_i$ in (\ref{isubkuranishi}) 
by essentially the same Galois group $G_{\varphi_i}$.
Moreover, each point $Q_j 
\in \{ Q_1, \cdots. Q_{s_i}, Q_{s_i+1}. \cdots, Q_{s_i+t_i} \}$ on $\widehat{W_i}$
is extended as a section $\overline{\cQ}_j:\overline{V_i} \rightarrow \cW$
of $\overline{\psi}_i$.

In fact, we set $\pi_{\cS}:\cS \rightarrow \cW=\cS/G_{\varphi_i}$ and 
$\pi_{\cS,\overline{b_0}}=\pi_{\cS}|_{(\pi_{\cS})^{-1}(\overline{b_0})}
:\widehat{S}_i \rightarrow \widehat{W_i}=\widehat{S}_i/G_{\varphi_i}$.
Then $Q_j$ satisfies one of the following two conditions:

\noindent {\rm (i)} \ There exists a point $P_{i,j}$ in the support of 
$\bP_i$ such that $Q_j=\pi_{\cS,\overline{b_0}}(P_{i,j})$.

\noindent {\rm (ii)} \ $Q_j$ is a branch point of the covering 
$\pi_{\cS,\overline{b_0}}$.

In case (i), since $\psi_i$ is a Kuranishi family of 
$(\widehat{S}_i, \bP_i)$, there exists a section $\bs:V_i \rightarrow \cS$ 
such that $\bs(V_i)$ passes through $P_{i,j}$.
Then the image $\pi_{\cS}(\bs(V_i))$ 
defines the desired section $\overline{\cQ}_j$.
In case  (ii),
since $G_{\varphi_i}$ and $G_{\varphi_{i, b}}$ have the same 
signature ($\lambda_1^{(i)}, \cdots, \lambda_{s_i}^{(i)}$) 
by \cite{H1971},   
the components of the discriminant locus of 
$\overline{\psi}_i: \cW \rightarrow \overline{V_i}$ are sections.
Then one of them is the desired $\overline{\cQ}_j$.

Thus the family $\overline{\psi}_i$ is a deformation of 
pointed Riemann surfaces 
($\widehat{W_i};Q_1, \cdots, Q_{s_i+t_i}$).
More strongly,  it follows from  (\ref{invsublogquaddiff}) that 
$\overline{\psi}_i$ is nothing but a standard Kuranishi family of 
($\widehat{W_i};Q_1, \cdots, Q_{s_i+t_i}$).

\medskip

{\it Step 3} \quad 
We fix a point $[S, w] \in T({\sigma}) \subset D_{\epsilon}(\sigma)$.
We may assume that the source space $\Sigma_g(\sigma)$ 
of the Weyl marking $w:\Sigma_g(\sigma) \rightarrow S$ 
is the topological model of a stable curve which has an automorphism 
$\varphi$ with the 
numerical data  (\ref{Eq0504no5b}), (\ref{defnumdata}).

We consider the Teichm\"{u}ller marking 
$w_i:(\hat{R}_i, \hat{\bP}_i) \longrightarrow (\hat{S}_i, \bP_i)$
in (\ref{partialmark}).
Then there exists a branched covering 
$\hat{\pi}_i:\hat{R}_i \rightarrow \Sigma_{\overline{g}_i}$ 
over a Riemann surface $\Sigma_{\overline{g}_i}$ of genus $\overline{g}_i$ 
such that the covering transformation group of $\hat{\pi}_i$ coincides with $G_{\varphi_i}$ 
in (\ref{smallcover}).
The cardinality of the set of points consisting of the branch points for $\hat{\pi}_i$ 
and $\hat{\pi}_i(\hat{\bP}_i)$ is $s_i+t_i$, and we write this set by $\{ \whQ_1, \cdots, \whQ_{s_i+t_i} \}$.
By the natural descent of $w_i$, we have a pointed oriented homeomorphism 
\begin{equation}
\label{decentpartialmark}
\hat{w}_i:(\Sigma_{\overline{g}_i}; \whQ_1, \cdots, \whQ_{s_i+t_i}) 
\longrightarrow 
(\widehat{W_i};Q_1, \cdots, Q_{s_i+t_i}).
\end{equation}
We may consider (\ref{decentpartialmark}) as a  Teichm\"{u}ller marking of 
the fiber $\overline{\psi}_i^{-1}(\overline{b}_0)$ of the family (\ref{quotfamimportant}).
Since  (\ref{quotfamimportant}) is a standard Kuranishi family, 
it follows from the discussion of Arbarello--Cornalba--Griffiths \cite[Chap.XV, \S2]{ACG} that 
this marking is extended to the whole family   (\ref{quotfamimportant}).
That is to say,  (\ref{quotfamimportant})  induces a family 
of Teichm\"{u}ller-marked pointed Riemann surfaces.

\medskip

{\it Step 4} \quad 
By applying 
the method  of the expression of pointed Teichm\"{u}ller 
spaces  via the patching of the base spaces of the standard Kuranishi families 
of pointed Riemann surfaces due to \cite[Chap.XV, \S2]{ACG} and also the fundamental property of the  action 
on  Kuranishi families, 
we globalize the discussion in Steps 1 $\sim$ 3. 

Let $T_{\sigma}^{\varphi}(p)$ be the connected component of 
$T_{\sigma}^{\varphi}$ containing the point $p=[S,w]$, 
where the marking $w$ is the one defined in Step 3.
This marking induces the marking of the $\varphi^j$-image 
of the normalized componet ($\hat{S}_i, \bP_i$) of $S$ as 
$\varphi^j \circ w_i:(\hat{R}_i, \hat{\bP}_i) 
\rightarrow (\varphi^j(\hat{S}_i), \varphi^j(\bP_i))$, where 
$w_i$ is also given in Step 3.
Then the point 
$$
p_{i,j}=[(\varphi^j(\hat{S}_i), \varphi^j(\bP_i)), \varphi^j \circ w_i],\:\:\:\:\:
(1 \leq i \leq \overline{r}, \:0 \leq j \leq m_i-1)
$$
is contained in the $\varphi_i(=\varphi^{m_i})$-invariant locus 
$T_{g_i, k_i}^{\varphi_i}$ of $T_{g_i, k_i}$.
For a fixed $i$, 
the connected components $T_{g_i, k_i}^{\varphi_i}(p_{i,j})$ of 
$T_{g_i, k_i}^{\varphi_i}$ containing $p_{i,j}$ ($0 \leq j \leq m_i-1$) 
are clearly isomorphic to each other, and we set 
$(T_{g_i, k_i}^{\varphi_i})^{(0)}:=T_{g_i, k_i}^{\varphi_i}(p_{i,0})$ 
for convenience.
Then, from (\ref{directsuminvsubsp}) and Lemma \ref{littlepatch},  
we have the 
composition of analytic embeddings
\begin{equation}
\label{firstembproof}
\prod_{i=1}^{\overline{r}} (T_{g_i, k_i}^{\varphi_i})^{(0)}
\longrightarrow 
\prod_{i=1}^{\overline{r}} 
T_{g_i, k_i}^{\varphi_i}(p_{i,0}) \times \cdots 
\times T_{g_i, k_i}^{\varphi_i}(p_{i,m_i-1})
\longrightarrow
\prod_{i=1}^{\overline{r}} 
\underbrace{T_{g_i, k_i}\times \cdots 
\times T_{g_i, k_i}}_{m_i}.
\end{equation}
From 
(\ref{quotfamimportant}) and the discussions in 
Step 3 and in \S \ref{autoriemannequi}, we have an isomorphism
\begin{equation}
\label{keyisom}
(T_{g_i, k_i}^{\varphi_i})^{(0)} \cong  
T_{\overline{g}_i, s_i+t_i}.
\end{equation}
Then the first assertion of (i) follows from  
(\ref{firstembproof}) and (\ref{keyisom}).
Since each of the set of connected components of the 
$T_{g_i, k_i}^{\varphi_i}$'s is countable by the same 
argument as in \cite{H1971},  
the set of connected components of 
$T_{\sigma}^{\varphi}$ is also countable.
Hence the assertion (i) holds.

Since the connected components of $T_{\sigma}^{\varphi}$ are isomorphic to 
each other so that these isomorphisms are given by the change of 
markings, 
they are mapped by 
the forgetting map  $T({\sigma}) \rightarrow \overline{M}_g$ of the 
markings  onto the same image, 
which is isomorphic to 
$\prod_{i=1}^{\overline{r}}M_{\overline{g}_i, s_i+t_i}$ 
from the assertion (i).
This is irreducible since each 
$\overline{M}_{\overline{g}_i, s_i+t_i}$ is irreducible.
By the same argument as in \cite[p.106]{Broughton1990} 
using the factorization of the forgetful map to the composition
of the infinite unramified covering and the finite Galois covering, 
it is a closed subvariety on $\overline{M}_g$.
Hence the assertion (ii)  holds.
\qed

\begin{coro}  {\rm (\cite{Terasoma1998})} \quad 
\label{coro0504no1}
$\overline{M}_g^{\varphi}$ is connected.
\end{coro}

\begin{remark}
\label{Rem0505no1}
The relation between the equisymmetric strata $T_{\sigma}^{\varphi}$ at the boundary 
and their limits on $M_g$ or $T_g$  seems to be interesting. See for instance 
\cite{DG2019}, etc.
\end{remark}


\subsection{Harris--Mumford   coordinates around equisymmmetric strata}
\label{HMcoord}

In this subsection, we define a special system of local coordinates 
on the controlled deformation space 
$D_{\epsilon}(\sigma)$ 
around an arbitrarily chosen point $p$ of the equisymmetric strata  
according to the method of \cite[\S1]{HM}. 

We fix a point $p=[S, w] \in T_{\sigma}^{\varphi} \subset T({\sigma}) 
\subset D_{\epsilon}(\sigma)$.
Since an open neighborhood of $p$ in $D_{\epsilon}(\sigma)$ 
is the base $B$ of a standard Kuranishi family of $S$, it follows form  
(\ref{locnbdbase}), (\ref{fundexactseq}) and (\ref{dualsp}) that 
the dual bases of $H^0(S, \Omega_{S}^1 \otimes \omega_{S})$ express
a system of  local coordinates at $p$ in $D_{\epsilon}(\sigma)$.
From (\ref{dualfundamentalseq}), it can be decomposed into the parameters 
for the variable defomations (I) and the smoothing deformations 
(II)  as in \S \ref{Subsec0301}.
By considering the action of $\varphi$ 
on  $H^0(S, \Omega_{S}^1 \otimes \omega_{S})$, 
we choose the basis as follows: 

(I) Let
$S=\sum_{i=1}^{\overline{r}}\sum_{j=0}^{m_i-1} \varphi^{j}(S_i)$ 
be the orbit-irreducible decomposition in (\ref{Eq0504no1}),
and
$\whS=\coprod_{i=1}^{\bar{r}} \coprod _{j=0}^{m_i-1} \widehat{\varphi^j(S_i)}$
be the natural decomposition of the normalization. 
As the dual to (\ref{defnonsmoothhere}), 
we consider 
the parameter space 
$\bigoplus_{i=1}^{\bar{r}} \bigoplus_{j=0}^{m_i-1} \varphi^j(V^*_i)$
of the variable deformations, where $V^*_i=H^0(\whS_i, 2K_{\whS_i}+\bP_i)$.
Now we choose the eigenbasis in Definition \ref{definition0502newno1} of the vector space $V^*_i$ 
with respect to the $\varphi^{m_i}$-action 
\begin{equation}
\label{equation0505no1}
\{ v_{i,1}, \cdots, v_{i, q_i} \} \in V^*_i, \:\:\:
(\varphi^{m_i})^*(v_{i, j})=\be\Big( \frac{\theta_{i,j}}{n_i} \Big) v_{i, j}, \:\:\:
(1 \leq j \leq q_i)
\end{equation}
where $q_i={\rm dim}\:V^*_i$ and $n_i$ is the order of the action of 
$\varphi^{m_i}$ on $\whS_i$.

With respect to the vector space $\varphi^j(V^*_i)$ for $1 \leq j \leq m_i-1$, 
we choose the basis  
$\{ \varphi^j(v_{i,1}), \cdots, \varphi^j(v_{i, q_i}) \}$. 
Then we have  
$(\varphi^{m_i})^*(\varphi^j(v_{i, j}))
=\be( \theta_{i,j}/{n_i}) \varphi^j(v_{i, j})$.

\medskip

(II) Let $\sum_{i=1}^{\bar{k}} \sum_{j=0}^{m(P_i)-1} \varphi^j(P_i)$  be the set of 
nodes of $S$ such that ${\rm Stab}_G(P_i)= \langle \varphi^{m(P_i)} \rangle$ 
and $k=\sum_{i=1}^{\bar{k}} m(P_i)$.
Let $v_{P_i}$ be the generator  of the torsion sheaf $\tau_{P_i}$ 
given in (\ref{torsionsheaf}).
Then the basis of the parameter space of the smoothing deformations is given by 
\begin{equation}
\label{equation0505no2}
\{  \varphi^j(v_{P_i}) \}_{1 \leq i \leq \bar{k}, 0 \leq j \leq m(P_i)-1} 
\:\: \in \:\:
\bigoplus_{i=1}^{\bar{k}} \bigoplus_{j=0}^{m(P_i)-1} \tau_{ \:\varphi^j(P_i)}.
\end{equation}
For the eigenvalues, we have the following. Here 
we write  $P=\varphi^j(P_i)$, $v_P= \varphi^j(v_{P_i})$, 
$(\delta^{(1)}/\lambda^{(1)})(P)=\delta^{(1)}/\lambda^{(1)}$ and so on for simplicity.
\begin{lemma}
\label{lemma0505no1} 
If $P$ is a non-amphidrome node such that 
the covalencies of both sides are $\delta^{(1)}/\lambda^{(1)}$ and $\delta^{(2)}/\lambda^{(2)}$, 
then 
$
(\varphi^{m(P)})^* v_P=\be \left( -\delta^{(1)}/\lambda^{(1)}-\delta^{(2)}/\lambda^{(2)}
\right) v_P.
$

If $P$ is an amphidrome node with  covalency $\delta/\lambda$, then 
$
(\varphi^{m(P)/2})^* v_P=\be \left( -\delta/\lambda \right) v_P.
$
\end{lemma}

{\it Proof} \quad Assume $P$ is non-amphidrome.
Note that $v_P$ is written as  $ydx^{\otimes 2}/x=xdy^{\otimes 2}/y$  
modulo $xy=0$ by  (\ref{torsionsheaf}). 
From (\ref{Eq0504no2b}), we have
$$
(\varphi^{m(P)})^*v_P=\frac{\be \left(-\frac{\delta^{(2)}}{\lambda^{(2)}} \right) y 
\cdot \be \left( -\frac{2\delta^{(1)}}{\lambda^{(1)}} \right)dx^{\otimes 2}}
{\be \left( -\frac{\delta^{(1)}}{\lambda^{(1)}} \right)x}
=\be \left( -\frac{\delta^{(1)}}{\lambda^{(1)}}-\frac{\delta^{(2)}}{\lambda^{(2)}}
\right) v_P.
$$
Siminarly in the case where $P$ is amphidrome, 
we have the desired result from  (\ref{Eq0504no2c}).
\qed

\bigskip

Based on this discussion, we define the local 
coordinates at $p$ 
of $D_{\epsilon}(\sigma)$ 
by noticing that 
$k=\sum_{i=1}^{\bar{k}}m(P_i)$ and 
 $3g-3-k=\sum_{i=1}^{\bar{r}}m_iq_i$.
 
\begin{definition}
\label{definition0505no1}
The system of Harris--Mumford  coordinates\index{Harris--Mumford  coordinates} 
$(z_1, \cdots, z_{3g-3})$ 
of $D_{\epsilon}(\sigma)$ 
at $p=[S, w] \in T_{\sigma}^{\wmu}$ is defined by the following: 

\noindent {\rm (i)} $p=\{ (z_1, \cdots,  z_{3g-3})=(0, \cdots, 0) \}$.

\noindent {\rm (ii)} We rewrite it by using the lexicographic order with respect to  
$i, j, \alpha,\beta, \gamma$ as 
\begin{equation}
\label{equation0505no3}
(z_1, \cdots,  z_{3g-3})=(z_1^{(0)}, \cdots, z_i^{(j)}, \cdots, z_{\bar{k}}^{(m(P_k)-1)}, 
z_{1,1}^{(0)}, \cdots, z_{\alpha, \beta}^{(\gamma)}, \cdots, z_{\bar{r}, q_{\bar{r}}}^{(m_{\bar{r}}-1)})
\end{equation}
for $1 \leq i \leq \bar{k}, \: 0 \leq j \leq m(P_i)-1, \: 1 \leq \alpha \leq \bar{r}, \: 1 \leq \beta  \leq q_{\alpha}, \: 
0 \leq \gamma \leq m_{\alpha}-1$.
Then the coordinate $z_i^{(j)}$ is the dual vector of $\varphi^j(v_{P_i})$ 
where $v_{P_i}$ is  given  in (\ref{equation0505no2}).
The coordinate $z_{\alpha, \beta}^{(\gamma)}$ is the dual vector of 
$\varphi^{\gamma}(v_{\alpha, \beta})$ which is given  in (\ref{equation0505no1}). 
\end{definition}

For the coordinates (\ref{equation0505no3}) on $B \subset D_{\epsilon}(\sigma)$ near $p$, 
the action of $\wmu$ on $B$ is written as follows. 
The proof of Lemma \ref{lemma0505no2} is obvious from 
Lemma \ref{lemma0505no1}  and (\ref{equation0505no1}).

\begin{lemma}
\label{lemma0505no2}
The action of $\wmu$ on $B$ near $p$ is written via the coordinates (\ref{equation0505no3}) as
$$
\wmu: z_i^{(j)} \longmapsto z_i^{(j+1)}
=\be \left( \frac{1}{N} \left( 
\frac{\delta^{(1)}}{\lambda^{(1)}}(P_i)+ \frac{\delta^{(2)}}{\lambda^{(2)}}(P_i) \right)
\right) z_i^{(j)}, \:\:
z_{\alpha, \beta}^{(\gamma)} \longmapsto z_{\alpha, \beta}^{(\gamma+1)}
=\be \left(  -\frac{\theta_{\alpha, \beta}}{N}
\right) z_{\alpha, \beta}^{(\gamma)}
$$
by identifying $z_i^{(m(P_i))}=z_i^{(0)}$ and  
$ z_{\alpha, \beta}^{(m_{\alpha})}= z_{\alpha, \beta}^{(0)}$.
\end{lemma}


\section{Monodromy and orbifold moduli maps of degenerations of Riemann surfaces}
\label{sec5}

 We discuss fundamental  properties of  the monodromy and the moduli maps 
of degeneration of Riemann surfaces of genus $g \geq 2$ from several points of view.

In \S6.1, we study pseudo-periodic maps of negative twist 
$\mu$, whose totality is denoted by $\bP_g^{(-)}$.
Our interest in  $\bP_g^{(-)}$ comes from the fact that 
the topological monodromy of a degeneration of a Riemann surface belongs to this class.
In this subsection, we show that $\mu$ is obtained from the lifting via the  real blow-up 
of an analytic automorphism $\mu^{{\rm an}}$ of a stable curve. 
It follows that the  fundamental invariants of $\mu$ (\cite{Nielsen2}, \cite{MM}) 
essentially come from 
those of $\mu^{{\rm an}}$ 
except for the {\it screw numbers} (\cite[Def.~2.4]{MM}), 
which express the fractional Dehn twists 
along the exceptional circles.
We also review the notion of {\it generalized quotient space} $\Sigma_g/\mu$  
consisting of {\it cores and non-core components} constructed  in \cite[Chap.~3]{MM}.
This gives  important 
information about the central fiber of a normally minimal 
degeneration whose topological monodormy coincides with $\mu$.

In \S6.2,
we define {\it the orbifold model} $f:S \rightarrow \Delta$ of a degeneration 
by contracting the non-core components of the central fiber 
which is topologically identified with the generalized quotient (\cite{MM}).  
Here $S$ is a normal complex space with at most quotient singularities such that 
the orbifold structure of $f$ is explicitly induced from the 
precise stable reduction $\wf:S \rightarrow \wDelta$
given in \cite[\S2]{A2010}.
Note that this type of  orbifold model historically originates from 
Imayoshi  \cite{Im81} from the viewpoint of Teichm\"uller theory.

The discussion in \S6.3 is the main part of this section.
We define the notion of  orbifold moduli map 
$J_{f}:\Delta \rightarrow \overline{M}_g^{orb}$ 
for an orbifold model $f$ of a degeneration, 
and show that $J_{f}$ has 
{\it the Kodaira-periodicity} property in the following sense.
For an elliptic fibration, 
Kodaira \cite{Kodaira} defined {\it the functional invariant} $J$ 
as the map from the base to the upper half-plane $\cH$ (=Teichm\H{u}ller space 
of $g=1$) by means of the elliptic modular function.
Since the monodormy in $SL(2, \bZ)$ (=the mapping class group of $g=1$) 
acts on $J$ around a degenerate fiber, the local expression of 
$J$ has a certain ``periodicity''.
This property is extended to $g \geq 2$ as   
the orbifold moduli maps $J_f$ and the action of the Weyl groups.  

In \S6.4,
we give two examples of degenerations of Riemann surfaces 
and 
their orbifold structures together with their orbifold 
moduli maps. 


\subsection{Pseudo-periodic maps and automorphisms of stable curves}
\label{ppsubsec}
We compare  
a pseudo-periodic map of negative twist\index{pseudo-periodic map of negative twist}
with an automorphism of a stable curve via the lifting given 
in \S \ref{LiftStableFN}.
For the basic terminologies, see \cite{MM}, \cite[\S1]{A2010}, \cite[\S5]{ImHand}.

The isotopy class of an oriented homeomorphism 
$\mu:\Sigma_g \rightarrow \Sigma_g$ is called a 
{\it pseudo-periodic map of negative twist 
with respect to a simplex $\sigma=\langle C_1, \cdots, C_k \rangle$  
(of Harvey's curve complex, \cite{Harvey})} if 

\noindent (i) \ $\mu$ preserves $\sigma$, and the restriction $\mu |_{B}$ to the 
complement $B=\Sigma_g \setminus \coprod_{1 \leq i \leq k} C_i$  is a periodic map, 
i.e. a certain power  $(\mu |_{B})^{N}$ is isotopic to the identity map ${\rm id}_B$. 

\noindent (ii) \ Let $m(C_i)$  ($1 \leq i \leq k$) be the minimal natural number such that 
$\mu^{m(C_i)}(\overrightarrow{C_i})=\overrightarrow{C_i}$ as  oriented curves.
Then $\mu^{m(C_i)}$ acts on an annular neighborhood $\cA_i$ of $C_i$ as 
a right-handed fractional Dehn twist.

Decomposing $B=\sum_{i=1}^{r} B_i$ into connected components, we call $B_i$ 
{\it a body component} (\cite[Def.4.8]{MM}), see also Fig.~4.1 of \cite{MM}.
We also call the minimal natural number $N$  with the property stated in (i) 
{\it the pseudo-period}\index{pseudo-period}
 of  $\mu$.
The map $\mu$ with   the property (i) is 
called a {\it pseudo-periodic map} (\cite[Def.1.1]{MM}), or 
in Bers'  terminology, a map of 
elliptic or parabolic type  (see \cite[\S3]{ImHand}).
We denote the subset of $\Gamma_{g}$  
consisting of pseudo-periodic maps of negative twist 
with respect to $\sigma$ by 
\begin{equation}
\label{ppmapsymbol}
{\bf P}_{g}^{-}(\sigma)
\: \subset \: \Gamma_{g}.
\end{equation}
The set of congujacy classes of  (\ref{ppmapsymbol}) is denoted by 
$\widehat{{\bf P}_{g}^{-}}(\sigma)
\: \subset \: \widehat{\Gamma_{g}}$, which is characterized as follows:


\begin{theorem} {\rm (\cite{Nielsen2}, \cite{MM})}
\label{nmm}
An element  $\mu \in \widehat{{\bf P}_{g}^{-}}(\sigma)$ 
is uniquely determined  by  

\noindent {\rm (a)} the Nielsen valencies\index{ Nielsen valency} at  
the multiple points  and  at the boundary curves 
of $\{ B_{i} \}$, 

\noindent {\rm (b)} the screw numbers\index{screw number} 
$\{ \bs(C _j) \}$ of  $\{ \cA _j \}$ 
where $\cA _j$ is a small annular neighborhood of 
$C_j$ such that $|\bs(C _j)|$ $(\bs(C _j) \leq 0)$ is  
the fractional Dehn twist of $\cA _j$, 

\noindent {\rm (c)} the action of $\mu$ on the extended partition 
graph\index{action on the extended partition graph} 
$\bGamma(\mu)$, i.e. the one-dimensional oriented graph 
whose vertices  correspond to $\{ B _i \}$ 
and whose edges correspond naturally to $\{ C _j \}$.

\end{theorem}

The invariants (a) and  (b) are due to \cite{Nielsen2}, and 
(c) is due to \cite{MM}.
We call them {\it the numerical data (a), (b), (c)} of $\mu$.
The pseudo-periodic map and the analytic automorphism of a stable curve are related 
as follows.
\begin{prop}
\label{autoppmapequiv}
Let $\mu: \Sigma_g \rightarrow \Sigma_g$ be an element of  ${\bf P}_{g}^{-}(\sigma)$.
Then there exists an analytic automorphism 
$\mu_\sigma: \Sigma_g(\sigma) \rightarrow  \Sigma_g(\sigma)$ 
with respect to some complex structure on the stable curve 
$\Sigma_g(\sigma)$ such that the following diagram is isotopically commutative;

\begin{picture}(195,71)(-35,-3)
\put(80,47){\makebox(6,6){$\Sigma_g$}}
\put(80,5){\makebox(6,6){$\Sigma_g(\sigma)$}} 
\put(165,48){\makebox(8,6){$\Sigma_g$}} 
\put(167,4){\makebox(8,6){$\Sigma_g(\sigma)$.}} 
\put(81,43){\vector(0,-1){24}} 
\put(170,43){\vector(0,-1){24}} 
\put(99,51){\vector(1,0){51}}
\put(99,7){\vector(1,0){51}}
\put(61,26){\makebox(6,6){${\rm cont}_{\sigma}$}} 
\put(185,26){\makebox(6,6){${\rm cont}_{\sigma}$}} 
\put(122,-4){\makebox(6,6){$\mu_\sigma$}} 
\put(121,40){\makebox(6,6){$\mu$}} 
\end{picture}

\noindent Conversely, any analytic automorphism $\mu_\sigma$ of a stable curve 
$\Sigma_g(\sigma)$ is 
lifted to an element $\mu$ of   ${\bf P}_{g}^{-}(\sigma)$.
\end{prop}

{\it Proof} \ ({\it Step 1}) \  If $\sigma=\emptyset$, then $\mu$ is a periodic map,  and it is known 
that $\mu$ is isotopic to  an analytic automorphism.
This fact is  a corollary of Kerchhoff's theorem \cite{Ke}, 
or this isotopy  is explicitly  constructed, see \cite{Nielsen3} or \cite[Chap.2]{MM}.

If $\sigma \not= \emptyset$, $\mu$ is isotopic on $B=\Sigma_g \setminus \coprod_{1 \leq i \leq k} \cA_i$ 
to a direct sum of  
analytic automorphisms of Riemann surfaces with  real boundaries, and 
the restrictions to  $\coprod \cA_i$
could shrink to point maps (in $\Sigma_g(\sigma)$),  
just as in  Prop.~\ref{isotopicfiniteangled}.
Since $\Sigma_g(\sigma)$ is  constructed  from $\Sigma_g$ 
by shrinking $\coprod \cA_i$ to nodes via isotopy, 
the analytic automorphism on $B$ descends to an analytic automorphism of 
$\Sigma_g(\sigma)$. 

 ({\it Step 2}) \ Conversely, let $\mu_\sigma: \Sigma_g(\sigma) \rightarrow  \Sigma_g(\sigma)$ 
be an analytic automorphism with respect to some complex structure.\par
Assume that $P_i={\rm cont}(C_i)$ is a non-amphidrome node with  co-valencies
$\delta^{(1)}/\lambda^{(1)}$ and $\delta^{(2)}/\lambda^{(2)}$ at the disk neighborhoods 
$U^{(1)}$  and $U^{(2)}$ of the local components of both sides, where 
$P_i=U^{(1)} \cap U^{(2)}$.
By definition, $\mu_\sigma^{m(P_i)}$ preserves $U^{(j)}$ and rotates it by the angle $2\pi \delta^{(j)}/\lambda^{(j)}$ 
clockwise, within the view from the insides of both components ($j=1,2$).
Let $\pi_P^{(j)}:\wU^{(j)} \rightarrow U^{(j)}$ be the real blowing up at $P_i$ with the exceptional 
circle $C^{(j)}=(\pi_P^{(j)})^{-1}(P)$, 
and $\cA$  a small annulus with  boundary $\partial \cA=\partial \cA^{(1)} \coprod \partial \cA^{(2)}$.

We construct a part of a Riemann surface 
$(\Sigma_g)_{{\rm loc}}=\wU^{(1)} \cup \cA \cup \wU^{(2)} $ from the building blocks 
$\wU^{(1)}, \:\wU^{(2)}$ and $\cA$ by pasting $C^{(j)}$ to $\partial \cA^{(j)}$ naturally, 
and define a homeomorphism  
$\mu_{{\rm loc}}:(\Sigma_g)_{{\rm loc}} \rightarrow (\Sigma_g)_{{\rm loc}}$ 
which is a local lift of $\mu_\sigma^{m(P_i)}$ as follows. 
The restriction 
$\mu_{{\rm loc}}|_{\wU^{(j)}}: \wU^{(j)} \rightarrow \wU^{(j)}$ is defined  to be  the lifting  of 
$\mu_\sigma^{m(P_i)}$ via $\pi_P^{(j)}$, and 
$\mu_{{\rm loc}}|_{\cA}: \cA \rightarrow \cA$ is defined to be the fractional Dehn twist, 
i.e. the  linear twist (\cite[\S2.3]{MM})  with  
{\it screw number} \cite[Def.~2.4]{MM} 
\begin{equation}
\label{screwnonamp}
\bs(C_i) =-\frac{\delta_i^{(1)}}{\lambda_i^{(1)}}-\frac{\delta_i^{(2)}}{\lambda_i^{(2)}}-\bK(C_i)\:\:\:\:\:\:\: 
(\bK(C_i) \in \bZ, \:\bK(C_i) \geq -1).
\end{equation}
Since $-\delta^{(2)}/\lambda^{(2)} \equiv \delta^{(1)}/\lambda^{(1)}+\bs(C_i)$ (${\rm mod} \: \bZ$), 
the map  $\mu_{{\rm loc}}$ is well-defined and expresses a right-handed fractional Dehn twist 
with $\bs(C_i) \leq 0$.

When $P_i$ is an amphidrome node with  co-valency
$\delta/\lambda$ of the local components of both sides, 
then by the same notations as above we define the  homeomorphism 
$\mu_{{\rm loc}}:(\Sigma_g)_{{\rm loc}} \rightarrow (\Sigma_g)_{{\rm loc}}$ 
which is the local lift of $\mu_\sigma^{m(P_i)/2}$ as follows.
The restrictions 
$\mu_{{\rm loc}}|_{\wU^{(1)}}: \wU^{(1)} \rightarrow \wU^{(2)}$ 
and $\mu_{{\rm loc}}|_{\wU^{(2)}}: \wU^{(2)} \rightarrow \wU^{(1)}$ 
are defined  to be the liftings  of 
$\mu_{\sigma}^{m(P_i)/2}$ via $\pi_P^{(1)}$ and $\pi_P^{(2)}$,
and 
$\mu_{{\rm loc}}|_{\cA}: \cA \rightarrow \cA$ is defined to be   
the special twist (\cite[\S2.4]{MM})  which interchanges the boundary 
components $\partial \cA^{(1)}$ and $\partial \cA^{(2)}$ 
with screw number $\bs(C_i)/2$, where 
\begin{equation}
\label{screwamp}
\bs(C_i)=-\frac{2\delta_i}{\lambda_i}-2\bK(C_i)  
\:\:\:\:\:\:\:  (\bK(C_i) \in \bZ, \bK(C_i) \geq 0).
\end{equation}

The restrictions of $\mu_\sigma$ to parts of $\Sigma_g(\sigma)$, other 
than the disk neighborhoods of the nodes, have trivial liftings.
Following the action of $\mu_\sigma$ 
on the dual graph of $\Sigma_g(\sigma)$, 
we patch these parts of Riemann surfaces and the homeomorphisms.
Then we obtain  $\Sigma_g$ and 
the desired homeomorphism  $\mu$. 
\qed

\bigskip
Note that (\ref{screwnonamp}) and  (\ref{screwamp}) 
are the screw numbers in the data (b) of Theorem~\ref{nmm}.
By the descent from $\mu$ to $\mu_\sigma$, 
the data (a) and (c) are preserved, 
while the data (b) 
vanish by this descent, i.e. 
the screw numbers 
are  not 
the data of the analytic 
automorphism of the stable curve.
(The screw numbers express essentially 
{\em the fractional Dehn twist coefficients}\index{fractional Dehn twist coefficient} 
along the exceptional circles, cf. Liu \cite{Liu2021}.)

With respect to this complex structure on $\Sigma_g(\sigma)$, 
there exists a quotient holomorphic map 
$\pi_{{\mu}_\sigma}: \Sigma_g(\sigma) 
\longrightarrow \Sigma_g(\sigma)/G_{\mu_\sigma}$
(Def.~\ref{Def0504no1}).
This quotient  map is topologically lifted 
to {\it the generalized quotient map}\index{generalized quotient map}
 \cite[Def.3.4]{MM}
\begin{equation}
\label{standardcoverlift}
\pi_{\mu}: \Sigma_g(\sigma) 
\longrightarrow \cW(\mu)=\Sigma_g/ \langle \mu \rangle.
\end{equation}
Here $\cW(\mu)$ is a non-reduced nodal Riemann surface in general 
(i.e.~a multiplicity is attached to each component) 
and $\pi_{\mu}$ is a pinched covering (\cite{MM}, Def.3.2), i.e. 
a finite unramified topological covering 
except over  nodes.
Each connected component of  
the inverse image of a node is homeomorphic to a
circle.
The space  $\cW(\mu)$ is decomposed into the parts (\cite{MM}, p.95, Fig.4.1);
\begin{equation}
\label{partdecomp}
\cW(\mu)=\sum_i {\rm Core}_i+
\sum_j {\rm Tail}_j+
\sum_j {\rm Arch}_j+
\sum_j {\rm Quasitail}_j.
\end{equation}
The orbit $\mu^{\alpha}(B_i)$ ($0 \leq \alpha \leq {m(B_i)-1}$) of 
$B_i$ is mapped to ${\rm Core}_i$ by $\pi_{\mu}$ 
except for small disk neighborhoods of multiple points.
The parts
${\rm Tail}_j,\:{\rm Arch}_j,\: {\rm Quasitail}_j$ are chains or trees of 
$\bP^{1}$'s which are the images under  
$\pi_{\mu}$ of the neighborhoods 
of multiple points, non-amphidrome annuli and amphidrome annuli, respectively,  
and these multiplicities 
are explicitly 
determined from the data (a) and (b).

Let ${\rm cont}^{(nc)}: \cW(\mu) \rightarrow 
\cW(\mu)^{\sharp}$ be the 
contraction of {\it the non-cores}\index{non-cores}
 (i.e. all the Tails\index{tail}, Archs\index{arch} and 
 Quasitails\index{quasitail}). 
Then $\cW(\mu)^{\sharp}$ can be identified with 
$ \Sigma_g(\sigma)/G_{\mu_{\sigma}}$ (see Remark \ref{nonreducedquot}), 
in other words, the following homotopically commutative diagram is the lifting 
of $\pi_{\mu_{\sigma}}$ 
to $\pi_{\mu}$;

\begin{picture}(195,71)(-35,-3)
\put(80,47){\makebox(6,6){$\Sigma_g$}}
\put(80,5){\makebox(6,6){$\Sigma_g(\sigma)$}} 
\put(165,48){\makebox(8,6){$\cW(\mu)$}} 
\put(183,4){\makebox(24,6){$\cW(\mu)^{\sharp} \cong \Sigma_g/G_{\mu_{\sigma}}$.}} 
\put(81,43){\vector(0,-1){25}} 
\put(165,43){\vector(0,-1){25}} 
\put(99,51){\vector(1,0){51}}
\put(99,7){\vector(1,0){51}}
\put(61,26){\makebox(6,6){${\rm cont}_{\sigma}$}} 
\put(188,28){\makebox(6,6){${\rm cont}^{(nc)}$}} 
\put(122,-4){\makebox(6,6){$\pi_{\mu_\sigma}$}} 
\put(121,40){\makebox(6,6){$\pi_{\mu}$}} 
\end{picture}

\vspace{-0.4cm}

\begin{center}
(Diagram I) \ The lifting of the analytic quotient to the generalized quotient
\end{center}


\subsection{Orbifold structures of  degenerations of Riemann surfaces}
\label{orbifoldmodeldeg}
In this subsection, we propose the notion of  
{\it orbifold model}\index{orbifold model}
of a degeneration\index{degeneration} of 
Riemann surfaces.

Let $M$ be a 2-dimensional normal complex space and 
$f:M \rightarrow \Delta=\{ t \in \bC \: | \: |t| \leq \epsilon_0 \}$
a proper surjective holomorphic map 
to a disc with a sufficiently small radius such that 
any fiber $f^{-1}(t)$ over $t \in \Delta^*=\Delta \setminus \{ 0 \}$ is a Riemann surface of genus $g$.
We call $f$ {\it a degeneration of genus $g$}.
We always assume $g \geq 2$ unless otherwise mentioned.

First we assume that $M$ is nonsingular.
Let $F=f^{-1}(0)=\sum_{i=1}^{r_0} m_iF_i$ be the irreducible decomposition of the central fiber.
If the reduced scheme $F^{{\rm red}}= \sum_{i=1}^{r_0} F_i$ has at most nodal singularities
such that 
any ($-1$)-spherical component of $F^{{\rm red}}$ has at least three intersection points with other components, 
then $f$ is called {\it the normally minimal model} which is uniquely determined 
in the local birational equivalence class of $f$.

Two degenerations of genus $g$ are {\it topologically (resp. analytically) equivalent} 
if, with  their normally minimal models 
$f:M \rightarrow \Delta$ and $f':M' \rightarrow \Delta'$, there exist 
oriented homeomorphisms (resp. analytic isomorphisms) 
$h_M:M \rightarrow M'$ and $h_{\Delta}: \Delta \rightarrow \Delta'$ 
such that $h_{\Delta} \circ f=f' \circ h_M$.
The equivalence class is called 
the {\it  topological type} (resp. the {\it analytic type}) of  
the degeneration.

Now we fix a point $t_0 \in \Delta^*$, and let $w:\Sigma_g \rightarrow f^{-1}(t_0)$ 
be a Teichm\"{u}ller marking.
We choose a smooth loop $L_{t_0t_0}$ in $\Delta^*$ 
which starts from $t_0$, goes around $0 \in \Delta$ once counterclockwise 
and comes back to $t_0$. 
Let $\mu(L_{t_0t_0}):f^{-1}(t_0) \rightarrow f^{-1}(t_0)$ be
the diffeomorphism along $L_{t_0t_0}$.
Then 
$w^{-1} \circ \mu(L_{t_0t_0}) \circ w: \Sigma_g \rightarrow \Sigma_g$ 
is an oriented homeomorphism. 
For another point $t'_0 \in \Delta^*$, we define the Teichm\"{u}ller marking of 
$f^{-1}(t'_0)$ coming from $w$, by 
$
w'=\mu(L_{t_0t'_0}) \circ w:\Sigma_g \longrightarrow f^{-1}(t'_0)
$
where $\mu(L_{t_0t'_0}):f^{-1}(t_0) \rightarrow f^{-1}(t'_0)$ is 
the diffeomorphism along a path $L_{t_0t'_0}$ in $\Delta^*$
connecting $t_0$ and $t'_0$.
Then the map $w'^{-1} \circ \mu(L_{t'_0t'_0}) \circ w'$ 
is conjugate to $w^{-1} \circ \mu(L_{t_0t_0}) \circ w$.
In other words, the conjugacy class of
\begin{equation}
\label{topmonodwithmark}
\mu_f(w)=[w^{-1} \circ \mu(L_{t_0t_0}) \circ w] 
\end{equation}
is well-defined, independently of the choices of $t_0$ and $L_{t_0t'_0}$. \par
The proof goes as follows: We take a ``straight''  path  $L'_{t_0t'_0}$ in 
$\Delta^*$ connecting $t_0$ and $t'_0$, and we suppose that the loop 
$L'_{t_0t'_0}L_{t'_0t'_0}L'_{t'_0t_0}$ is isotopic to $L_{t_0t_0}$. We will denote the diffeomorphisms 
$\mu(L_{t_0t'_0})$ and $\mu(L'_{t_0t'_0})$ by $a$ and $b: f^{-1}(t_0) \to f^{-1}(t'_0)$, 
respectively. Then by the asumption on $L'_{t_0t'_0}$, we have $b^{-1}\mu(L_{t'_0t'_0})b 
=\mu(L_{t_0t_0})$, and 
\begin{equation}\label{conjugacyclass}
\begin{split}
&[w'^{-1}\circ \mu(L_{t'_0t'_0})\circ w']\\
&=[w^{-1}a^{-1}\circ \mu(L_{t'_0t'_0})\circ aw]\\
&=[w^{-1}a^{-1}bb^{-1}\circ \mu(L_{t'_0t'_0})\circ bb^{-1}aw]\\
&=[w^{-1}a^{-1}b\circ (b^{-1}\mu(L_{t'_0t'_0})b)\circ b^{-1}aw]\\
&=[w^{-1}a^{-1}b\circ \mu(L_{t_0t_0})\circ b^{-1}aw]\\
&=[w^{-1}(a^{-1}b)ww^{-1}\circ\mu(L_{t_0t_0})\circ ww^{-1}(b^{-1}a)w]\\
&=[c^{-1}w^{-1}\circ\mu(L_{t_0t_0})\circ wc]
\end{split}
\end{equation}
where we put $c=w^{-1}(b^{-1}a)w : \Sigma_g \to \Sigma_g$. The equation (\ref{conjugacyclass}) shows that the cojugacy 
classes $[w'^{-1}\circ\mu(L_{t'_0t'_0})\circ w']$ and $[w^{-1}\circ\mu(L_{t_0t_0})\circ w]$ coincide. 
We call this conjugacy class $\mu_f(w)$ 
{\it the topological monodromy of $f$ with respect to the marking $w$}
\index{topological monodromy with respect to marking}.
Note  that $\mu_f(w)$ is determined by the fibering structure of 
$f^{-1}(\Delta^*) \rightarrow \Delta^*$.

For the same complex structure of $f^{-1}(t_0)$, 
we choose another Teichm\"{u}ller marking 
$\ww:\Sigma_g \rightarrow f^{-1}(t_0)$
and consider the topological monodormy $\mu_f(\ww)$ with respect to $\ww$.
Then $\mu_f(w)$ and $\mu_f(\ww)$ are obviously conjugate to each other in $\Gamma_g$, i.e.
$$
\mu_f(\ww)=[w^{-1}  \circ \ww]^{-1} \mu_f(w) [w^{-1} \circ \ww] \in \Gamma_g.
$$
Therefore, the conjugacy class $\widehat{\mu_f(w)}$ of 
$\mu_f(w)$ is uniquely determined by $f$, independently of 
the choice of the marking $w$. We write
\begin{equation}
\label{topmonodwithoutmark}
\mu_f=\widehat{\mu_f(w)} \in \widehat{\Gamma}_g \:\:\: 
({\rm conjugacy}\: {\rm classes} \: {\rm of}\: \Gamma_g),
\end{equation}
and call $\mu_f$ 
{\it the topological monodromy of $f$}\index{topological monodromy}.
Then:

\begin{theorem}
\label{monodpseudo}
{\rm (i)} (\cite{Clemens}, \cite{Ac}, \cite{Im81}, \cite{ST}, \cite{ES})  
The topological monodromy belongs to the conjugacy classes of pseudo-periodic maps of 
negative twist.

\noindent {\rm (ii)} (\cite[Th.7.2]{MM})
The topological structure of a degeneration of genus $g$ is 
uniquely determined by its topological monodromy.
The set of  topological monodromies and the set of conjugacy classes of pseudo-periodic maps of 
negative twist correspond bijectively.
\end{theorem}

The relation between the normally minimal model 
and the generalized quotient  (\ref{standardcoverlift}) 
of  the topological 
monodormy $\mu_{f}$ 
is the following:


\begin{theorem} (\cite[Chap.9]{MM}) 
\label{centralfibquotient}
The central fiber $F=f^{-1}(0)$ of the normally minimal model
$f: M \rightarrow \Delta$ 
of a degeneration of genus $g$ and the generalized quotient space $\cW(\mu_f)$ of 
$\mu_f$ topologically coincide with each other 
as  non-reduced nodal Riemann surfaces.
In paticular,  $F$ is decomposed into cores and non-cores as in 
(\ref{partdecomp}). 
\end{theorem}

Now we change the model in the birational equivalence class of $f$.
Since any proper subset of irreducible components of 
the central fiber $F$  of the normally minimal model $f$ has 
negative intersection form, there exists the analytic contraction map 
$\tau: M \rightarrow M^{\sharp}$ of the  \underline{non-cores} of $F$ 
by Grauert's theorem \cite{Gr}.
Let $f^{\sharp}:M^{\sharp} \rightarrow \Delta$ be the natural map.
The restriction of $\tau$ to $F$ coincides with the contraction map 
$$
\tau |_{F}={\rm cont}^{(nc)}: F \longrightarrow F^{\sharp}=(f^{\sharp})^{-1}(0)
$$
given in Diagram I of \S6.1.
The total space $M^{\sharp}$ is a normal complex space with cyclic and dihedral 
quotient singularities whose supports are on the contraction points of the non-cores,  
and the types of singularities are explicitly given by the data 
of $\mu_f$ (\cite[Lem.~2.1]{A2010}) such that 
the non-cores in $F$ are the exceptional set of the 
minimal resolution of these singularities.
Since $M^{\sharp}$ has at most quotient singularities, we may consider it as 
a complex orbifold in the sense of Satake \cite{Sa}.
We call $f^{\sharp}$ {\it the orbifold model}
of the degeneration $f$\index{orbifold model of degeneration}.

The explicit orbifold structure of $f^{\sharp}$ is described as follows  
(\cite[\S\S 2,3]{A2010}).
Let $\pi: \widetilde{\Delta} \rightarrow \Delta$ be the covering map of disks 
defined by $u \mapsto t=u^N$ where $N$ is the pseudo-period of $\mu_f$, 
and let $\widetilde{M}$ be the normalization of the fiber product 
of $M^{\sharp}$ and $\widetilde{\Delta}$ over $\Delta$.
Let $\wf:\wM \rightarrow \wDelta$ be the natural map.
Then the covering transformation group 
$G:=\langle \widetilde{\mu} \rangle \simeq \bZ/N\bZ$ 
acting on 
$\wf^{-1}(\wDelta^*) \rightarrow \wDelta^*=\wDelta \setminus \{ 0 \}$ 
is extended holomorphically over $\wf$ by the normality of $\wM$ 
such that the following commutative diagram holds:

\setlength{\unitlength}{0.90mm}
\begin{picture}(160,29)(-14,3)

\put(50,0){\makebox(6,6)[l]{$\Delta \simeq \wDelta/G$}} 
\put(50,20){\makebox(6,12)[l]{$M^{\sharp} \simeq \wM/G$}} 
\put(53,21){\vector(0,-1){14}} 
\put(19,23){\makebox(6,6){$M$}} 
\put(26,26){\vector(1,0){22}} 
\put(24,23){\vector(4,-3){25}} 
\put(31,7){\makebox(6,6){$f$}} 
\put(97,23){\makebox(6,6){$\wM$}} 
\put(96,25){\vector(-1,0){21}} 
\put(97,1){\makebox(6,6){$\wDelta$}} 
\put(96,3){\vector(-1,0){25}} 
\put(100,21){\vector(0,-1){13}} 
\put(102,11){\makebox(6,6){$\wf$}} 
\put(83,0){\makebox(6,12)[l]{$\pi$}} 
\put(83,22){\makebox(6,12)[l]{$\widetilde{\pi}$}} 
\put(33,20){\makebox(6,6){$\tau$}} 
\put(54,12){\makebox(6,6){$f^{\sharp}$}} 

\vspace{-0.3cm}

\end{picture}
\begin{center}
(Diagram II) \ The process of the precise stable reduction
\end{center}
Since the topological monodromy $\mu_{\wf} \simeq \mu_f^N$ 
is trivial on the {\it body} $B$ (as defined before Th.6.1), the non-cores of 
the generalized quotient $\cW(\mu_{\wf})$ consists of {\it archs} (\cite[Def.4.7]{MM})
of multiplicity $1$, 
i.e.  $\cW(\mu_{\wf})$ is a semi-stable curve of genus $g$.
Then the central fiber $\wF=\wf^{-1}(0)$ 
is a stable curve with the topological type 
$\Sigma_g(\sigma)$ such that 
$\wM$ has rational double points of type A at the nodes of $\wF$.
From the viewpoint of the semi-stable curve $\cW(\mu_{\wf})$, 
$\wF$ is  the image of the contraction of the archs on it.
The restriction to the fiber 
$\widetilde{\pi}|_{\wF}:\wF\rightarrow  F^{\sharp}$ 
of the map $\widetilde{\pi}$ in Diagram II 
essentially coincides with the quotient map $\pi_{\mu_{\sigma}}$ in Diagram I.
We call $\wf$ {\it the precise stable reduction}\index{precise stable reduction}
 of $f$ (\cite{A2010}).
This stable reduction is minimal in the sense that the degree $N$ of 
the base change 
is minimal among all the stable reductions of $f$, because the generalized quotient 
for $\mu_f^n$ with $n < N$ cannot be a semi-stable curve of genus $g$ 
by the algorithms given in \cite[Chap.3]{MM}.

The action of the generator $\wmu$ of $G$ on $\wf$ near $\wF$ 
is an extension of the automorphism 
$\wmu|_{\wF}: \wF \rightarrow \wF$ of the stable curve $\wF$ 
to that of the family $\wf$,
and is locally described as follows:

(i)
Assume that $P$ is a nonsingular point of $\wF$ which 
belongs to an irreducible component 
$\wF_i$ with  ${\rm Stab}_G(\wF_i)= \langle \wmu^{m(\wF_i)} \rangle$.
The automorphism 
$\wmu^{m(\wF_i)}|_{\wF_i}:\wF_i \rightarrow \wF_i$ of order $n(\wF_i)$ 
and the associated 
$n(\wF_i)$-fold cyclic covering 
$\wpi|_{\wF_i}:\wF_i \rightarrow \wpi(\wF_i) \subset F^{\sharp}$ 
are given in 
(\ref{Eq0504no2}) and (\ref{smallcover}).
Let ($x, u$) be the local coordinates at 
$P=\{ (x,u)=(0,0) \}$ 
of $\wM$ such that $x$ is the fiber coordinate and $u$ is the lift 
of the base coordinate 
$\wDelta=\{ u \in \bC \:|\: |u| \leq \epsilon^{1/N} \}$.
Then the map $\wmu^{m(\wF_i)}$ locally acts as
\begin{equation}
\label{eq0602x00}
\wmu^{m(\wF_i)}: \:(x,u) \longmapsto
\left( \wmu^{m(\wF_i)}|_{\wF_i}(x), \:
\be \left( \frac{m(\wF_i)}{N} \right) u \right).
\end{equation}
Moreover if $P$ is the ramification point of $\wpi|_{\wF_i}$ 
with  co-valency $\delta/\lambda$, then 
${\rm Stab}_G(P)= \langle \wmu^{m(\wF_i)n(\wF_i)/\lambda} \rangle$ 
and the map $\wmu^{m(\wF_i)n(\wF_i)/\lambda}$ acts near $P$ as 
\begin{equation}
\label{eq0602x01}
\wmu^{m(\wF_i)n(\wF_i)/\lambda}: \:(x,u) \longmapsto
\left( \be \left( \frac{\delta}{\lambda} \right) x, \:
\be \left( \frac{m(\wF_i)n(\wF_i)}{\lambda N} \right) u \right).
\end{equation}

(ii) 
Assume $P$ is a non-amphidrome node of $\wF$ with 
${\rm Stab}_G(P)= \langle \wmu^{m(P)} \rangle$.
By (\ref{screwnonamp}), 
the screw number at the curve $C_P={\rm Cont_{\sigma}}^{-1}(P)$ 
via the map 
${\rm Cont}_{\sigma}: \Sigma_g \rightarrow \Sigma_g(\sigma) \simeq \wF$ 
is given by 
$\bs(C_P)=-\delta^{(1)}/\lambda^{(1)}-\delta^{(2)}/\lambda^{(2)}-\bK$.
The monodromy map $\mu_{\wf}\simeq \mu_f^N$ behaves on an 
anular neighborhood of $C$ as the right-handed integral Dehn twist of times 
\begin{equation}
\label{equation0602x2}
n(P):=\frac{N|\bs(C_P)|}{m(P)}=
\frac{\ell
(\delta^{(1)} \lambda^{(2)}+ 
\delta^{(2)}\lambda^{(1)} + 
\bK\lambda^{(1)} \lambda^{(2)})}
{{\rm gcd}(\lambda^{(1)}, \lambda^{(2)})},
\end{equation}
where $\ell=N/\bigr({\rm lcm}(\lambda^{(1)},\lambda^{(2)}) m(P) \bigl) \in \bZ$.
The restriction of this map to the tubular neighborhood of $P$ 
gives the Milnor fibration  of the singularity $P$ 
whose monodormy map is the one described above.
This means that $\wM$ is defined locally near $P$ by the equation 
$xy=u^{n(P)}$ where $x$ and $y$ are local parameters of the 
components of both sides of $P$.
Here $n(P)$ coincides with the Milnor number of the singularity $P$ of $\wM$.
The map $\wmu^{m(P)}$ is given locally near $P$ 
from (\ref{Eq0504no2b})
by
\begin{equation}
\label{equation0602x3}
\wmu^{m(P)}:\:(x,y,u) \longmapsto 
\left( \be \left( \frac{\delta^{(1)}}{\lambda^{(1)}} \right) x, \:
\be \left( 
\frac{\delta^{(2)}}{\lambda^{(2)}} \right) y, \:
\be \left( 
\frac{m(P)}{N} 
\right) u
\right).
\end{equation}

(iii) 
Assume $P$ is an amphidrome node with ${\rm Stab}_G(P)
= \langle \wmu^{m(P)/2} \rangle$. 
The screw number is given by $\bs(C_P)=-2(\delta/\lambda)-2\bK$ 
from (\ref{screwamp}).
Then $\wM$ is defined locally near $P$ by the equation 
$xy=u^{n(P)}$ such that the Milnor number\index{Milnor number}
 $n(P)$ is given by
\begin{equation}
\label{equation0602x4}
n(P)=\frac{N|\bs(C_P)|}{m(P)}=2\ell \bigr( \delta + \bK \lambda \bigl),
\end{equation}
where $\ell=N/(\lambda m(P)) \in \bZ$.
The map $\wmu^{m(P)/2}$ is locally given from (\ref{Eq0504no2c}) by
\begin{equation}
\label{equation0602x5}
\wmu^{m(P)/2}: \: (x,y,u) \longmapsto 
\left( \be \left( \frac{\delta}{2\lambda} \right) y, \:
\be \left( \frac{\delta}{2\lambda} \right) x, \:
\be \left( \frac{m(P)}{2N} \right) u \right).
\end{equation}

Depending  on the above arguments, we define the following;
\begin{definition}
\label{definition0602no1}
The {\em marked orbifold structure} of a degeneration $f$ is defined by;  

\noindent {\rm (i)} the set consisting of the orbifold model, the precise stable 
reduction and the group 
\begin{equation}
\label{equation0602x6}
\{ 
 \wf:\wM \rightarrow \wDelta, \: f^{\sharp}: M^{\sharp} \rightarrow \Delta, \:
G=\bZ/N\bZ
\}
\end{equation}
satisfying the Diagram II 
so that  the action of $G$ on $\wf$ 
is locally induced from  (\ref{eq0602x00}) $\sim$ (\ref{equation0602x5}), and

\noindent {\rm (ii)} the lifting of the Weyl marking of the stable curve 
$\wF=\wf^{-1}(0)$ given by
\begin{equation}
\label{equation0602x7}
\ww={\rm cont}^{({\rm nc})} \circ \pi_{\mu_{\wf}}  \: : \:
\Sigma_g \longrightarrow \cW(\mu_{\wf}) \longrightarrow \wF.
\end{equation}
We sometimes write this marked orbifold structure by 
$f^{{\rm orb}}= \{ \wf, f^{\sharp},G, \ww \}$ for simplicity.
\end{definition}

\begin{remark}
\label{correctorb}
{\rm (i)} 
Since $\wM$ has singularities, ($\wf, G$) does not define directly the complex orbifold structure on $f^{\sharp}$. 
However this point is easily recovered by Takamura's method \cite{Tak}:
the actions (\ref{equation0602x3}) and (\ref{equation0602x5}) 
near the nodes are lifted to the local linear actions of type A 
of $\bC^2$ explicitly.
Hence the orbifold charts of $M^{\sharp}$ in the neighborhoods of the points of the contraction of the non-cores 
are defined as the open neighborhoods at the origin of 
$\bC^2$ and the above lifted actions (see also \cite[\S3.2]{A2010}).
For other orbifold charts  of $M^{\sharp}$, we could use the open sets of $\wM$ and 
the restricted actions of $G$ on them.
In conclusion, $f^{\sharp}:M^{\sharp} \longrightarrow \Delta$ has the structure 
of an orbifold fibration 
(see Definition \ref{deforbfib}). 

 {\rm (ii)} The marking $\ww$ (\ref{equation0602x7}) trivially decends to the Weyl marking
\begin{equation}
\label{equation0602x8}
w: \Sigma_g(\sigma) \longrightarrow \wF.
\end{equation}
If one compares the marking $\ww$
with $w$, $\ww$ is determined through the lifting $\Sigma_g \rightarrow \Sigma_g(\sigma)$ 
which includes the data of the screw numbers. 
This is the reason why we use $\ww$ instead of $w$.
Note that this type of marking for a stable curve is used 
by Hubbard--Koch \cite[Def.2.1]{HK} 
for the construction of the augmented Teichm\"{u}ller 
space (see also  \cite[p.490]{ACG}).
\end{remark}


\subsection{Local orbifold moduli maps and Kodaira-periodicity}
\label{ppsubsection}

In this subsection, we define the 
{\it local orbifold moduli map}\index{local orbifold moduli map} 
of a degenetaion and show that it has the 
Kodaira-periodicity\index{Kodaira-periodicity}.

For a given degeneration $f: M \rightarrow \Delta$ of genus $g$ 
with topological monodormy $\mu_f \in  \cP_g^{(-)}(\sigma)$, 
the restricted holomorphic 
family $f^{-1}(\Delta^*)  \rightarrow \Delta^*$ 
induces the moduli map $ \Delta^* \rightarrow M_g$.
By the classical stable reduction theorem and the valuative criterion algebraically, 
or by Imayoshi \cite[Th.2,Th.4]{Im81} analytically, 
this map has a holomorphic extension to the 
Deligne--Mumford compactification\index{Deligne--Mumford compactification}
\begin{equation}
\label{equation0603no1}
J_f:\Delta \longrightarrow \overline{M}_g,
\end{equation}
which is called the 
canonically extended  moduli map\index{canonically extended  moduli map} (\cite{Namikawa1974}).
This map is determined by $f^{-1}(\Delta^*)  \rightarrow \Delta^*$, and is 
independent of the choice of the local birational model of $f$.

Now we consider the marked orbifold structure 
$f^{{\rm orb}}= \{ \wf,  f^{\sharp}, G, \ww \}$ of $f$ in Definition \ref{definition0602no1}
and the orbifold 
$\overline{M}_g^{orb} = 
\{ ( D_{\epsilon}(\sigma), W(\sigma), \varphi_{\sigma}, M_{\epsilon}(\sigma)) \}_{\sigma \in \cC_g/ \Gamma_g}$ 
in (\ref{mgbarchars}).
The map $J_f$ is lifted to the map of these spaces as follow.
By the descent of  the marking $\ww$ to $w$ by (\ref{equation0602x8}), 
we consider the $\sigma$-Weyl marked stable curve $[\wF, w]$ ($\wF=\wf^{-1}(0)$).
Since the generator $\wmu$ of $G$ induces an analytic automorphism of $\wF$,
the point $p=[\wF, w]$ is contained in the equisymmetric strata $T_{\sigma}^{\wmu}$ in Theorem \ref{stableequithm}:
\begin{equation}
\label{equation0603no2}
p=[\wF, w] \in T_{\sigma}^{\wmu} \subset T({\sigma}) \subset D_{\epsilon}(\sigma).
\end{equation}
We consider the universal  family 
$\overline{\pi}: \overline{Y}_g^{orb} \longrightarrow  \overline{M}_g^{orb}$ 
of (\ref{uinvdegfam}), (\ref{totalunivchars}).
By Theorem \ref{KuranishioverCont}, there exists 
an open neighborhood $B$ of $p$ in $D_{\epsilon}(\sigma)$ such that 
the  restricted family 
$
\overline{\pi}_X=\pi_{\sigma}|_X :  X=\pi_{\sigma}^{-1}(B)  \subset X_{\epsilon}(\sigma)
\longrightarrow B \subset D_{\epsilon}(\sigma)
$
is a standard Kuranishi family\index{standard Kuranishi family}
 of $\wF$ with  marking $w$.
Since $\wf:\wM \rightarrow \wDelta$ is a local deformation of $\wF$, 
it follows from the universality of the Kuranishi family that 
there uniquely exists a holomorphic map
\begin{equation}
\label{equation0603no3}
\wJ_{\wf} \: : \: \wDelta \longrightarrow B
\end{equation}
such that $\wf$ is the pull back of $\overline{\pi}_X$ by $\wJ_{\wf} $.
Then we have the commutative diagram 

\setlength{\unitlength}{0.40mm}
\begin{picture}(270,71)(7,-3)
\put(80,50){\makebox(6,6){$\wDelta$}}
\put(80,5){\makebox(6,6){$\Delta$}} 
\put(213,49){\makebox(12,6){$\wJ_{\wf} (\wDelta)\: \subset \:\:\: B \:\:\: 
\subset \:\: D_{\epsilon}(\sigma)$}} 
\put(227,4){\makebox(12,6){$J_f(\Delta)\subset \varphi_{\sigma}(B) 
\subset M_{\epsilon}(\sigma) \subset \overline{M}_g$}} 
\put(81,43){\vector(0,-1){24}} 
\put(181,43){\vector(0,-1){24}} 
\put(254,43){\vector(0,-1){24}} 
\put(93,51){\vector(1,0){64}}
\put(93,7){\vector(1,0){65}}
\put(68,26){\makebox(6,6){$\pi$}} 
\put(190,26){\makebox(6,6){$\varphi_{\sigma,\wDelta}$}} 
\put(260,26){\makebox(6,6){$\varphi_{\sigma}$}} 
\put(121,39){\makebox(6,6){$\wJ_{\wf} $}} 
\put(121,-5){\makebox(6,6){$J_f$}} 
\end{picture}
\vspace{-0.2cm}

\begin{center}
(Diagram III) \ The local orbifold moduli map
\end{center}

\noindent where $\varphi_{\sigma,\wDelta}$ is the restriction to 
$\wJ_{\wf} (\wDelta)$ of the orbifold structure map 
$\varphi_{\sigma}:  D_{\epsilon}(\sigma) \longrightarrow 
M_{\epsilon}(\sigma)=D_{\epsilon}(\sigma) / W(\sigma)$.
More precisely, the space $B$ in Diagram III is really 
the Kuranishi space with  Weyl marking ($B, w$) and 
the image of the restriction of the forgetting map 
$(B,w) \rightarrow B \rightarrow \varphi_{\sigma}(B)$ 
is nothing but the quotient of $B$ by $G$; 
$\varphi_{\sigma}(B)=B/G$.
By the construction of $\wf$ in Diagram II in \S \ref{orbifoldmodeldeg}, 
the group $G$ also acts on the analytic subspace $\wJ_{\wf}(\wDelta)$ of $B$ 
such that $J_f(\Delta)=\wJ_{\wf}(\wDelta)/G$.
In this sense, $\wJ_{\wf}$ is the lifting of $J_f$.
Considering these points, we define the following:

\begin{definition}
\label{definition0603no1}
For a degeneration $f: M \rightarrow \Delta$ of Riemann surfaces with  topological 
monodormy $\mu_f \in \cP_g^{(-)}(\sigma)$ and  orbifold structure 
$f^{{\rm orb}}= \{ \wf, f^{\sharp},G, \ww \}$, we define 
the {\em local orbifold moduli map} by 
\begin{equation}
\label{equation0603no3b}
\{ \wJ_{\wf}: \wDelta \rightarrow D_{\epsilon}(\sigma), \:J_f: \Delta \rightarrow M_{\epsilon}(\sigma), \:
G \}
\end{equation}
which satisfies Diagram III.
For simplicity, we sometimes write (\ref{equation0603no3b})  merely by 
$\wJ_{\wf}: \wDelta \rightarrow D_{\epsilon}(\sigma)$.
\end{definition}

In order to describe $\wJ_{\wf}: \wDelta \rightarrow D_{\epsilon}(\sigma)$ explicitly, 
we  will define a special class of holomorphic functions.

\begin{definition}
\label{definition0603no3}
A holomorphic function $\varphi(t)$ of a variable $t$ 
at the origin in $\bC$ is called a 
{\em pseudo-periodic function\index{pseudo-periodic function} with multiplicity $\gamma \in \bN$ and 
period $L \in \bN$} if there exists a holomorphic function 
$\widetilde{\varphi}(t)=\sum_{i=0}^{\infty} c_{i}t^{i}$ 
$(c_i \in \bC)$ with 
\begin{equation}
\label{defpff}
\varphi(t)=t^{\gamma}\widetilde{\varphi}(t^{L})
=c_{0}t^{\gamma}+c_{1}t^{\gamma+L}+c_{2}t^{\gamma+2L}+ \cdots  \:\:\:\:\:\: (c_{0} \not= 0).
\end{equation}
\end{definition}

\begin{remark}
\label{remark_kodaira}
The notion of pseudo-periodic function in this framework  is inspired by 
Kodaira \cite[\S8]{Kodaira}.
As we already explained,  
the functional invariant of a degeneration of elliptic curves defined in \cite[\S7]{Kodaira} is 
the orbifold moduli map in our terminology. 
The functional  invariant 
belongs to the class of pseudo-periodc functions in the present terminology 
around each degenerate fiber germ in  \cite[\S8]{Kodaira}. 
\end{remark}

Pseudo-periodic  functions are  characterized as follows.

\begin{lemma}
\label{lemma0603no1}
Let $\varphi(t)$ be 
a holomorphic function at the origin with $\varphi(0)=0$ and $\varphi(t) \not\equiv 0$.
Then the following conditions {\rm (i)} and {\rm (ii)} are equivalent;

\noindent
{\rm (i)} $\varphi(t)$  is a 
pseudo-periodic function with  multiplicity $\gamma$ and  period $L$,

\noindent
 {\rm (ii)} $\varphi(t)$ admits an action 
$\bC \rightarrow \bC$ given by $t \mapsto \be(1/L)t$ 
with the charactor $\be(\gamma/L)$ such that 
the derivations of $\varphi(t)$ at the origin satisfy the following;
\begin{equation}
\label{equation0603no4}
\varphi \Big(\be \big( \frac{1}{L} \big) t \Big)= \be \big(  \frac{\gamma}{L} \big) \varphi(t), \:\:\:
\varphi^{(\gamma)}(0) \not= 0, \: \varphi^{(j)}(0) = 0 \:\: 
{\rm for} \:\: \forall j < \gamma.
\end{equation}
\end{lemma}

{\it Proof} \quad  
We assume (i).
From the expression (\ref{defpff}) of  $\varphi(t)$, we have 
$$
\varphi \Big(\be \big( \frac{1}{L} \big) t \Big)=c_0 \be \big( \frac{\gamma}{L} \big) t^{\gamma}
+\sum_{i \geq 1} c_i \be \big( \frac{\gamma+iL}{L} \big) t^{\gamma+iL}
= \be \big(  \frac{\gamma}{L} \big) \varphi(t).
$$
Since $c_{0} \not= 0$, the conditions of the derivatives in (ii) are also satisfied.

Conversely we assume (ii).
We decompose $\varphi(t)$ into 
$$
\varphi(t)=\varphi_1(t)+\varphi_2(t), \:\:\:\: {\rm where} \:\: 
\varphi_1(t)=\sum_{i \equiv \gamma \: {\rm mod}\:L} c_it^i, \:\:\:
\varphi_2(t)=\sum_{j \not\equiv \gamma \: {\rm mod}\:L} c_jt^j.
$$
Then we have 
$$
\varphi \Big(\be \big( \frac{1}{L} \big) t \Big)=
\be \big(  \frac{\gamma}{L} \big) \varphi_1(t)
+\sum_{j \not\equiv \gamma \: {\rm mod}\:L} c_j\be \big(  \frac{j}{L} \big)  t^j.
$$
Thus the condition (ii) says that 
$$
c_j\be \big(  \frac{j}{L} \big) =c_j\be \big(  \frac{\gamma}{L} \big), \:\: 
{\rm for}\:\:\forall j \not\equiv \gamma \:\: ({\rm mod}\:L).
$$
This means that $c_j=0$, $\forall j \not\equiv \gamma \:\: ({\rm mod}\:L)$, 
i.e. $\varphi_2(t) \equiv 0$.
Hence we have $\varphi(t)=\varphi_1(t)$. 
By the conditions of the derivatives in (ii), 
the leading term of $\varphi_1(t)$ should be the non-zero constant multiple of 
$t^{\gamma}$.
Thus we obtain (i).
\qed

\begin{definition}
\label{DefKgammadash}
Let $\varphi(t)$ be a pseudo-periodic function with  multiplicity 
$\gamma$ and the period $L$.
Let $\gamma'$ and $\bK$ be  
integers which satify 
$\gamma=\gamma'+\bK L$,  $\bK \geq 0, \:\:0 \leq \gamma' \leq L-1$, 
$\gamma' \equiv \gamma \:\: ({\rm mod}\:L)$.
We define the {\em analytic screw number}\index{analytic screw number}
 of $\varphi(t)$ by 
\begin{equation}
\label{Eqanalyticscrew}
\bs(\varphi)=\frac{\gamma}{L}=\frac{\gamma'}{L}+\bK,
\end{equation}
and call  $\bK$  and $\gamma'/L$  the  {\em integral term} and the {\em fractional term} 
of $\bs(\varphi)$ respectively.
\end{definition}

The geometric meaning of Def.~\ref{DefKgammadash} will be 
clarified by Th.~\ref{theorem0603no1} (iv) and Rem.~\ref{RemGeoAnSc}.

Now we express the map $\wJ_{\wf}: \wDelta \rightarrow D_{\epsilon}(\sigma)$ in (\ref{equation0603no3b}) 
around $p=[\wF,w] \in T_{\sigma}^{\wmu}$
explicitly.
From the discussions in \S \ref{HMcoord},  
the orbit-irreducible decomposition $\wF=\sum_{i=1}^{\bar{r}} \sum_{j=0}^{m(\wF_i)-1}\wmu^j(\wF_i)$  
and the decomposition  $\sum_{i=1}^{\bar{k}} \sum_{j=0}^{m(P_i)-1} \wmu^j(P_i)$ of nodes of $\wF$ 
induce the Harris--Mumford coordinates 
$(z_1^{(0)}, \cdots, z_i^{(j)},  \cdots, 
 z_{\alpha, \beta}^{(\gamma)}, \cdots, z_{\bar{r}, q_{\bar{r}}}^{(m(\wF_{\bar{r}})-1)})$ 
 of $D_{\epsilon}(\sigma)$  around $p$ as in (\ref{equation0505no3}).
 Then the map $\wJ_{\wf}$ is expressed by $3g-3$ holomorphic functions 
 $\varphi_i^{(j)}(u)$ ($1 \leq i \leq \bar{k}, \: 0 \leq j \leq m(P_i)-1$), 
 $\psi_{\alpha, \beta}^{(\gamma)}(u)$ ($1 \leq \alpha \leq \bar{r}, \:1 \leq \beta \leq q_{\alpha}, \:
 0 \leq \gamma \leq m(\wF_{\alpha})-1$)
with $\varphi_i^{(j)}(0)=\psi_{\alpha, \beta}^{(\gamma)}(0)=0$ 
 sending to each of these coordinates as 
\begin{equation}
\label{equation0603no5}
\wJ_{\wf}: \:\:\:z_i^{(j)}=\varphi_i^{(j)}(u), \:\:\: z_{\alpha, \beta}^{(\gamma)}
=\psi_{\alpha, \beta}^{(\gamma)}(u)\:\:\:\:\:
(u \in \wDelta).
\end{equation}
The following theorem is an extension 
of Kodaira's result from the viewpoint of Remark \ref{remark_kodaira}.

\begin{theorem} 
\label{theorem0603no1} 
Let $(\wJ_{\wf}: \wDelta \rightarrow D_{\epsilon}(\sigma), G$)
be the chart of the orbifold moduli map (\ref{equation0603no3b}) of 
the marked orbifold structure (\ref{equation0602x6}),(\ref{equation0602x7}) of 
a degeneration $f: M \rightarrow \Delta$ of genus $g \geq 2$.
Let $\varphi_i^{(j)}(u), \psi_{\alpha, \beta}^{(\gamma)}(u)$ 
be the holomorphic functions in (\ref{equation0603no5}) 
expressing  $\wJ_{\wf}$.
Then 

{\rm (i)} $\varphi_i^{(0)}(u)$ is a pseudo-periodic 
function with  period $N/m(P_i)$ and 
with multiplicity $n(P_i)$.
Here $n(P_i)$ is the Milnor number 
described in (\ref{equation0602x2}) 
(resp. (\ref{equation0602x4})) in the case where $P_i$ is a non-amphidrome 
node (resp. an amphidrome node).

{\rm (ii)} If $\psi_{\alpha, \beta}^{(0)}(u) \not\equiv 0$, then  
$\psi_{\alpha, \beta}^{(0)}(u)$ is a pseudo-periodic 
function with  period $N/m(\wF_{\alpha})$ and 
 multiplicity 
\begin{equation}
\label{equation0603no5b}
\gamma_{\alpha, \beta}:=\left( 
\frac{n(\wF_{\alpha})-\theta_{\alpha, \beta}}{n(\wF_{\alpha})} 
+\bK_{\alpha, \beta}
\right)  \frac{N}{m(\wF_{\alpha})} 
\end{equation}
where $\bK_{\alpha, \beta}$ is a non-negative integer  and 
$\theta_{\alpha, \beta}/n(\wF_{\alpha})$ is given in 
(\ref{equation0505no1}).

{\rm (iii)} We have $\varphi_i^{(j)}(u)=\varphi_i^{(0)}(u)$ for any $j$,
and $\psi_{\alpha, \beta}^{(\gamma)}(u)= \psi_{\alpha, \beta}^{(0)}(u)$ 
for any $\gamma$.

{\rm (iv)} The analytic screw numbers of the pseudo-periodic functions 
$\varphi_i^{(j)}(u)$ and $\psi_{\alpha, \beta}^{(\gamma)}(u)$ are given by
\begin{equation}
\label{EqAnSc}
\bs(\varphi_i^{(j)})=|\bs(C_{P_i})|,\:\:\:
\bs(\psi_{\alpha, \beta}^{(\gamma)})
=
\frac{n(\wF_{\alpha})-\theta_{\alpha, \beta}}{n(\wF_{\alpha})} 
+\bK_{\alpha, \beta}
\end{equation}
for any $j$ and $\gamma$, where $|\bs(C_{P_i})|$ is  the absolute value of the 
screw number of the cut curve $C_{P_i}$  given in  (\ref{screwnonamp}) and (\ref{equation0602x2}),  or
(\ref{screwamp}) and (\ref{equation0602x4}).
\end{theorem}

{\it Proof} \quad   We prove {\rm (i)}. First we assume that $P_i$ is a non-amphidrome node of $\wF$.
We set $L=N/m(P_i)$.
Since $z_i^{(0)}$ is the dual vector of the generator $v_{P_i}$ 
of the torsion sheaf $\tau_{P_i}$,
it follows from Lemma \ref{lemma0505no1}  and (\ref{equation0602x2}) 
that $\wmu^{m(P_i)}$ acts on $z_i^{(0)}$ by 
\begin{equation}
\label{equation0603no6}
\wmu^{m(P_i)}: \: z_i^{(0)} \longmapsto \be \Big( 
\frac{\delta^{(1)}}{\lambda^{(1)}}+\frac{\delta^{(2)}}{\lambda^{(2)}} \Big)z_i^{(0)}
=\be \Big(\frac{n(P_i)}{L} \Big) z_i^{(0)}.
\end{equation}
Since $\wmu(\wDelta)$ is naturally the $G$-invariant subspace in 
$H^0(\wF, \Omega_{\wF}^1 \otimes \omega_{\wF})^* \subset B \subset 
D_{\epsilon}(\sigma)$ and $\wmu^{m(P_i)}$ acts on $\wDelta$ by 
$u \mapsto \be(1/L)u$, 
the function $\varphi_i^{(0)}(u)$ satisfies 
$\varphi_i^{(0)}(\be(1/L)u)=\be\big( n(P_i)/L \big) \varphi_i^{(0)}(u)$.
Therefore by Lemma \ref{lemma0603no1}, $\varphi_i^{(0)}(u)$ is a 
pseudo-periodic function 
with period $L$, and its multiplicity is congruent to 
$n(P_i)$ modulo $L$.

On the other hand, as we explained in \S3.1 for the comment (II) of (\ref{dualfundamentalseq}), the Kuranishi 
family of $\wF$ has a local structure $\{ xy=t | (x,y,t) \in \bold{C}^3, |x|, |y|, |t| < 1 \}$ at a node 
(\cite[pp.184--186]{ACG}), i.e. it is a plumbing variety (cf.\cite{Kra}). Since $\wM$ has singularity 
$xy=u^{n(P_i)}$ at $P_i$ and $\wf : \wM \to \wDelta$ should be pulled back from the Kuranishi family 
by $\wJ_{\wf}$, the map $\wJ_{\wf}$ is locally  given by 
$$\wJ_{\wf} : u \mapsto t= u^{n(P_i)}\psi(u)$$
where $\psi(u)$ is a holomorphic function with $\psi(0)\neq 0$. \par
Therefore the multiplicity of $\varphi_i^{(0)}(u)$ exactly coincides with $n(P_i)$.\par

Hence we obtain the assertion (i) for non-amphidrome case.
The case where $P_i$ is an amphidrome node, the discussion is similar and is omitted.

We prove {\rm (ii)}.
We consider a component $\wF_{\alpha}$ and set $L'=N/m(\wF_{\alpha})$.
Since $ z_{\alpha, \beta}^{(0)}$ is the dual vector of $v_{\alpha, \beta}$ in (\ref{equation0505no1}),
the map $\wmu^{m(\wF_{\alpha})}$ acts by
$$
\wmu^{m(\wF_{\alpha})}: \:  z_{\alpha, \beta}^{(0)} \longmapsto 
\be \left(  -\frac{\theta_{\alpha, \beta}}{n(\wF_{\alpha})} \right)  z_{\alpha, \beta}^{(0)}
=\be \left(  \frac{1}{L'} \left( -\frac{\theta_{\alpha, \beta}L'}{n(\wF_{\alpha})}  \right) \right)  z_{\alpha, \beta}^{(0)}.
$$
Therefore we similarly obtain the assertion (ii) by Lemma  \ref{lemma0603no1}.
 
 Since the iterations of the map $\wmu$ permute cyclically and isomorphically each neighborhood 
 of the orbit of a node and also each neighborhood of the orbit of a component of $\wF$, 
 the assertion (iii) follows. 
 
 The assertion (iv) is clear from (i), (ii), (iii) and Definition \ref{DefKgammadash}.
\qed

 
\begin{remark}
\label{remarkmodulimap}
In the property (ii) of Th.\ref{theorem0603no1}, the trivial function  
$\psi_{\alpha, \beta}^{(0)}(u) \equiv 0$ may occure.
In particular, in the case where $\bar{k}=0$ and 
 $\psi_{1, \beta}^{(0)}(u) \equiv 0$  for all $1 \leq \beta \leq 3g-3$, 
$\wJ_{\wf}$ is the constant map $\wJ_{\wf}(\wDelta)=p$, i.e.
$\wM \simeq \wF \times \wDelta$  and $\wf$ is the projection 
$\wF \times \wDelta \longrightarrow \wDelta$ 
such that $f$ is obtained from the resolution of $(\wF \times \wDelta)/G$.
\end{remark}


\begin{remark}
\label{RemGeoAnSc}
The fractional term 
$(n(\wF_{\alpha})-\theta_{\alpha, \beta})/n(\wF_{\alpha})$ of the analytic screw number 
$\bs(\psi_{\alpha, \beta}^{(\gamma)})$ in (\ref{EqAnSc})  is a topological 
invariant, because it  is determined by the total valency using Prop. \ref{eigenspaceformula}.
On the other hand, the integral term  $\bK_{\alpha, \beta}$ of 
$\bs(\psi_{\alpha, \beta}^{(\gamma)})$ is not a topological invariant, and is 
determined purely by the analytic structure of $f$.
By comparison with the usual screw numbers  (\ref{screwnonamp}) and (\ref{screwamp}), 
the number  $\bK_{\alpha, \beta}$ seems to be an ``analytic analog of the 
 number of integral Dehn twists". See also Kuno's paper \cite[\S4.3]{Kuno}.

In a certain situation, $\bK_{\alpha, \beta}$ 
is related to  {\em the modular invariant}\index{modular invariant} 
(\cite{Reid},\cite{Tan2010}) of the fiber germ $(f, F)$ 
from the viewpoint of \cite{AY}.
This point  will be discussed in a forthcoming paper.
\end{remark}


\subsection{Examples of degenerations and their invariants}
\label{Examples}

Here we give two examples of degenerations of Riemann surfaces 
and show 
their orbifold structures and the properties of their orbifold 
moduli maps.

\begin{example}
\label{example0604no2}
Let $\mu: \Sigma_{8} \longrightarrow \Sigma_{8}$ be the element of $\bP_8^{(-)}(\sigma)$ 
with $\sigma=\langle C_1, \cdots, C_5 \rangle$ and the pseudo-period 
$N=84$ described in Figure III.
The total valencies of 
$\cB= \coprod_{i=1}^{4} \cB_{i}=\Sigma_{8} \setminus \coprod_{j=1}^{5} C_{j}$ 
are 
$(g=3,\bar{g}=0, n=7;  2/7+6/7+{\bf 6/7})$ on $\cB_{1}$, 
$(g=1,\bar{g}=0, n=4; {\bf 1/4}+1/4+{\bf 1/2})$ on $\cB_{2}$, 
$(g=1,\bar{g}=0, n=3; {\bf 2/3}+{\bf 2/3}+{\bf 2/3})$ on $\cB_{3}$ and $\cB_{4}$.
Here the valencies written by boldface are attached to the boundary curves 
and those by  roman are to the multiple points assigned by the cross symbols 
inside the body.

The action on the graph 
has order $2$ with $\mu(\cB_{3})=\cB_{4}$, $\mu(C_4)=C_5$ 
such that $C_1$, $C_4$, $C_5$ are non-amphidrome and 
$C_2$, $C_3$ are amphidrome.
The screw numbers are 
$\bs(C_{1})=-6/7-1/4-\bK_{1}$  $(\bK_{1} \geq -1)$, 
$\bs(C_{2})=-2(2/3+\bK_{2})$ $(\bK_{2} \geq 0)$,
$\bs(C_{3})=-2(2/3+\bK_{3})$ $(\bK_{3} \geq 0)$, 
$\bs(C_{4})=\bs(C_{5})=-1/2-2/3-\bK_{4}$ $(\bK_{4} \geq -1)$.

\begin{figure}[b]

\setlength{\unitlength}{0.70mm}
\begin{picture}(185,60)(-5,-5)

\begingroup
\linethickness{1.8pt}
\cbezier(60,40)(-5,29)(-5,14)(60,2)
\endgroup
\begingroup
\linethickness{1.0pt}
\cbezier(60,40)(54,29)(54,14)(60,2)
\endgroup
\begingroup
\linethickness{1.8pt}
\qbezier(16,20)(20,16)(24,20)
\qbezier(29,20)(33,16)(37,20)
\qbezier(42,20)(46,16)(50,20)
\qbezier(17,19)(20,22)(23,19)
\qbezier(30,19)(33,22)(36,19)
\qbezier(43,19)(46,22)(49,19)
\endgroup
\cbezier(60,40)(64,29)(64,14)(60,2)
\put(30,37){\makebox(6,6)[b]{$\cB_{1}$}} 
\put(57,42){\makebox(6,6)[b]{$C_{1}$}} 
\put(20,8){\vector(3,2){8}}
\put(9,0){\makebox(16,9){$\frac{2}{7}+\frac{6}{7}+{\bf \frac{6}{7}}$}} 
\put(50,-1){\vector(1,2){5}}
\put(39,-8){\makebox(16,9){$-\frac{6}{7}-\frac{1}{4}-\bK_{1}$}} 
\put(40,29){\line(1,-1){2}} 
\put(40,27){\line(1,1){2}}
\put(42,13){\line(1,-1){2}} 
\put(42,11){\line(1,1){2}}

\begingroup
\linethickness{1.8pt}
\qbezier(60,40)(82,44)(104,42)
\qbezier(60,2)(82,-2)(104,0)
\endgroup
\begingroup
\linethickness{1.0pt}
\qbezier(104,42)(98,31)(104,25)
\endgroup
\begingroup
\linethickness{1.8pt}
\cbezier(104,25)(92,23)(92,19)(104,17)
\endgroup
\begingroup
\linethickness{1.0pt}
\qbezier(104,17)(98,11)(104,0)
\endgroup
\qbezier(104,42)(107,31)(104,25)
\qbezier(104,17)(107,11)(104,0)
\begingroup
\linethickness{1.8pt}
\cbezier(104,25)(116,23)(116,19)(104,17)
\qbezier(73,20)(79,16)(85,20)
\qbezier(74,19)(79,22)(84,19)
\endgroup
\put(79,45){\makebox(6,6)[b]{$\cB_{2}$}} 
\put(94,31){\makebox(6,6)[b]{$C_{4}$}} 
\put(94,6){\makebox(6,6)[t]{$C_{5}$}} 
\put(80,-3){\vector(0,1){7}}
\put(72,-12){\makebox(16,9){${\bf \frac{1}{4}}+\frac{1}{4}+{\bf \frac{1}{2}}$}} 
\put(110,-3){\vector(-1,2){4}}
\put(117,-10){\makebox(16,9){$-\frac{1}{2}-\frac{2}{3}-\bK_{4}$}} 
\put(110,44){\vector(-1,-2){4}}
\put(117,43){\makebox(16,9){$-\frac{1}{2}-\frac{2}{3}-\bK_{4}$}} 
\put(78,12){\line(1,-1){2}} 
\put(78,10){\line(1,1){2}}

\begingroup
\linethickness{1.8pt}
\qbezier(104,42)(129,42)(154,32)
\qbezier(154,32)(164,26)(162,21)
\qbezier(162,21)(164,15)(154,10)
\qbezier(104,0)(129,0)(154,10)
\qbezier(126,20)(134,14)(144,20)
\qbezier(126,20)(134,26)(144,20)
\endgroup
\begingroup
\linethickness{1.0pt}
\qbezier(113,20)(119,17)(126,20)
\qbezier(144,20)(153,17)(162,21)
\endgroup
\qbezier(113,20)(119,22)(126,20)
\qbezier(144,20)(153,23)(162,21)
\begingroup
\linethickness{1.8pt}
\qbezier(115,31)(120,28)(126,31)
\qbezier(116,31)(120,33)(125,31)
\qbezier(115,10)(120,7)(126,10)
\qbezier(116,10)(120,12)(125,10)
\endgroup
\put(135,26){\makebox(6,6)[b]{$\cB_{3}$}} 
\put(135,8){\makebox(6,6)[t]{$\cB_{4}$}} 
\put(117,22){\makebox(6,6)[b]{$C_{3}$}} 
\put(150,23){\makebox(6,6)[b]{$C_{2}$}} 
\cbezier(180,26)(188, 22)(188,18)(180,14)
\put(180,14){\vector(-1,-1){0}}
\put(188,18){\makebox(6,6)[l]{$\mu$}} 
\put(141,2){\vector(-6,5){18}}
\put(145,-5){\makebox(16,9){$-2(\frac{2}{3}+\bK_{3})$}} 
\put(142,39){\vector(-6,-5){7}}
\put(146,37){\makebox(16,9){${\bf \frac{2}{3}}+{\bf \frac{2}{3}}+{\bf \frac{2}{3}}$}} 
\put(152,7){\vector(-6,5){6}}
\put(158,2){\makebox(16,9){${\bf \frac{2}{3}}+{\bf \frac{2}{3}}+{\bf \frac{2}{3}}$}} 
\put(164,30){\vector(-6,-5){8}}
\put(170,27){\makebox(16,9){$-2(\frac{2}{3}+\bK_{2})$}} 

\end{picture}

\medskip

\begin{center}
(Figure III) The data of the topological monodromy of Ex.~\ref{example0604no2}
\end{center}

\end{figure}

Let $f:M \longrightarrow \Delta$ be the  normally minimal model of a degeneration with 
topological monodormy  $\mu_{f}=\mu$.
The central  fiber $F=f^{-1}(0)$, i.e. the generalized quotient $\cW(\mu_f)$  
is given as in Figure IV by the algorithm in \cite{MM}.
Here the circles mean $\bP^1$'s 
and the numbers in the circles mean their multiplicities.
(The case of $\bK_{i}=-1$ for some $i$ is omitted in this figure.)
$F$ has three core components of multiplicities $7, 4$ and $6$.

By the contraction $M \longrightarrow M^{\sharp}$  of the non-core 
of  $F$, we obtain the orbifold model 
$f^{\sharp}:M^{\sharp} \longrightarrow \Delta$.
The  fiber $F^{\sharp}=(f^{\sharp})^{-1}(0)$ is described in Figure V.
Here $M^{\sharp}$ has five isolated cyclic quotient singularities and 
two dihedral 
quotient singularities whose supports  are indicated 
by the cross symbols in this figure.

\begin{figure}[t]

\setlength{\unitlength}{0.80mm}
\begin{picture}(180,52)(5,11)

\begingroup
\linethickness{1.0pt}
\put(35,30){\circle{8}}
\endgroup
\put(34,29){\makebox(2,2)[l]{$7$}} 
\multiput(30.5,25.5)(-3.4,-3.4){2}{\circle{5}}
\put(29.7,24.7){\makebox(2,2)[l]{$2$}}
\put(26.3,21.3){\makebox(2,2)[l]{$1$}}
\multiput(30.5,34.5)(-3.4,3.4){6}{\circle{5}}
\put(29.6,34.1){\makebox(2,2)[l]{$6$}}
\put(26.2,37.5){\makebox(2,2)[l]{$5$}}
\put(22.8,40.9){\makebox(2,2)[l]{$4$}}
\put(19.4,44.2){\makebox(2,2)[l]{$3$}}
\put(16.0,47.6){\makebox(2,2)[l]{$2$}}
\put(12.6,51.0){\makebox(2,2)[l]{$1$}}
\multiput(41.5,30)(5,0){6}{\circle{5}}
\put(40.5,29.3){\makebox(2,2)[l]{$6$}}
\put(45.5,29.3){\makebox(2,2)[l]{$5$}}
\put(50.5,29.3){\makebox(2,2)[l]{$4$}}
\put(55.5,29.3){\makebox(2,2)[l]{$3$}}
\put(60.5,29.3){\makebox(2,2)[l]{$2$}}
\put(65.5,29.3){\makebox(2,2)[l]{$1$}}
\multiput(72.1,29.3)(3.5,0){2}{\circle*{1}}
\cbezier(65.7,32.9)(68.5,37.0)(71.5,32.9)(73.7,37.0)
\cbezier(73.7,37.0)(76.4,32.9)(79.4,37.0)(82.2,32.9)
\put(73.7,38.3){\makebox(3,3)[b]{$\bK_{1}$}}
\multiput(81,30)(5,0){2}{\circle{5}}
\put(80.0,29.3){\makebox(2,2)[l]{$1$}}
\put(85.0,29.3){\makebox(2,2)[l]{$1$}}

\begingroup
\linethickness{1.0pt}
\put(92.3,30){\circle{8}}
\endgroup
\put(91.3,29){\makebox(2,2)[l]{$4$}} 
\multiput(98.8,30)(5,0){2}{\circle{5}}
\put(97.8,29.3){\makebox(2,2)[l]{$2$}}
\put(102.8,29.3){\makebox(2,2)[l]{$2$}}
\multiput(109.4,29.3)(3.5,0){2}{\circle*{1}}
\cbezier(102.2,32.9)(104.9,37.0)(107.7,32.9)(110.7,37.0)
\cbezier(110.7,37.0)(113.4,32.9)(114.4,37.0)(119.2,32.9)
\put(110.7,38.3){\makebox(3,3)[b]{$\bK_{4}$}}
\multiput(118.3,30)(5,0){2}{\circle{5}}
\put(117.0,29.3){\makebox(2,2)[l]{$2$}}
\put(122.0,29.3){\makebox(2,2)[l]{$4$}}
\put(92.3,36.5){\circle{5}}
\put(91.0,35.8){\makebox(2,2)[l]{$1$}}

\begingroup
\linethickness{1.0pt}
\put(129.7,30){\circle{8}}
\endgroup
\put(128.7,29){\makebox(2,2)[l]{$6$}} 
\multiput(134.2,25.5)(3.4,-3.4){2}{\circle{5}}
\multiput(134.2,34.5)(3.4,3.4){2}{\circle{5}}
\put(133.1,24.7){\makebox(2,2)[l]{$4$}}
\put(136.5,21.3){\makebox(2,2)[l]{$2$}}
\put(133.1,33.8){\makebox(2,2)[l]{$4$}}
\put(136.5,37.2){\makebox(2,2)[l]{$2$}}
\multiput(140.6,17.9)(2,-2){2}{\circle*{1}}
\cbezier(140.9,24.0)(142.7,21.5)(143.9,21.0)(145.8,22.2)
\cbezier(145.8,22.2)(146.7,19.0)(148.4,16.9)(149.4,15.8)
\put(145.8,22.6){\makebox(3,3)[b]{$\bK_{2}$}}
\put(149.9,9.6){\circle{5}}
\put(146.5,13.0){\circle{5}}
\put(148.9,8.9){\makebox(2,2)[l]{$2$}}
\put(145.5,12.3){\makebox(2,2)[l]{$2$}}
\multiput(140.6,42.1)(2,2){2}{\circle*{1}}
\cbezier(134.5,39.2)(136.8,44.5)(138.8,40.3)(138.9,47.3)
\cbezier(138.9,47.3)(141.4,45.0)(144.3,49.7)(144.6,50.8)
\put(137.1,48.3){\makebox(3,3)[b]{$\bK_{3}$}}
\put(149.9,50.4){\circle{5}}
\put(148.9,49.7){\makebox(2,2)[l]{$2$}}
\put(146.5,47.0){\circle{5}}
\put(145.5,46.3){\makebox(2,2)[l]{$2$}}
\put(153.3,13.0){\circle{5}}
\put(146.5,6.2){\circle{5}}
\put(153.3,47.0){\circle{5}}
\put(146.5,53.8){\circle{5}}
\put(152.0,12.4){\makebox(2,2)[l]{$1$}}
\put(145.2,5.6){\makebox(2,2)[l]{$1$}}
\put(152.0,46.4){\makebox(2,2)[l]{$1$}}
\put(145.2,53.2){\makebox(2,2)[l]{$1$}}

\end{picture}

\begin{center}
(Figure IV) The central  fiber $F$ of the normally minimal model of Ex.~\ref{example0604no2}
\end{center}

\end{figure}

\begin{figure}
\includegraphics[scale=0.6,bb=5 5 500 250]{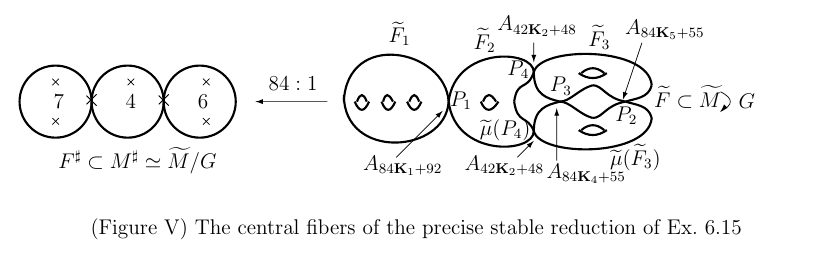}
\end{figure}

Let $\widetilde{\Delta} \longrightarrow \Delta$ 
be the cover $u \mapsto t=u^{84}$, and 
$\wM$ be the normalization of 
$M^{\sharp} \times_{\Delta} \wDelta$.
We obtain the precise stable reduction $\wf:\wM \longrightarrow \wDelta$, 
and the orbifold structure $(f^{\sharp}, \wf, G=\langle \wmu \rangle  \simeq \bZ/84\bZ)$.
The stable  fiber 
$\wF=\wf^{-1}(0)$ is as in Figure V.
The topological monodromy  $\mu_{\wf} \simeq \mu_f^{84}$ 
is trivial on $\cB$ and acts  as $n_i=84 |\bs(C_{i})|$-right Dehn twist at an annular 
neighborhood of $C_j$, i.e.
\begin{equation}
\label{equationexampleno1}
n_1=84\bK_1+93, \:
n_2=84\bK_2+56, \:
n_3=84\bK_3+56, \:
n_4=n_5=42\bK_4+49.
\end{equation}
The marking is defined by the composition 
$w:\Sigma_8 \rightarrow \Sigma_8(\sigma) \simeq \cW(\mu_{\wf}) \rightarrow \wF$.
Let $\wF=\sum_{i=1}^2 \wF_i+\sum_{j=0}^1 \wmu^j(\wF_3)$ be the orbit-irreducible 
decomposition such that 
$\overline{w(\cB_i)}=\wF_i$ ($1 \leq i \leq 3$), $\overline{w(\cB_4)}=\wmu(\wF_3)$.
Let $\bP=\sum_{i=1}^3 P_i+\sum_{j=0}^1 \wmu^(P_4)$ be the decomposion of the set of nodes 
on $\wF$ such that 
 $w(C_i)=P_i$  ($1 \leq i \leq 4$) and $w(C_5)=\wmu(P_4)$.
$\wM$ has the singularity of type $xy=u^{n_i}$ 
for  (\ref{equationexampleno1}) at each node.
The action of $G$ to $\wM$ is not effective.
In fact, the order of ${\rm Stab}_G(\wF_1) =G$ is $86$ while 
the order of the automorphism $\wmu |_{\wF_1}$ of $\wF_1$ is $7$.

We set $V_1=H^0(\wF_1, 2K_{\wF_1}+P_1)$, 
$V_2=H^0(\wF_2, 2K_{\wF_2}+P_1+P_4+\wmu(P_4))$, 
$V_3=H^0(\wF_3, 2K_{\wF_3}+P_2+P_3+P_4)$ and 
$\wmu(V_3)=H^0(\wmu(\wF_3), 2K_{\wF_3}+P_2+P_3+\wmu(P_4))$.
By Example \ref{example0502no1} and the similar argument using Prop.~\ref{eigenspaceformula}, 
the log-quadratic characters  are 
\begin{equation}
\label{equation0604no2}
{\rm Ch}_{\wmu}V_1= \left\{ \frac{1}{7}, \frac{2}{7}, \frac{2}{7}, 
\frac{3}{7}, \frac{4}{7}, \frac{5}{7}, \frac{6}{7} \right\}, \:
{\rm Ch}_{\wmu}V_2=\left\{ \frac{1}{4}, \frac{1}{2}, \frac{3}{4} \right\},\:
{\rm Ch}_{{\wmu}^2}\wmu^j(V_3)=
\left\{ \frac{1}{3}, \frac{1}{3}, \frac{2}{3} \right\}, 
\end{equation}
for $j=0,1$.
The little Teichm\"{u}ller space $T({\sigma})$ is locally 
an open set of $\bigoplus_{i=1}^2 V_i \bigoplus_{j=0}^1 \wmu^j(V_3)$ 
and ${\rm dim}\:T({\sigma})=16$.
The equisymmetric strata $T_{\sigma}^{\wmu}$ 
consists of a unique point $p=[\wF, w]$, since 
the $0$-eigen space  
is $0$-dimensional by (\ref{equation0604no2})
and $T_{\sigma}^{\wmu}$ is connected 
by Cor.~\ref{coro0504no1}.
Let 
\begin{equation}
\label{equation0604no2b}
(z_1, z_2,z_3, z_4^{(0)}, z_4^{(1)}, z_{1,1}, \cdots, z_{1,7},
z_{2,1}, z_{2,2},z_{2,3}, z_{3,1}^{(0)}, z_{3,2}^{(0)},z_{3,3}^{(0)},
z_{3,1}^{(1)},z_{3,2}^{(1)},z_{3,3}^{(1)})
\end{equation}
be the system of Harris--Mumford 
coordinates at $p$ on $D_{\epsilon}(\sigma)$ 
which are ordered as the dual vectors 
of the ordered torsion sheaf at the nodes and the eigenvectors corresponding 
to the characters in (\ref{equation0604no2}).
From (\ref{equation0604no2}), the multiplicities 
given in (\ref{equation0603no5b}) are 
\begin{equation}
\label{equation0604no3}
\gamma_{1,j}=84\bK_{1,j}+12\gamma_{1,j}', \:\:
\gamma_{1,j}'=6,5,5,4,3,2,1 \:(j=1,2,3,4,5,6,7), \: \bK_{1,j} \in \bZ_{\geq 0},
\end{equation}
\vspace{-0.9cm}
\begin{equation}
\label{equation0604no4}
\gamma_{2,j}=84\bK_{1,j}+21\gamma_{2,j}', \:\:
\gamma_{2,j}'=3,2,1 \:(j=1,2,3), \:\bK_{2,j} \in \bZ_{\geq 0},
\:\:\:\:\:\:\:\:\:\:\:\:\:\:\:\:\:\:\:\:\:\:\:\:\:\:\:\:\:\:\:\:\:
\end{equation}
\vspace{-0.75cm}
\begin{equation}
\label{equation0604no5}
\gamma_{3,j}=42\bK_{3,j}+14\gamma_{3,j}', \:\:
\gamma_{3,j}'=2,2,1 \:(j=1,2,3), \:\bK_{3,j} \in \bZ_{\geq 0},
\:\:\:\:\:\:\:\:\:\:\:\:\:\:\:\:\:\:\:\:\:\:\:\:\:\:\:\:\:\:\:\:
\end{equation}
where the notations mean that $\gamma_{1,1}'=6, \gamma_{1,2}'=5, \cdots, 
\gamma_{3,3}'=1$.
From  Th.~\ref{theorem0603no1}, (\ref{equation0604no2b}), 
(\ref{equationexampleno1}) and 
(\ref{equation0604no3})$\sim$(\ref{equation0604no5}),
the the orbifold moduli map 
$\wJ_{\wf}:\wDelta \rightarrow D_{\epsilon}(\sigma)$ $(u \in \Delta)$ 
has the Kodaira-periodicity given by
$$
z_i=\sum_{k=0}^{\infty} c_{i,k}u^{n_i+84k}\:\:\:\:\:
(i=1,2,3, c_{i,0} \not=0, c_{i,k} \in \bC), \:\:\:\:\:\:\:\:\:
$$
\vspace{-0.3cm}
$$
z_4^{(j)}=\sum_{k=0}^{\infty} c_{4,k}u^{n_4+42k}\:\:\:\:\:
(j=0,1, c_{4,0}^{(k)} \not=0, c_{4,k} \in \bC), 
\:\:\:\:\:\:\:\:
$$
\vspace{-0.3cm}
$$
z_{i,j}=\sum_{k=0}^{\infty} c_{i,j,k}u^{\gamma_{i,j}+84k}\:\:\:\:\:
(i=1,2, j=1,2,3, c_{i,j,k} \in \bC), 
$$
\vspace{-0.3cm}
$$
z_{3,i}^{(j)}=\sum_{k=0}^{\infty} c_{3,i,k}u^{\gamma_{3,i}+42k}\:\:\:\:\:
(i=1,2,3, j=0,1, c_{3,i,k} \in \bC).
$$

\end{example}


\begin{example}
\label{Ex0604no2b}
Let $\mu: \Sigma_{5} \longrightarrow \Sigma_{5}$ be the element of $\bP_5^{(-)}(\sigma)$ 
with $\sigma=\langle C_1, \cdots, C_5 \rangle$ and 
$N=15$ as in Figure VI.
The total valencies of 
$\cB= \coprod_{i=0}^{5} \cB_{i}=\Sigma_{5} \setminus \coprod_{j=1}^{5} C_{j}$ 
are 
$(g=0,\bar{g}=0, n=5, 2/5+3/5+{\bf 1})$ on $\cB_{0}$ and 
$(g=1,\bar{g}=0, n=3,2/3+2/3+{\bf 2/3})$ on $\cB_{i}$ $(1 \leq i \leq 5)$.
The action on the graph 
has order $5$ with  permutations 
$(\cB_{1}, \cB_{3}, \cB_{5}, \cB_{2}, \cB_{4})$, 
$(C_{1}, C_{3}, C_{5}, C_{2}, C_{4})$, i.e. 
the $2/5$-turn in the terminology of \cite[p.148]{MM} such that 
$C_{j}$ is non-amphidrome with  the screw 
number 
$\bs(C_{j})=-2/3-\bK$  $(\bK \geq 0)$.


\begin{figure}[t]
\includegraphics[scale=0.6,bb=5 5 500 250]{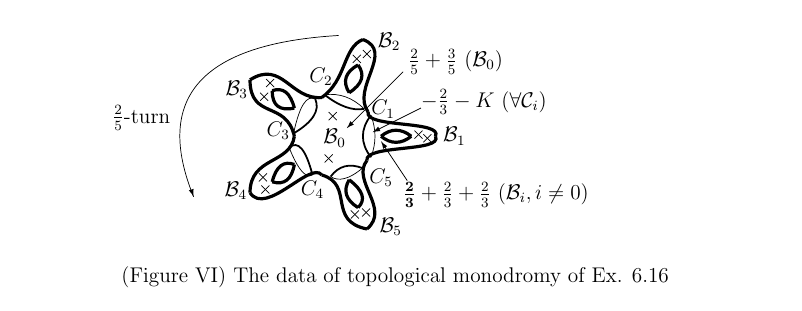}
\end{figure}

Let $f:S \longrightarrow \Delta$ be the degeneration with 
topological monodormy  $\mu_{f}=\mu$.
The  fiber $F=f^{-1}(0)$ and 
the process to obtain its precise stable reduction 
$\wf:\wS \longrightarrow \wDelta$ is shown in Figure VII.
Here 
$\widetilde{F}=\wf^{-1}(0)=\sum_{i=0}^5 \widetilde{F}^{(i)}$ 
is the irreducible decomposition 
such that $\widetilde{F}^{(i)}=\widetilde{\mu}^{i-1}(\widetilde{F}^{(1)})$ 
$(1 \leq i \leq 5)$ for the induced automorphism
$\widetilde{\mu}: \wF \longrightarrow \wF$.
Then $g(\widetilde{F}^{(0)})=0,\: g(\widetilde{F}^{(i)})=1$ and  $\wS$ has 
$A_{3\bK+1}$-singularity at the node $P_i$ 
($1 \leq i \leq 5$). 


\begin{figure}[h]

\setlength{\unitlength}{0.80mm}
\begin{picture}(180,50)(20,10)

\begingroup
\linethickness{1.0pt}
\put(34,30){\circle{10}}
\endgroup
\put(31.7,28.8){\makebox(2.5,2.5)[l]{$15$}} 
\multiput(27.8,24.9)(-4.2,-4,2){2}{\circle{6}}
\put(25.5,23.9){\makebox(2.5,2.5)[l]{$10$}}
\put(22.5,19.8){\makebox(2.5,2.5)[l]{$5$}}
\multiput(27.8,35.1)(-4.2,4.2){2}{\circle{6}}
\put(25.5,34.1){\makebox(2.5,2.5)[l]{$10$}}
\put(22.5,38.2){\makebox(2.5,2.5)[l]{$5$}}
\multiput(42.0,30)(6,0){2}{\circle{6}}
\put(39.8,29.1){\makebox(2.5,2.5)[l]{$10$}}
\put(46.6,29.3){\makebox(2,2)[l]{$5$}}
\multiput(53.5,29.3)(3.5,0){2}{\circle*{1}}
\put(63.5,30.0){\circle{6}}
\put(62.5,29.1){\makebox(2.5,2.5)[l]{$5$}}
\cbezier(46.6,33.5)(48.5,37.5)(52.5,33.4)(54.8,37.5)
\cbezier(54.8,37.5)(56.7,33.4)(60.2,37.5)(63.2,33.7)
\put(54.0,38.6){\makebox(3,3)[b]{$K$}}
\begingroup
\linethickness{1.0pt}
\put(71.5,30){\circle{10}}
\endgroup
\put(70.6,29.1){\makebox(2.5,2.5)[l]{$5$}} 
\multiput(77.7,24.9)(4.2,-4,2){2}{\circle{6}}
\put(76.7,23.9){\makebox(2.5,2.5)[l]{$3$}}
\put(80.8,19.8){\makebox(2.5,2.5)[l]{$1$}}
\multiput(77.7,35.1)(4.2,4.2){2}{\circle{6}}
\put(76.7,34.1){\makebox(2.5,2.5)[l]{$2$}}
\put(80.8,38.2){\makebox(2.5,2.5)[l]{$1$}}

\put(85,30){\vector(1,0){16}}
\put(88.0,29.8){\makebox(9,9){$contr.$}} 

\begingroup
\linethickness{1.0pt}
\multiput(108,30)(10,0){2}{\circle{10}}
\endgroup
\put(105.6,28.6){\makebox(2.5,2.5)[l]{$15$}}
\put(117.1,28.6){\makebox(2.5,2.5)[l]{$5$}}
\put(106.7,34.2){\line(1,-1){2}} 
\put(106.7,32.2){\line(1,1){2}}
\put(106.8,27.8){\line(1,-1){2}} 
\put(106.8,25.8){\line(1,1){2}}
\put(117.1,34.2){\line(1,-1){2}} 
\put(117.1,32.2){\line(1,1){2}}
\put(117.2,27.8){\line(1,-1){2}} 
\put(117.2,25.8){\line(1,1){2}}
\put(111.4,31.5){\line(1,-1){3}} 
\put(111.4,28.7){\line(1,1){3}}

\put(141,30){\vector(-1,0){16}}
\put(129.2,29.8){\makebox(9,9){$15:1$}} 

\begingroup
\linethickness{1.0pt}
\put(168,30){\circle{16}}
\put(181,30){\circle{10}}
\qbezier(178.5,30)(181.5, 26.8)(183.5,30)
\qbezier(179.0,30)(180.8, 32.8)(183.0,30)
\put(172.1,42.2){\circle{10}}
\qbezier(171.3, 40.0)(175.5, 42.1)(172.9, 44.4)
\qbezier(171.5, 40.0)(169.0, 43.1)(172.7, 44.4)
\put(172.1,17.8){\circle{10}}
\qbezier(172.9, 15.6)(168.5, 17.5)(171.3, 20.0)
\qbezier(172.7, 15.6)(175.2, 17.5)(171.5, 20.0)
\put(157,37.5){\circle{10}}
\qbezier(155, 39.0)(160.0, 39.5)(159, 36.0)
\qbezier(155.2, 39.0)(153.5, 36.5)(158.8, 36.0)
\put(157,22.5){\circle{10}}
\qbezier(155, 21.0)(155.7, 26.1)(159, 24.0)
\qbezier(155.2, 21.0)(157.6, 18.7)(158.8, 24.0)
\endgroup

\put(171.8,28.3){\makebox(3,3){$P_{1}$}}
\put(168.3,33.0){\makebox(3,3){$P_{2}$}}
\put(168.6,23.9){\makebox(3,3){$P_{5}$}}
\put(162.8,30.8){\makebox(3,3){$P_{3}$}}
\put(162.8,25.0){\makebox(3,3){$P_{4}$}}

\put(190.3,29.5){\makebox(3,3){$\wF^{(1)}$}}
\put(181.0,42.7){\makebox(3,3){$\wF^{(2)}$}}
\put(181.0,14.3){\makebox(3,3){$\wF^{(5)}$}}
\put(159.5,44.5){\makebox(3,3){$\wF^{(3)}$}}
\put(158.8,12.7){\makebox(3,3){$\wF^{(4)}$}}

\put(186,20.9){\vector(-5,2){11.0}}
\put(190.0,20.1){\makebox(3,3){$\wF^{(0)}$}}
\put(188.5,37.3){\vector(-5,-3){11.8}}
\put(190.5,37.0){\makebox(18,9)[b]{$A_{3\bK+1} (\forall P_{i}$)}} 
\qbezier(170.5,52.5)(139.5,45)(149.5,25)
\put(149.5,25){\vector(1,-3){0}}
\put(138,46){\makebox(24,9){$\frac{2}{5}$-turn}} 

\put(52,12){\makebox(6,12)[b]{$F \subset S$}}
\put(110,12){\makebox(6,12)[b]{$F^{\sharp} \subset S^{\sharp}$}}
\put(198.5,12){\makebox(6,12)[b]{$\wF \subset \wS$}}

\end{picture}

\begin{center}
(Figure VII) The central fibers of the precise stable reduction  of Ex.~\ref{Ex0604no2b}
\end{center}

\end{figure}

We set $V_0=H^0(\wF_0, 2K_{\wF_0}+ \sum_{i=1}^5 P_i)$, 
$V_1=H^0(\wF_1, 2K_{\wF_1}+P_1)$. 
A calculation by  using Prop.~\ref{eigenspaceformula} shows that
the log-quadratic characters  are 
\begin{equation}
\label{equation0604ex2no1}
{\rm Ch}_{\wmu}V_0= \left\{ \frac{1}{5}, \frac{4}{5} \right\}, \:
{\rm Ch}_{\wmu}V_1=\left\{ \frac{2}{3} \right\}.
\end{equation}
Then  $T({\sigma})$ is locally 
an open set of $V_0 \bigoplus_{j=0}^4 \wmu^j(V_1) \simeq \bC^7$,  
and  $T_{\sigma}^{\wmu}$ 
consists of a unique point $p=[\wF, w]$   
by (\ref{equation0604ex2no1}).
From (\ref{equation0604ex2no1}), the multiplicities 
given in (\ref{equation0603no5b}) are 
\begin{equation*}
\gamma_{0,1}=15\bK_{0,1}+12, \:\:
\gamma_{0,2}=15\bK_{0,2}+3, \:\:
\gamma_{1,1}=3\bK_{1,1}+1 \:\:\:
(\bK_{0,1}, \bK_{0,2}, \bK_{1,1} \in \bZ_{\geq 0}).
\end{equation*}
Let 
$(z_1^{(0)},z_1^{(1)},z_1^{(2)}, z_1^{(3)}, z_1^{(4)},  
z_{0,1}, z_{0,2},
z_{1,1}^{(0)}, z_{1,1}^{(1)}, z_{1,1}^{(2)}, z_{1,1}^{(3)}, z_{1,1}^{(4)})$ 
be the system of Harris--Mumford 
coordinates at $p$ on $D_{\epsilon}(\sigma)$.
Then the orbifold moduli map 
$\wJ_{\wf}:\wDelta \rightarrow D_{\epsilon}(\sigma)$ $(u \in \Delta)$ 
has the Kodaira-periodicity as
$$
z_1^{(j)}=\sum_{k=0}^{\infty} c_{1,k}u^{2+3\bK+3k}\:\:\:
(0 \leq j \leq 4, c_{1,0} \not=0, c_{1,k} \in \bC), \:\:\:
\:\:\:\:\:\:\:\:\:\:\:\:\:\:\:\:\:\:\:\:\:\:\:\:\:\:\:\:\:\:\:\:\:\:\:\:\:\:\:\:\:\:\:\:\:\:\:\:\:\:\:\:\:\:\:\:\:\:\:\:\:\:\:\:\:\:\:\:
$$
\vspace{-0.3cm}
$$
z_{0,j}=\sum_{k=0}^{\infty} c_{0,j,k}u^{\gamma_{0,j}+15k}\:\:
(j=1,2, c_{0,j,k} \in \bC), \:\:
z_{1,1}^{(j)}=\sum_{k=0}^{\infty} c_{1,1,k}u^{\gamma_{1,1}+3k}\:\:
(0 \leq j \leq 4, c_{1,1,k}\in \bC).
$$

\end{example}


\section{Recovery of fibered complex surfaces from the universal degenerating family}
\label{secrecovery}

The goal of this section is to show that 
any fibered complex surface 
can be pulled back from the universal degenerating family of Riemann surfaces 
$\overline{\pi}: \overline{Y}_g^{orb} \longrightarrow  \overline{M}_g^{orb}$ 
given in \S4. 
By globalizing the notions discussed in \S6, 
we define the global monodromy $\bmu$ 
and the global orbifold moduli map $\bJ^{orb}$ for a 
fibered complex surface 
${\bf f}:\bM \rightarrow \bB$ of genus $g \geq 2$.
Conversely, starting from the invariants ($\bmu, \bJ^{orb}$), 
we construct a fibered complex surface $\bff$ 
realizing these invariants, 
by pulling back the universal family $\overline{\pi}$ via $\bJ^{orb}$. 

Note that Kodaira \cite{Kodaira} constructed an elliptic surface 
with  given ($\bmu, \bJ^{orb}$) (the homological invariant and the 
functional invariant in his terminology), and called it 
{\it the basic member of the elliptic surface} (\cite[p.603]{Kodaira}).
Our Theorem \ref{Th0703no1} may be considered as the construction of 
{\it the basic members of the fibered complex surfaces of genus $g \geq 2$}.

In \S7.1, our discussion is from the local point of view.
For a given pseudo-periodic map $\mu_{loc} \in \bP_g^{-}(\sigma)$, 
we define the set 
$\bP^*(\wDelta, \Delta;  D_{\epsilon}(\sigma), M_{\epsilon}(\sigma); 
T_{\sigma}^{\mu_{\sigma}})$ 
of {\it the dual pseudo-periodic maps} of $\mu_{loc}$, 
which consists of the orbifold maps $\bJ_{loc}^{orb}$ 
from the orbifold disk to the  chart  $D_{\epsilon}(\sigma)$ 
of $\overline{M}_g^{orb}$ at the equisymmetric strata 
$T_{\sigma}^{\mu_{\sigma}}$ which have the Kodaira-periodicity 
given in Theorem \ref{theorem0603no1}.
Concisely speaking, for a given ($\mu_{loc}, \bJ_{loc}^{orb}$), 
we construct the degeneration $f: M \rightarrow \Delta$ 
which realizes these invariants,  
by  pulling back in the orbifold theoretic sense via $J_{loc}^{orb}$  the 
local structure of the universal family $\overline{\pi}$.
Moreover,  we show that 
$\bP^*(\wDelta, \Delta;  D_{\epsilon}(\sigma), M_{\epsilon}(\sigma); 
T_{\sigma}^{\mu_{\sigma}})$  is the 
classifying space of analytic structures over the fixed topological structure of $f$.

In \S7.2, 
we define the global monodormy and the global orbifold moduli map 
of a fibered complex surface ${\bf f}:\bM \rightarrow \bB$.
Over the restricted holomorphic family 
${\bf f}^{(0)}:\bM^{(0)} \rightarrow \bB^{(0)}$ outsides the 
discriminant locus of $\bff$, 
these notions are the usual ones, 
namely, the monodromy representation to the mapping class group 
$\bmu: \pi_1(\bB^{(0)}, b_0) \rightarrow \Gamma_g$ and 
the pair $(\bJ^{(0)})^{orb}=(\wbJ^{(0)}, \bJ^{(0)})$ 
consisting of the moduli map $\bJ^{(0)}: \bB^{(0)} \rightarrow M_g$ 
and its lift $\wbJ^{(0)}: \wbB^{(0)} \rightarrow T_g$ from the 
universal cover of $\bB^{(0)}$ to the Teichm\H{u}ller space.
We prove that ($\bmu, (\bJ^{(0)})^{orb}$) and 
($\mu_{loc}, \bJ_{loc}^{orb}$)'s  around the critical set given in \S7.1 
are well-patched globally.

In \S7.3, we achieve our main purpose by proving Theorem \ref{Th0703no1}.


\subsection{Local recovery of degenerations from the universal family}
\label{subsec7.1}

We discuss the recovery of a degeneration\index{recovery of degeneration}
 $f: M \rightarrow \Delta$ from the local structure of 
$\overline{\pi}$.

We fix an element $\mu \in \bP_g^{-}(\sigma)$ with 
pseudo-period $N$ (see (\ref{ppmapsymbol})).
We set $G=\bZ/N\bZ$, and let 
($\wDelta, G, \pi_{\wDelta},\Delta$) be the orbifold disk defined by 
$\pi_{\wDelta}:\wDelta \ni u \mapsto t=u^N \in \Delta$.
Let $\mu_{\sigma}: \Sigma_g(\sigma) \rightarrow \Sigma_g(\sigma)$ 
be the analytic automorphism which is a descent of $\mu$ 
(see Prop.~\ref{autoppmapequiv}). 
We consider the equisymmetric strata 
$T_{\sigma}^{\mu_{\sigma}} \subset D_{\epsilon}(\sigma)$. 
We  fix  a point $p=[S, w] \in T_{\sigma}^{\mu_{\sigma}}$, 
and let
($\cdots, z_i^{(j)}, \cdots, z_{\alpha, \beta}^{(\gamma)},\cdots $) 
be the system of Harris--Mumford coordinates\index{Harris--Mumford coordinates} 
at $p$ on $D_{\epsilon}(\sigma)$.

\begin{definition}
\label{Def0701no1}
An orbifold map\index{orbifold map} 
$$
(\wJ,J,G): (\wDelta, G, \pi_{\wDelta},\Delta) 
\longrightarrow 
(D_{\epsilon}(\sigma), W(\sigma), \varphi_{\sigma}, M_{\epsilon}(\sigma))
$$
to the chart (\ref{mgbarchars}) of $\overline{M}_g^{orb}$ is 
said to be  {\em a dual  pseudo-periodic map}\index{dual  pseudo-periodic map}
 of $\mu$ at $p$ 
if the following conditions are satisfied: 
A holomorphic map $\wJ: \wDelta  \rightarrow D_{\epsilon}(\sigma)$  with 
$\wJ(0)=p$ is
expressed  by ($3g-3$) pseudo-periodic  functions 
$z_i^{(j)}=\varphi_i^{(j)}(u)$, $z_{\alpha, \beta}^{(\gamma)}
=\psi_{\alpha, \beta}^{(\gamma)}(u)$
such that 

\noindent 
 {\rm (d-i)}  $\varphi_i^{(0)}(u)$ has 
period $N/m(P_i)$ and 
multiplicity $n(P_i)$  as in (\ref{equation0602x2}) 
or  (\ref{equation0602x4}), 

\noindent 
 {\rm (d-ii)}  $\psi_{\alpha, \beta}^{(0)}(u)$ has 
period $N/m(\wF_{\alpha})$ and 
multiplicity $\gamma_{\alpha, \beta}$ as in (\ref{equation0603no5b}), 
or  $\psi_{\alpha, \beta}^{(0)}(u) \equiv 0$, 

\noindent 
 {\rm (d-iii)} 
  $\varphi_i^{(j)}(u)=\varphi_i^{(0)}(u)$ for $j$,  $0 \leq j \leq m(P_i)-1$, and 
 $\psi_{\alpha, \beta}^{(\gamma)}(u)=
\psi_{\alpha, \beta}^{(0)}(u)$ for $\gamma$, $0 \leq \gamma \leq m(\wF_{\alpha})-1$.

The holomorphic map $J=\varphi_{\sigma} \circ \wJ \circ (\pi_{\wDelta})^{-1}:
\Delta \longrightarrow M_{\epsilon}(\sigma)$ is well-defined.
\end{definition}

The motivation of this definition comes from  Th.~\ref{theorem0603no1}.
All the invariants in  {\rm (d-i)} $\sim$  {\rm (d-iii)} except for  $\bK_{\alpha, \beta}$
 are numerically determined from the data (a),(b),(c) of $\mu$ 
 in Th.~\ref{nmm}.
 But  the $\bK_{\alpha, \beta}$'s 
can be any non-negative integers provided that 
 $\gamma_{\alpha, \beta}  > 0$ (see (\ref{equation0603no5b})).
 Since there are infinitely many choices of such $\bK_{\alpha, \beta}$'s (cf. Th.6.12 (ii)),  
 there are infinitely many choices of ($\wJ,J,G$) for  a given $\mu$ and $p$.
Let $\bP^*(\wDelta, \Delta; D_{\epsilon}(\sigma), M_{\epsilon}(\sigma); p)$ 
be the set of the dual pseudo-periodic maps  ($\wJ,J,G$) of $\mu$ at $p$.
By varying $p$ over $T_{\sigma}^{\mu_{\sigma}}$, we have the following definition: 

\begin{definition}
\label{Def0701no2}
The set of  dual pseudo-periodic maps of  $\mu$ for $T_{\sigma}^{\mu_{\sigma}}$ is defined by
$$
\bP^*(\wDelta, \Delta; D_{\epsilon}(\sigma), M_{\epsilon}(\sigma); T_{\sigma}^{\mu_{\sigma}})
= \bigcup_{p \in T_{\sigma}^{\mu_{\sigma}}} 
\bP^*(\wDelta, \Delta; D_{\epsilon}(\sigma), M_{\epsilon}(\sigma); p).
$$
\end{definition}

\begin{theorem}
\label{theorem0603no2}
{\rm (i)}  For any $\mu \in \bP_g^{-}(\sigma)$ and 
$(\wJ,J,G) \in 
\bP^*(\wDelta, \Delta; D_{\epsilon}(\sigma), M_{\epsilon}(\sigma); T_{\sigma}^{\mu_{\sigma}})$, 
there  uniquely exists a degeneration $f: M \rightarrow \Delta$ 
of Riemann surfaces of genus $g$ such that 
the marked topological monodormy of $f$ coincides with $\mu$, and 
the chart map of the orbifold moduli map\index{orbifold moduli map}
 of $f$ coincides with $(\wJ,J,G)$.

{\rm (ii)} 
For elements 
$(\wJ^{(i)},J^{(i)},G) \in 
\bP^*(\wDelta, \Delta; D_{\epsilon}(\sigma), M_{\epsilon}(\sigma); T_{\sigma}^{\mu_{\sigma}})$ 
$(i=1,2)$, 
let $f^{(i)}: M^{(i)} \rightarrow \Delta$ be the degenerations 
given in {\rm (i)} .
Then $f^{(1)}$ is analytically equivalent to $f^{(2)}$ if and 
only if $(\wJ^{(1)},J^{(1)},G)$ coincides with $(\wJ^{(2)},J^{(2)},G)$.

\end{theorem}

{\it Proof}  \quad We prove {\rm (i)}.
We consider an element 
$(\wJ,J,G) \in 
\bP^*(\wDelta, \Delta; D_{\epsilon}(\sigma), M_{\epsilon}(\sigma); p)$
for $p=[S, w] \in  T_{\sigma}^{\mu}$, 
Harris-Mumford coordinates 
($\cdots, z_i^{(j)}, \cdots, z_{\alpha, \beta}^{(\gamma)},\cdots $) 
of $p$ on $ D_{\epsilon}(\sigma)$, and the expression of $\wJ$ via the 
pseudo-periodic functions $\varphi_i^{(j)}(u)$, and 
$\psi_{\alpha, \beta}^{(\gamma)}(u)$.

Let $\pi_{\sigma}: X_{\epsilon}(\sigma) \rightarrow D_{\epsilon}(\sigma)$ be 
the family in Th.~\ref{KuranishioverCont}, and let 
$\wf:\wM \rightarrow \wDelta$ 
be the pulled back family (in the sense of Def. 4.10) from $\pi_{\sigma}$ by  the map 
$\wJ: \wDelta \rightarrow D_{\epsilon}(\sigma)$, i.e.
$\wf$ is the second projection of the fiber product 
$\wM=\pi_{\sigma}^{-1}(\wJ(\wDelta)) \times_{\wJ(\wDelta)} \wDelta \longrightarrow \wDelta$.
Then the central fiber  $\wf^{-1}(0)$ is isomorphic to $S$ by  construction. 

Here we prove the following claim:

\bigskip

\noindent {\bf Claim} 
$\wM$ is normal with $A$-type singularities and any fiber ${\wf}^{-1}(u)$ 
for $u \not= 0$ is non-singular.

\bigskip

In fact, since $z_i^{(j)}$'s  are the smoothing coordinates at the nodes and 
the $\varphi_i^{(j)}(u)$'s  
are not identically zero, all the nodes of $S$ vanish via the deformation $\wf$ of $S$, 
i.e. the fiber $\wf^{-1}(u)$ for $u \not=0$ is non-singular.
Next, the total space of the restricted family   $\pi_{\sigma}^{-1}(\wJ(\wDelta)) \longrightarrow \wJ(\wDelta)$
has non-normal singularities along the central fiber 
which is the non-isolated cusp of the complex space curve locally given 
in $\bC^{3g-2}$ by the image of the map 
\begin{equation}
\label{equation0603no8}
\wJ:\wDelta \in u \longmapsto (\cdots, z_i^{(j)}, \cdots,  z_{\alpha, \beta}^{(\gamma)},\cdots, x)
= (\cdots, u^{n(P_i)}, \cdots,  u^{\gamma_{\alpha, \beta}},\cdots, x)
\end{equation}
modulo units under the conditions $n(P_i) \geq 2$ and $\gamma_{\alpha, \beta} \geq 2$ for any 
$i, \alpha,\beta$, 
where $x$ is a local fiber coordinate of $\pi_{\sigma}$ at a nonsingular point 
of $\pi_{\sigma}^{-1}(0)$.
Fortunately, these codimension-one singularities 
are resolved along the open locus of the central fiber  automatically 
by the pull back by $\wJ$, 
for  $u$ may be considered via  (\ref{equation0603no8}) as the local uniformization parameter of this singularity.
On the other hand, the germ $(\widetilde{Y}_g, P_i^{(j)})$ at the node $P_i^{(j)}$ in 
$\widetilde{Y}_g$ coincides with the germ 
at the origin given by $xy=z_i^{(j)}$ in the ambient coordinates $( \cdots, z_i^{(j)}, \cdots, z_{\alpha, \beta}^{(\gamma)}, 
\cdots, x, y)$ of $\bC^{3g-1}$ (cf. \cite[Chap.XI, \S3]{ACG}). 
Then the germ $(\wM, P_i^{(j)})$ is the $A$-type singularity $xy=u^{n(P_i)}$. 
The above claim is proved.

\bigskip

Now by the definitions (d-i) $\sim$ (d-iii), Lemma \ref{lemma0603no1} and (\ref{equation0505no2}), 
the complex curve $\wJ(\wDelta) \subset B \subset D_{\epsilon}(\sigma)$ 
is invariant under the action of $G= \langle \mu \rangle$ on $B$ compatible 
with respect to 
the action on $\wDelta$, i.e.
$$
\wJ \left( \be \left( \frac{1}{N} \right) u \right)=\wJ \left( \mu(u) \right).
$$
Therefore $G$ acts relatively on the family $\wf$, which is a natural extension 
of the automorphism of the central fiber $S$.
The explicit descriptions of this action on $\wf$ are essentially the same as 
(\ref{eq0602x00}) $\sim$ (\ref{equation0602x5}) given for Diagram II.
Let $f^{\sharp}:M^{\sharp}=\wM/G \longrightarrow \Delta=\wDelta/G$ 
be the quotient family.
By the composition of $f^{\sharp}$ and the minimal resolution 
$M \rightarrow M^{\sharp}$ of the singularities on $M^{\sharp}$, 
we obtain the normally minimal model $f: M \rightarrow \Delta$ 
of a degeneration.

The marked topological monodormy of $f$ coincides with 
the preassigned $\mu \in \bP_g^{-}(\sigma)$, 
since it is determined by the action of $G$ as in \cite[Th.3.1.1]{A2010}.
The orbifold structure of $f$ is  nothing but the above 
$\{ f^{\sharp}, \wf, G \}$.
Let $J_f: \Delta \rightarrow M_{\sigma}$ be the canonically 
extended moduli map of $f$. Then $\wJ$ is clearly the lift of 
$J_f$ by $G$, i.e. $\wJ$ is the chart map of the orbifold moduli map 
of $f$.
Thus we obtain the desired unique degeneration $f$.

We prove {\rm (ii)}.
Assume that $f^{(i)}$ are analytically equivalent to each other for $i=1,2$.
Then 
$J_{f^{(i)}}: \Delta \longrightarrow M_{\sigma} \subset \overline{M}_g$ coincides with 
each other (cf. \cite{Namikawa1974}).
Therefore there exists an element $g \in G$ such that 
$\wJ_{\wf^{(2)}}= g \circ \wJ_{\wf^{(1)}}$.
If $g \not= {\rm id}$, 
then $\mu_{f^{(1)}} \not= \mu_{f^{(2)}}$, 
contradicting  the assumption. 
Hence $g={\rm id}$ and $(\wJ^{(1)},J^{(1)},G)=(\wJ^{(2)},J^{(2)},G)$.
The converse is obvious.
\qed

\begin{coro}
\label{Cor0701no1}
Let $f: M \rightarrow \Delta$ be a degeneration of genus $g$ with a 
marking $w:\Sigma_g \rightarrow f^{-1}(t_0)$ $(t_0 \in \partial \Delta)$, and 
$\mu_f  \in \bP_g^{-}(\sigma)$ be the marked topological monodromy of 
$(f, w)$.

Let $(\wJ_{\wf},J_f,G) : (\wDelta, G, \pi_{\wDelta},\Delta) 
\longrightarrow 
(D_{\epsilon}(\sigma), W(\sigma), \varphi_{\sigma}, M_{\epsilon}(\sigma))$ 
be the orbifold moduli map of $(f, w)$, 
and  $\pi_{\sigma}: X_{\epsilon}(\sigma) \rightarrow D_{\epsilon}(\sigma)$ be 
the family in Theorem~\ref{KuranishioverCont}.
Then the orbifold model of $f$ is isomorphic to the orbifold pull back 
(in Def.\ref{definition0403no3}) of $\pi_{\sigma}$ 
via $(\wJ_{\wf},J_f,G)$.
\end{coro}

\begin{coro}
\label{Cor0701no2}
Let ${\bf AS}(\mu)$ be the set of (complex) analytic structures of degenerations 
$f: M \rightarrow \Delta$ of genus $g$ under a fixed marking $w$ 
whose topological monodromies $\widehat{\mu}_f$ 
coincide with  a fixed $\mu \in \bP_g^{-}(\sigma)$.
Then ${\bf AS}(\mu)$ is in a bijective correspondence with 
$\bP^*(\wDelta, \Delta; D_{\epsilon}(\sigma), M_{\epsilon}(\sigma); T_{\sigma}^{\mu_{\sigma}})$. 
\end{coro}


\subsection{Global orbifold moduli maps for fibered complex surfaces}
\label{subsec7.2}

Here we define the notion of 
orbifold moduli map for a fibered complex surface.

Let ${\bf f}:\bM \rightarrow \bB$ be a proper surjective holomorphic map from 
a 2-dimensional complex manifold $\bM$ to a compact Riemann surface  $\bB$ 
of genus $h \geq 0$.
Let ${\rm Disc}_{{\bf f}}(\bB)=\{ Q_1, \cdots, Q_s \} \subset \bB$ 
be the discriminant locus, and set 
$\bB^{(0)}=\bB \setminus {\rm Disc}_{{\bf f}}(\bB)$.
We call ${\bf f}$ a 
{\it fibered complex surface of genus}\index{fibered complex surface}
 $g \geq 2$
if any fiber of $\bff$ over  $\bB^{(0)}$ 
is a Riemann surface of genus $g$.
Each degenerate fiber $F_i={\bf f}^{-1}(Q_i)$ 
is assumed to be normally minimal in the sense of \S \ref{orbifoldmodeldeg}.
Let  
${\bf f}^{(0)}:\bM^{(0)} \rightarrow \bB^{(0)}$ 
be the restricted holomorphic family of ${\bf f}$ over $\bB^{(0)}$.
We fix a point $b_0 \in \bB^{(0)}$, and consider the fiber 
$\bff^{-1}(b_0):=\Sigma_g$ as the base Riemmann surface of the marking.

The fundamental group\index{fundamental group} 
$\pi_1(\bB^{(0)}, b_0)$ is generated by 
the loops 
$\beta_1, \cdots, \beta_{2h}, \alpha_1, \cdots, \alpha_s$ 
on $\bB^{(0)}$ starting and ending at $b_0$ with the unique relation 
$$
\beta_1 \beta_2 \beta_1^{-1} \beta_2^{-1} \beta_3 \beta_4 \cdots 
 \beta_{2h-1}^{-1} \beta_{2h}^{-1} \alpha_1 \cdots \alpha_s=1.
$$
Here each $\alpha_i$ is a loop which goes around $Q_i$ once counterclockwise
and $\beta_1, \cdots, \beta_{2h}$ are canonical generators 
of $\pi_1(\bB, b_0)$.
Since ${\bf f}^{(0)}$ is a differentiable fiber bundle, we have 
{\it the monodromy representation} to the mapping class group
\begin{equation}
\label{eq0702newno1}
\mu_{\bff}: \pi_1(\bB^{(0)}, b_0) \longrightarrow \Gamma_g.
\end{equation}

Let $\Delta_i= \{ t_i \in \bC \:|\: |t_i| < \epsilon_0 \}$ be a local 
disk coordinate at $Q_i$ on $\bB$, and 
$f_i=\bff|_{M_i}: M_i=\bff^{-1}(\Delta_i) \longrightarrow \Delta_i$ 
be the degeneration over $\Delta_i$. 
Then we may identify $\mu_{\bff}(\alpha_i)$ 
as the marked topological monodromy $\mu_{f_i}$ of $f_i$
defined in \S6.2.
Therefore, there exists an  element 
$\sigma^{(i)} \in \cC_g$ such that 
$\mu_{\bff}(\alpha_i)=\mu_{f_i} \in \bP_g^{-}(\sigma^{(i)})$.
Let $N_i$ denote the pseudo-period of $\mu_{f_i}$.

Let $\varphi_{\wbB^{(0)}}:\wbB^{(0)} \longrightarrow \bB^{(0)}$ 
be the universal covering of $\bB^{(0)}$.
Then $\pi_1(\bB^{(0)}, b_0)$ acts on $\wbB^{(0)}$, and we have 
$\bB^{(0)} \cong \wbB^{(0)}/\pi_1(\bB^{(0)}, b_0)$. 
On the other hand, 
the cyclic group $\bZ/N_i\bZ$ acts on 
$\wDelta_i=\{ u_i \in \bC \:|\: |u_i| < \epsilon^{1/N_i} \}$ 
by $u_i \mapsto \be(1/N_i)u_i$,  
and the projection 
$\varphi_{\wDelta_i}: \wDelta_i \rightarrow \Delta_i$ is 
defined by $u_i \mapsto t_i=u^{N_i}$.
We define the complex orbifold structure over $\bB$ by 
\begin{equation}
\label{eq0702no1}
\bB^{{\rm orb}}=
\left( \wbB^{(0)}, \pi_1(\bB^{(0)}, b_0) , \varphi_{\wbB^{(0)}}, \bB^{(0)} \right)
\bigcup_{1 \leq i \leq s} 
\left( \wDelta_i, \bZ/N_i\bZ , \varphi_{\wDelta_i}, \Delta_i \right).
\end{equation}
Note that, although $\pi_1(\bB^{(0)}, b_0)$ is an infinite group 
and $\bZ/N_i \bZ$ is a finite group, 
the orbifold structure $\bB^{{\rm orb}}$ is well-defined 
in the sense of \S4.3.

Let $\wbf^{(0)}:\wbM^{(0)}= \bM^{(0)} \times_{\bB^{(0)}} \wbB^{(0)}
\longrightarrow \wbB^{(0)}$ be the pull back of $\bff^{(0)}$ 
over $\wbB^{(0)}$.
Let $\bM \rightarrow \bM^{\sharp}$ be the contraction map of all the non-cores of 
$F_i$ ($1 \leq i \leq s$), 
and ${\bf f}^{\sharp}:\bM^{\sharp} \rightarrow \bB$ be the 
natural holomorphic map.
Let 
$\{ f_i^{\sharp}:M_i^{\sharp} \rightarrow \Delta_i, \: \wf_i:\wM_i \rightarrow \wDelta_i, \: G_i=\bZ/N_i\bZ \}$
be the orbifold structure of $f_i$ in Definition~\ref{definition0602no1} 
and Remark~\ref{correctorb}.
We define the complex orbifold structure over $\bM^{\sharp}$ by
\begin{equation}
\label{eq0702no2}
(\bM^{\sharp})^{{\rm orb}}=
\left( \wbM^{(0)}, \pi_1(\bB^{(0)}, b_0) , \varphi_{\wbM^{(0)}}, \bM^{(0)} \right)
\bigcup_{1 \leq i \leq s} 
\left(  \wM_i, \bZ/N_i\bZ , \varphi_{\wM_i}, M_i^{\sharp} \right)
\end{equation}
where $\varphi_{\wbM^{(0)}}$ and $\varphi_{\wM_i}$ are natural projections.
Then we have the  orbifold fibration\index{orbifold fibration}
 of  
Definition~\ref{deforbfib} by 
\begin{equation}
\label{eq0702no3}
({\bf f}^{\sharp})^{{\rm orb}}: (\bM^{\sharp})^{{\rm orb}} 
\longrightarrow \bB^{{\rm orb}}.
\end{equation}

Now let
\begin{equation}
\label{eq0702no4}
\bJ^{(0)}: \bB^{(0)} \longrightarrow M_g
\end{equation}
be the moduli map. The lifting of $\bJ^{(0)}$ over $\wbB^{(0)}$ is 
defined as follows.
For a point $x \in \bB^{(0)}$, let  $\gamma(b_0,x)$ be an arc on $\bB^{(0)}$ 
starting from $b_0$ and ending at $x$, and $[\gamma(b_0,x)]$ be its 
homopoty class.
Then $\wbB^{(0)}$ consists of all the pairs ($x, [\gamma(b_0,x)]$).
Let
$$
w_{\gamma(b_0,x)}: \Sigma_g=\bff^{-1}(b_0) \longrightarrow \bff^{-1}(x)
$$
be the oriented homeomorphism of fibers of $\bff^{(0)}$ along $\gamma(b_0,x)$.
Since the isotopy class of $w_{\gamma(b_0,x)}$ 
defines a Teichm\"{u}ller marking of $\bff^{-1}(x)$, 
the map 
\begin{equation}
\label{eq0702no5}
\wbJ^{(0)}: \wbB^{(0)} \longrightarrow T_g
\end{equation}
is defined by 
$\wbJ^{(0)} \left( (x, [\gamma(b_0,x)]) \right)
=[\bff^{-1}(x), w_{\gamma(b_0,x)}]$.
Then $\wbJ^{(0)}$ is a holomorphic map which satisfies 
$\varphi_{T_g} \circ \wbJ^{(0)}=\bJ^{(0)}$, 
i.e. it is a lifting of $\bJ^{(0)}$ in (\ref{eq0702no4}).

On the other hand, it follows from (\ref{equation0603no1}) that 
$\bJ^{(0)}$ is holomorphically extended to 
\begin{equation}
\label{eq0702no6}
\bJ: \bB \longrightarrow \overline{M}_g.
\end{equation}
By (\ref{equation0603no3b}), the restricted map 
$J_{f_i}=\bJ|_{\Delta_i}: \Delta_i \longrightarrow \overline{M}_g$ 
is lifted to
\begin{equation}
\label{eq0702no7}
\wJ_{\wf_i}:\wDelta_i \longrightarrow D_{\epsilon}(\sigma^{(i)}).
\end{equation}
We assure the compatibility condition for the patching of 
(\ref{eq0702no5}) and (\ref{eq0702no7}) as an orbifold map.
If $\bJ$ is a constant map, i.e. $\bJ(\bB)=b_0$, 
then the maps (\ref{eq0702no5}) and (\ref{eq0702no7}) are clearly well-patched.
Therefore we assume that $\bJ(\bB)$ is an analytic curve on $\overline{M}_g$ .
Let $x \in \bB^{(0)} \cap (\Delta_i \backslash Q_i)$ be a point.
Let $\wx \in \wbB^{(0)}$ and $x^{\sharp} \in \wDelta_{i}$ be the points 
with $\varphi_{ \wbB^{(0)}}(\wx)=\varphi_{\wDelta_{i}}(x^{\sharp})=x$.
Let $V_x$ be a sufficiently small open neighborhood of $x$ of 
$\bB^{(0)} \cap (\Delta_i \backslash Q_i)$.
Let $\wV_{\wx}$ (resp.~$V^{\sharp}_{x^{\sharp}}$) be the open 
neighborhood of $\wx$ of $\wbB^{(0)}$ (resp.~of $x^{\sharp}$ of $\wDelta_{i}$)  
with $\varphi_{ \wbB^{(0)}}(\wV_{\wx})=\varphi_{\wDelta_{i}}(V^{\sharp}_{x^{\sharp}})=V_x$.
We set 
$\wy=\wbJ^{(0)}(\wx) \in T_g=D_{\epsilon}(\emptyset)$ and 
$y^{\sharp}=\wJ_{\wf_i}(x^{\sharp}) \in D_{\epsilon}(\sigma^{(i)})$.

\begin{lemma}
\label{Lemma0702no1}
There exist open neighborhoods $\wU_{\wy}$ of $\wy$ in $T_g$, and 
$U^{\sharp}_{y^{\sharp}}$ of $y^{\sharp}$ in $D_{\epsilon}(\sigma^{(i)})$ respectively 
such that

\noindent 
{\rm (i)} $\wbJ^{(0)}(\wV_{\wx})$ (resp.~$\wJ_{\wf_i}(V^{\sharp}_{x^{\sharp}})$) 
 is an analytic curve in $\wU_{\wy}$ (resp. $U^{\sharp}_{y^{\sharp}}$), 

\noindent 
{\rm (ii)} there exists a biholomorphic map 
$\psi:\wU_{\wy} \longrightarrow U^{\sharp}_{y^{\sharp}}$ 
so that the restriction map of $\psi$ to $\wbJ^{(0)}(\wV_{\wx})$ induces an 
isomorphism onto $\wJ_{\wf_i}(V^{\sharp}_{x^{\sharp}})$.
\end{lemma}

{\it Proof} \quad 
By Arbarello--Cornalba's theorem \cite{AC2009}, 
we may assume that $\wU_{\wy}$ is the base of a sufficiently small Kuranishi 
family of the Riemann surface $\bff^{-1}(x)$.
On the other hand, $D_{\epsilon}(\sigma^{(i)})$ is a union of the bases of 
 standard Kuranishi families\index{standard Kuranishi family}
  of stable curves by Theorem \ref{KuranishioverCont}.
The Kodaira--Spencer map\index{Kodaira--Spencer map} is an isomorphism at any point, 
particularly at $y^{\sharp}$, of the base of the standard Kuranishi family 
by the property (iii) in \S3.1.
Therefore we also may assume that $U^{\sharp}_{y^{\sharp}}$ 
is the base of a sufficiently small Kuranishi 
family of $\bff^{-1}(x)$.

Hence we may assume that there exists a biholomorphic map 
$\psi':\wU_{\wy} \longrightarrow U^{\sharp}_{y^{\sharp}}$ 
such that $\psi'$ is induced from the extension $\overline{g_0}$ 
to the Kuranishi family of an automorphism 
$g_0 \in {\rm Aut}(\bff^{-1}(x))$ by the property (v) in \S3.1.

Since $\varphi_{T_g}:T_g \rightarrow M_g$ and 
$\varphi_{D_{\epsilon}(\sigma^{(i)})}:D_{\epsilon}(\sigma^{(i)}) 
\rightarrow M_{\epsilon}(\sigma^{(i)})$ are the forgetting maps of 
the Teichm\"{u}ller marking\index{Teichm\"{u}ller marking}
 and the Weyl marking\index{Weyl marking} respectively, 
we have 
$\varphi_{T_g}(\wU_{\wy})
=\varphi_{D_{\epsilon}(\sigma^{(i)})}(U^{\sharp}_{y^{\sharp}})$ 
as an open set (in the classical topology) of $M_g$ containing $\bJ^{(0)}(x)$.

Now we consider the analytic curves 
$\wbJ^{(0)}(\wV_{\wx})$ and $\wJ_{\wf_i}(V^{\sharp}_{x^{\sharp}})$ 
on the above spaces 
$\wU_{\wy}$ and $U^{\sharp}_{y^{\sharp}}$ respectively.
It is clear that 
$\varphi_{T_g}(\wbJ^{(0)}(\wV_{\wx}))=
\varphi_{D_{\epsilon}(\sigma^{(i)})}(\wJ_{\wf_i}(V^{\sharp}_{x^{\sharp}}))
=\bJ^{(0)}(V_x)$.
Hence there exists an element $g_0 \in {\rm Aut}(\bff^{-1}(x))$ 
such that the restriction of the 
biholomorphic map 
$\overline{g_0} \circ \psi':\wU_{\wy} \longrightarrow U^{\sharp}_{y^{\sharp}}$ 
to $\wbJ^{(0)}(\wV_{\wx})$ sends isomorphically onto 
$\wJ_{\wf_i}(V^{\sharp}_{x^{\sharp}})$.
\qed

\bigskip

We extend  the canonically extended moduli map 
$\bJ: \bB \longrightarrow \overline{M}_g$ in 
(\ref{eq0702no6}) to the orbifold map as follows.
Lemma \ref{Lemma0702no1} assures its well-definedness.

\begin{definition}
\label{Definition0702no1}
The maps 
$\wbJ^{(0)}: \wbB^{(0)} \longrightarrow T_g$ in (\ref{eq0702no5})
and  
$\wJ_{\wf_i}:\wDelta_i \longrightarrow D_{\epsilon}(\sigma^{(i)})$ 
in (\ref{eq0702no7}) ($1 \leq i \leq s$) 
are well-patched and define the orbifold map 
$$
\bJ_{\bff}^{{\rm orb}}:\bB^{{\rm orb}} \longrightarrow \overline{M}_g^{orb}. 
$$
We call $\bJ_{\bff}^{{\rm orb}}$ 
{\em the (global) orbifold moduli map}\index{orbifold moduli map}
of a fibered complex surface\index{fibered complex surface} 
${\bf f}:\bM \rightarrow \bB$.
\end{definition}


\subsection{Global recovery of basic members of fibered complex surfaces}
\label{subsec7.3}

We construct all the  fibered complex surfaces 
for  possible monodromy representations and orbifold moduli maps 
by pulling back from our universal degenerating family. 

We let $D(\bB)=\{ Q_1, \cdots, Q_s \}$ 
be a finite set of a compact Riemann surface $\bB$, 
and set
$\bB^{(0)}=\bB \setminus D(\bB)$.
For a fixed point $b_0 \in \bB^{(0)}$, we consider a 
representation of 
the fundamental group to the mapping class group of genus $g \geq 2$
\begin{equation}
\label{eq0703newno1}
\bmu: \pi_1(\bB^{(0)}, b_0) \longrightarrow \Gamma_g.
\end{equation}
We assume that $\bmu$ satisfies the following condition.
Let $\alpha_i$ be a loop  
 which goes around $Q_i$ once counterclockwise 
for $1 \leq i \leq s$.
Then there exists an  element 
$\sigma^{(i)} \in \cC_g$ such that 
$\bmu(\alpha_i) \in \bP_g^{-}(\sigma^{(i)})$.
We call such a $\bmu$ 
{\it a pseudo-periodic representaion}\index{pseudo-periodic representaion} 
on $\bB$.

Let $N_i$ be the pseudo-period of $\bmu(\alpha_i)$ and 
$\bB^{{\rm orb}}$ be the orbifold structure over $\bB$ as in 
(\ref{eq0702no1}).
Let $\bJ: \bB \rightarrow \overline{M}_g$ be a holomorphic map, and 
\begin{equation}
\label{eq0703newno2}
\bJ^{{\rm orb}}: \bB^{{\rm orb}} \longrightarrow \overline{M}_g^{{\rm orb}}
\end{equation}
be an orbifold map over $\bJ$.
We assume that $\bJ^{{\rm orb}}$ satisfies the following conditions:

\noindent {\rm (i)} \ 
the restricted chart map on $\wbB^{(0)}$ maps $\wbB^{(0)}$ to $T_g$, i.e.
$\wbJ: {\wbB^{(0)}} \longrightarrow T_g=D_{\epsilon}(\emptyset)$,  

\noindent {\rm (ii)} \ 
the restricted chart map on each $\wDelta_i$ belongs 
(in Def.~\ref{Def0701no1}) to 
$$
\wbJ|_{\wDelta_i} \in 
\bP^*(\wDelta_i, \Delta_i; 
D_{\epsilon}(\sigma^{(i)}), M_{\epsilon}(\sigma^{(i)}); 
T_{\sigma^{(i)}}^{\mu_{\sigma^{(i)}}}).
$$
We call such $\bJ^{{\rm orb}}$ 
{\it a dual pseudo-periodic map}\index{dual pseudo-periodic map} 
associated with $\bmu$.

\begin{theorem}
\label{Th0703no1}
 {\rm (i)} \ For any 
 $(\bmu, \bJ^{{\rm orb}}) \in {\cal F}(\bB)$, 
 there uniquely exists (modulo analytc equivalence) a fibered 
 complex surface\index{fibered 
 complex surface} of genus $g \geq 2$ 
\begin{equation}
\label{eq0703neqno3}
\bff: \bM \longrightarrow \bB
\end{equation}
such that the monodromy representation of $\bff$ coincides with 
$\bmu$, and the orbifold moduli map of $\bff$ coincides with 
$\bJ^{{\rm orb}}$.
 
\noindent {\rm (ii)} \ 
The  orbifold fibration\index{orbifold fibration}  
$({\bf f}^{\sharp})^{{\rm orb}}: (\bM^{\sharp})^{{\rm orb}} 
\longrightarrow \bB^{{\rm orb}}$ 
associated with  (\ref{eq0703neqno3}) (see (\ref{eq0702no3})) coincides with 
the orbifold pull-back\index{orbifold pull-back}
 (Def.~\ref{definition0403no3}) 
from the universal degenerating family\index{universal degenerating family}
 (Def.~\ref{defunivdeg})\:
$
\overline{\pi}: \overline{Y}_g^{orb} \longrightarrow  \overline{M}_g^{orb}.
$
\end{theorem}

{\it Proof} \quad 
The restricted family of $\overline{\pi}$ to $D_{\epsilon}(\emptyset)=T_g$ 
is nothing but the universal family (\cite{Bers1973}) 
of  Teichm\"{u}ller-marked 
Riemann surfaces 
$Y_g \rightarrow T_g$.
Since $\wbJ({\wbB^{(0)}})$  is contained in $T_g$ by  assumption, 
the family over $\wbB^{(0)}$ obtained by the pull-back via 
$\wbJ$ 
from this universal family  is a holomorphic family 
$\widetilde{\bff}^{(0)}: \wbM^{(0)} \rightarrow \wbB^{(0)}$ of 
Teichm\"{u}ller-marked 
Riemann surfaces.
Since $\bB^{(0)}$ is the quotient space of $\wbB^{(0)}$ by 
$\pi_1(\bB^{(0)}, b_0)$, the forgetting map of the markings induces 
a holomorphic family 
$\bff^{(0)}: \bM^{(0)} \rightarrow \bB^{(0)}$ of 
Riemann surfaces.

On the other hand, it follows from Theorem \ref{theorem0603no2} that 
the pull-back via $\wbJ|_{\wDelta_i}$ of $\overline{\pi}$ 
induces a family of stable curves 
$\widetilde{f}_i :\wM_i \rightarrow \wDelta_i$ and its quotient 
family $f_i^{\sharp}:M_i^{\sharp} \rightarrow \Delta_i$.

By the same argument as that of the well-definedness of 
Definition \ref{Definition0702no1}, 
the maps ($\widetilde{\bff}^{(0)}, \bff^{(0)}$) and 
($\widetilde{f}_i, f_i^{\sharp}$) ($1 \leq i \leq s$) are 
well-patched as orbifold maps and define an 
 orbifold fibration  
$({\bf f}^{\sharp})^{{\rm orb}}: (\bM^{\sharp})^{{\rm orb}} 
\rightarrow \bB^{{\rm orb}}$. 
Hence by the argument in \S6.2, we have a 
fibered complex surface $\bff: \bM \rightarrow \bB$.

It is clear from the construction that 
the monodromy representation of $\bff$ coincides with 
$\bmu$, and the orbifold moduli map of $\bff$ coincides with 
$\bJ^{{\rm orb}}$.
The uniqueness of $\bff$ modulo analytic equivalence is also 
clear since it has the fixed orbifold moduli map $\bJ^{{\rm orb}}$.
\qed

\begin{remark}
\label{RemBasicmem}
The construction of the fibered complex surface 
 (\ref{eq0703neqno3}) 
 of genus $g \geq 2$ 
 is analogous to that of 
the basic members of elliptic surfaces (basic elliptic surfaces, for short) 
due to Kodaira \cite[\S8]{Kodaira}; 
nevertheless they are different in the following points.

{\rm (i)} \ Basic elliptic surfaces\index{basic elliptic surface}
 are assumed to have no multiple fibers 
(i.e fibers of types $_mI_b$, $m \geq 2$ given in \cite[p.565]{Kodaira}), 
while our fibered surfaces   (\ref{eq0703neqno3})  are admitted 
to have any singular fibers including multiple fibers.

{\rm (ii)} \ Let ${\cal F}({\cal J}, G)$ be the set of elliptic surfaces without multiple fibers 
which have J-invariant ${\cal J}$ and homological invariant 
$G$ (see \cite[Def.~8.1]{Kodaira}).
Then the basic elliptic surface  in ${\cal F}({\cal J}, G)$ is 
characterized as the unique member which admits a global section 
(\cite[Th.~10.2]{Kodaira}), 
and other members in  ${\cal F}({\cal J}, G)$ are obtained 
from the basic elliptic surface by  twisting defined in 
\cite[Def.~9.2]{Kodaira} (see \cite[Th.~10.1]{Kodaira}).

On the other hand, a  fibered complex surface of genus $g \geq 2$ 
is uniquely determined by the data $(\bmu, \bJ^{{\rm orb}})$ by Th.\ref{Th0703no1}.

This difference essentially comes from the fact that 
an elliptic curve has translation automorphisms, while a 
Reimann surface of genus $g \geq 2$ has no such infinite automorphisms.
\end{remark}

\vskip 5mm
\noindent{\it Acknowledgments}.  
We thank Prof. Sampei Usui for his advice and discussions on log geometry, 
Prof. Toshiyuki Akita for his advice on representations of automorphisms of Riemann surfaces, 
Prof. Kunio Obitsu for his advice on augmented Teichm\"uller spaces, 
Prof. Kazuhiro Konno for discussions on fibered complex surfaces. 
Finally we thank Prof. Athanase Papadopoulos for 
his  interest in our 
results and for  many comments which improved our presentation  
very much. 
The second named author is partially supported by 
JSPS Grant KAKENHI 17H01091.


\begin{small}

Tadashi Ashikaga, Faculty of Engineering, Tohoku-Gakuin University, Tagajo,\\
Miyagi 985-8537, Japan; e-mail: ashikaga@mail.tohoku-gakuin.ac.jp\par
\vskip 3mm
\noindent
Yukio Matsumoto, Department of Mathematics, Gakushuin University, Mejiro,\\
Toshima-ku, Tokyo 171-8588, Japan; e-mail: yukiomat@math.gakushuin.ac.jp\par

\end{small}

\end{document}